\documentstyle{article}

\newcommand{\N}{{\rm I\kern-.5ex N}}
\newcommand{\Z}{{\sf \vrule height 1.55ex depth-1.2ex width.03em\kern-.11em Z
\kern-.9ex Z\kern-.11em\vrule height 0.3ex depth0ex width.03em}}
\newcommand{\Q}{{\rm\kern.2ex\vrule height1.55ex depth-.05ex width.03em\kern-
.7ex Q}}
\newcommand{\R}{{\rm I\kern-.5ex R}}
\newcommand{\Rvar}{{\rm I\kern-.5ex R}}
\newcommand{\C}{{\rm\kern.3ex\vrule height1.55ex depth-.05ex width.03em\kern-
.7ex C}}
\newcommand{\nabp}{\nab \hspace{-1.05ex}
\rule[.5ex]{.2ex}{.8ex}   \hspace{1.05ex}}
\newcommand{\hoed}{\hat{\rule{0ex}{1.7ex}}\,}
\newcommand{\spat}{\hspace{4ex}}

\newcommand{\flip}{ \chi }
\newcommand{\op}{B(H)}
\newcommand{\cop}{B_0(H)}
\newcommand{\cst}{C$^*$-algebra}
\newcommand{\ar}{A_r}
\newcommand{\arh}{\hat{A}_r}
\newcommand{\ah}{\hat{A}}
\newcommand{\deh}{\hat{\Delta}}
\newcommand{\psih}{\hat{\psi}}
\newcommand{\lah}{\hat{\la}}
\newcommand{\pih}{\hat{\pi}}
\newcommand{\sh}{\hat{S}}
\newcommand{\cu}{{\cal U}}
\newcommand{\cp}{{\cal P}}
\newcommand{\cj}{{\cal J}}
\newcommand{\nab}{\nabla}
\newcommand{\nah}{\nabla^{\frac{1}{2}}}
\newcommand{\hde}{H_\sde}
\newcommand{\lade}{\la_\sde}
\newcommand{\cuh}{\hat{\cu}}
\newcommand{\Jh}{\hat{J}}
\newcommand{\Th}{\hat{T}}
\newcommand{\nabh}{\hat{\nab}}

\newcommand{\od}{\odot}
\newcommand{\ot}{\otimes}

\newcommand{\la}{\Lambda}

\newcommand{\om}{\omega}
\newcommand{\io}{\iota}
\newcommand{\vfi}{\varphi}
\newcommand{\vep}{\varepsilon}

\newcommand{\al}{\alpha}
\newcommand{\be}{\beta}
\newcommand{\ga}{\Gamma}
\newcommand{\sde}{\delta}
\newcommand{\de}{\Delta}
\newcommand{\et}{\eta}
\newcommand{\th}{\theta}

\newcommand{\si}{\sigma}
\newcommand{\Mfi}{{\cal M}_{\vfi}}
\newcommand{\Nfi}{{\cal N}_{\vfi}}
\newcommand{\Mps}{{\cal M}_{\psi}}
\newcommand{\Nps}{{\cal N}_{\psi}}
\newcommand{\lafi}{\la_\vfi}
\newcommand{\laps}{\la_\psi}

\newcommand{\text}[1]{\mbox{#1}}

\newcommand{\qed}{\ \hfill \rule{2mm}{2mm}}
\newenvironment{demo}{\smallskip\medskip\noindent\bf Proof :\ \  
\rm}{\qed\bigskip\par }

\newtheorem{definition}{Definition}[section]
\newtheorem{proposition}[definition]{Proposition}
\newtheorem{lemma}[definition]{Lemma}
\newtheorem{corollary}[definition]{Corollary}

\newtheorem{theorem}[definition]{Theorem}

\begin{document}
\begin{center}
\large\bf
$\protect\boldmath C^*$-algebraic quantum groups arising
from algebraic quantum groups.

\end{center}

\begin{center}
\rm J. Kustermans  \footnote{Research Assistant of the
National Fund for Scientific Research (Belgium)} \& A.
Van Daele

Departement Wiskunde

Katholieke Universiteit Leuven

Celestijnenlaan 200B

B-3001 Heverlee

Belgium
\bigskip

June 1996
\end{center}

\subsection*{Abstract}
We associate to an algebraic quantum group a C$^*$-algebraic quantum group
and show that this \newline C$^*$-algebraic quantum group essentially 
satisfies
an upcoming definition of Masuda, Nakagami \& \newline Woronowicz.

\section*{Introduction.}
In 1979, Woronowicz proposed the use of the C$^*$-algebra language in the 
field of 
quantum groups \cite{Wor4}. Since then, a lot of work has been done in this 
area but there is still
no satisfactory definition of a quantum group in the C$^*$-algebra framework.

However, the following subjects are better understood:
\begin{itemize}
\item compact \& discrete quantum groups and their duality theory 
\item some examples: quantum $SU(2)$, quantum Heisenberg, quantum $E(2)$, 
quantum Lorentz, duals of locally compact groups,...
\end{itemize}

A good definition of a quantum group should satisfy the following 
requirements:
\begin{itemize}
\item It should incorporate all the known examples and well understood parts 
of the theory.
\item There will have to be a good balance between the theory which can be 
extracted from
      the definition and the non complexity of the definition.
\item The definition has to allow a consistent duality theory.
\end{itemize}
The most difficult part of finding a satisfactory axiom scheme for quantum 
groups seems
to be the ability to prove the existence and uniqueness of a left Haar weight
from the proposed definition.

At the moment, Masuda, Nakagami and Woronowicz are working on a quasi-final 
definition of 
a quantum group in the C$^*$-algebra framework (see also \cite{MasNak}). 
It is quasi-final in the sense that there is some hope to find simpler axioms 
which 
imply their current definition. Later, we will give a preview of this 
definition.

\medskip

In \cite{VD6}, A. Van Daele introduced the notion of Multiplier Hopf-
algebras. They are natural 
generalizations of Hopf algebras to the case of non-unital algebras. 
Recently, A. Van Daele has looked into the case where such a Multiplier Hopf 
algebra
possesses a non-zero left invariant functional (Haar functional) and found 
some very 
interesting properties \cite{VD1}. It should be noted 
that everything in this theory is of an algebraic nature.

This category of Multiplier Hopf algebras with a Haar functional behaves very 
well in 
different ways:
\begin{itemize}
\item This category includes the compact \& discrete quantum groups.
\item It is possible to construct the dual within this category.
\item The category is closed under the double construction of Drinfeld (see 
\cite{Drab}).
\end{itemize}
The last and the first property imply that this category contains the Lorentz 
quantum group
\cite{PW}. 
However, this new category will not exhaust all quantum groups. It is not so 
difficult
to find many classical groups which are not Multiplier Hopf algebras. Also, 
quantum $E(2)$
will not fit in this scheme. Nevertheless, we still have a nice class of 
algebraic
quantum groups.

\medskip

The main purpose of this paper is the construction of a C$^*$-algebraic 
quantum group
(in the sense of Masuda, Nakagami \& Woronowicz) out of a Multiplier 
Hopf$\,^*$-algebra
which possesses a positive left invariant functional. In a first section, we 
wil give an 
overview of the results of A. Van Daele about such Multiplier Hopf$\,^*$-
algebras.
In a second section, we introduce the C$^*$-algebra together with the 
comultiplication.
From there on, we gradually prove that this C$^*$-algebra fits almost in the 
scheme
of Masuda, Nakagami \& Woronowicz.

\medskip

First, we introduce some notations and conventions. We will always use the 
minimal tensor
product between C$^*$-algebras and use the symbol $\ot$ for this complete 
tensor product. For any C$^*$-algebra $A$, we denote the Multiplier C$^*$-
algebra by $M(A)$. The flip map between two C$^*$-algebras will be denoted by 
$\flip$. 

For the algebraic tensorproducts of vectorspaces and linear mappings, we use 
the symbol $\od$.
The algebraic dual of a vectorspace $V$ will be denoted by $V'$.

Let $H$ be a Hilbert space. Then $B(H)$ will denote the C$^*$-algebra of 
bounded 
operators on $H$, whereas $B_0(H)$ will denote the C$^*$-algebra of compact 
operators on $H$.
Consider vectors $v,w \in H$, then $\om_{v,w}$ is the element in $B_0(H)^*$
such that $\om_{v,w}(x)=\langle x v, w \rangle$ for all $x \in B_0(H)$.

The domain of an unbounded operator $T$ on $H$ is denoted by $D(T)$.

The domain of an element $\al$ which is affiliated with some \cst\ , will be 
denoted
by ${\cal D}(\al)$ (we will use the same notation for closed mappings in a 
\cst\ which arise
from one-parameter groups).

Whenever we say that an unbounded operator is positive, it is included that 
this operator
is also selfadjoint. The same rules apply to elements affiliated with a C$^*$-
algebra.
Let $\al$ be an element affiliated with a \cst\ $A$, then $\al$ is called 
strictly positive
if $\al$ is positive and has dense range (it is then automatically injective).

A one-parameter group $\si$ on a \cst\ $A$ is called norm-continuous if and 
only if
for every $a \in A$, the mapping $\R \rightarrow A:t \mapsto \si_t(a)$ is 
norm-continuous.

We refer to the appendix for some notations and results about weights.

\medskip

As promised, we give now a preview of the definition of a C$^*$-algebraic 
quantum group
according to Masuda, Nakagami \& Woronowicz:

\medskip

Let $B$ be a C$^*$-algebra and $\de$ a non-degenerate $^*$-homomorphism from 
$B$
into $M(B\ot B)$ such that 
\begin{enumerate}
\item $\de$ is coassociative, i.\ e.\ $(\de \ot \io)\de = (\io \ot \de)\de$.
\item $\de$ satisfies the following density conditions:
      $\de(B)(B \ot 1)$ and $\de(B)(1 \ot B)$ are dense  subsets of $B \ot B$.
\end{enumerate}
Furthermore, we assume the existence of the following objects:
\begin{enumerate}
\item a KMS-weight $\vfi$ on $B$ with modular group $\si$,  
\item a norm continuous one parameter group $\tau$ on $B$,
\item an involutive $^*$-anti-automorphism $R$ on $B$,
\end{enumerate}
which satisfy the following properties:
\begin{enumerate}
\item For every $a \in \Mfi$, we have that $\de(a)$ belongs to 
${\overline{\cal M}}_{\io \ot       \vfi}$ and $(\io \ot \vfi)\de(a) = 
\vfi(a) 1$.
\item Consider $a,b \in \Nfi$. Let $\om \in B^*$ such that
      $\om R \tau_{-\frac{i}{2}}$ is bounded and call $\th$ the unique 
element in $B^*$ which 
      extends $\om R \tau_{-\frac{i}{2}}$. Then
      $$\vfi(\,b^* \, (\om \ot \io)\de(a)\,) = \vfi(\,(\th  \ot 
\io)(\de(b^*)) \,a \,) \ .$$
\item \begin{itemize}
      \item $\vfi$ is invariant under $\tau$.
      \item $\vfi$ commutes with $\vfi R$.
      \end{itemize}
\item \begin{itemize}
      \item For every $t \in \R$, we have that $R \tau_t = \tau_t R$.
      \item We have that $\de \tau_t = (\tau_t \ot \tau_t)\de$ for all $t \in 
\R$. 
      \item $\de R = \flip (R \ot R) \de$
      \end{itemize}
\end{enumerate}
Then, we call $(B,\de,\vfi,\tau,R)$ a C$^*$-algebraic quantum group.

\medskip

We call $\vfi$ the left Haar weight, $R$ the anti-unitary antipode and $\tau$ 
the scaling group
of our quantum group. We put $\kappa = R \tau_{-\frac{i}{2}}$, then $\kappa$ 
plays the role
of the antipode of our quantum group.

This definition seems to be a C$^*$-version of a definition of a quantum 
group in the von
Neumann algebra setting (see \cite{MasNak}), which in turn was a 
generalization of the (too restrictive) 
definition of a Kac-algebra (see \cite{E}). We are not sure that this will be 
the ultimate definition
of a C$^*$-algebraic quantum group proposed by Masuda, Nakagami \& 
Woronowicz, but we 
expect that this one gives a fairly good idea of it. 

A possible drawback of this definition is the complexity of the axioms. 
However, we will show
that the C$^*$-algebraic versions of Van Daele's algebraic objects fit almost 
in this scheme. The only difference lies in the fact that we only
can prove that $\vfi$ is relatively invariant with respect to $\tau$
in stead of invariant. It is not clear at the moment whether this 
definition of Masuda, Nakagami \& Woronowicz should be modified in this
respect.

\section{Algebraic quantum groups.}

In this first section, we will introduce the notion of an algebraic quantum 
group as can be 
found in \cite{VD1}. Moreover, we will give an overview of the properties of 
this algebraic quantum 
group. The proofs of these results can be found in the same paper \cite{VD1}. 
After this section, 
we will construct a C$^*$-algebraic quantum group out of this algebraic one, 
thereby heavily
depending on the material gathered in this section. We will first introduce 
some terminology.

\medskip

We call a $^*$-algebra $A$ non-degenerate if and only if we have for every 
$a\in A$
that: 
$$(\forall b \in A: a b =0) \Rightarrow a=0 \text{\ \ \ \ and \ \ \ \ }
(\forall b \in A: b a =0) \Rightarrow a=0 .$$

For a non-degenerate $^*$-algebra $A$, you can define the multiplier algebra 
$M(A)$,
this is a unital $^*$-algebra in which $A$ sits as a selfadjoint ideal (the 
definition of this
multiplier algebra is the same as in the case of C$^*$-algebras).

If you have two non-degenerate $^*$-algebras $A,B$ and a multiplicative 
linear mapping
$\pi$ from $A$ to $M(B)$, we call $\pi$ non-degenerate if and only if the 
vectorspaces
$\pi(A) B$ and $B \pi(A)$ are equal to $B$. Such a non-degenerate 
multiplicative linear map
has a unique multiplicative linear extension to $M(A)$, this extension will 
be denoted by the 
same symbol as the original mapping. Of course, we have similar definitions 
and results for
antimultiplicative mappings. If we work in an algebraic setting, we will 
always use this form
of non degeneracy as opposed to the non degeneracy of $^*$-homomorphisms 
between C$^*$-algebras!

For a linear functional $\om$ on a non-degenerate $^*$-algebra $A$ and any $a 
\in M(A)$ we 
define the linear functionals $\om a$ and $a \om$ on $A$ such that
$(a \om)(x) = \om(x a)$ and $(\om a)(x) = \om(a x)$ for every $x \in  A$.

You can find some more information about non-degenerate algebras in the 
appendix of \cite{VD6}.

\medskip

Now, let $\om$ be a linear functional on a $^*$-algebra $A$, then: 
\begin{enumerate}
\item $\om$ is called positive if and only if
$\om(a^* a)$ is positive for every $a \in A$.
\item We say that $\om$ is faithful if and only if for every $a \in A$, we 
have that
$$(\forall b \in A: \om(a b)=0) \Rightarrow a=0 \text{\ \ \ \ and \ \ \ \ }
(\forall b \in A: \om(b a) =0) \Rightarrow a=0 .$$
\item If $\om$ is positive, then $\om$ is faithful if and only if for every 
$a \in A$, we have that
$$\om(a^* a)=0 \Rightarrow a=0.$$
\end{enumerate}

Let $\om$ be a positive linear functional on a $^*$-algebra $A$. A GNS-pair  
for  $\om$ is by definition a pair $(K,\Gamma)$, where $K$ is a Hilbert space 
and $\Gamma$ is a linear map
from $A$ into $K$ such that $\Gamma(A)$ is dense in $K$ and $\langle 
\Gamma(a) , \Gamma(b) \rangle = \om(b^* a)$ for every $a,b \in A$. It is 
clear that such a pair exist and that it
is unique up to a unitary transformation.

\medskip

We have now gathered the necessary information to understand the following 
definition

\begin{definition}
Consider a non-degenerate $^*$-algebra $A$ and a non-degenerate $^*$-
homomorphism $\de$
from $A$ into $M(A \od A)$ such that
\begin{enumerate}
\item $(\de \od \io)\de = (\io \od \de)\de$.
\item The linear mappings $T_1$, $T_2$ from $A \od A$ into $M(A \od A)$
      such that 
      $$T_1(a \ot b) = \de(a)(b \ot 1) \text{\ \ \ \ and \ \ \ \ } T_2(a \ot 
b) = \de(a)(1 \ot       b)$$ 
      for all $a,b \in A$, are bijections from $A \od A$ to $A \od A$. 
\end{enumerate}
Then we call $(A,\de)$ a Multiplier Hopf$\,^*$-algebra.
\end{definition}

In \cite{VD6}, A.\ Van Daele proves the existence of a unique non-zero $^*$-
homomorphism $\vep$
from $A$ to $\C$ such that $(\vep \od \io)\de = (\io \od \vep)\de =\io$.
Furthermore, he proves the existence of a unique anti-automorphism $S$ on $A$
such that 
$$m(S \od \io)(\de(a)(1 \ot b)) = \vep(a) b \text{\ \ \ \ and \ \ \ \ }
m(\io \od S)((b \ot 1)\de(a)) = \vep(a) b $$
for every $a,b \in A$  (here, $m$ denotes the multiplication map from $A \od 
A$ to $A$).
As usual, $\vep$ is called the counit and $S$ the antipode of $(A,\de)$.
Moreover, $S(S(a^*)^*) = a$ for all $a \in A$. Also, $\flip(S \od S)\de = \de 
S$.

\medskip

Let $\om$ be a linear functional on $A$. We call $\om$ left invariant (with 
respect to
$(A,\de)$), if and only if $(\io \od \om)(\de(a)(b \ot 1)) = \om(a) b$ for 
every $a,b \in A$.
Right invariance is defined in a similar way.

\begin{definition}
Consider a Multiplier Hopf$\,^*$-algebra $(A,\de)$ such that there exists a
non-zero positive linear functional $\vfi$ on $A$ which is left invariant.  
Then we call
$(A,\de)$ an algebraic quantum group.
\end{definition}

For the rest of this paper, we will fix an algebraic quantum group $(A,\de)$ 
together with
a non-zero left invariant positive linear functional $\vfi$ on it.

An important feature of such an algebraic quantum group is the faithfulness 
and unicity
of left invariant functionals:
\begin{enumerate}
\item Consider a left invariant linear functional $\om$ on $A$, then there 
exists a unique
      element $c \in \C$ such that $\om = c \, \vfi$.
\item Consider a non-zero left invariant linear functional $\om$ on $A$, then 
$\om$ is 
      faithful.
\end{enumerate}
In particular, $\vfi$ is faithful.

\medskip

A first application of this unicity result concerns the antipode: 
Because $\vfi S^2$ is left invariant, there exists a unique complex number 
$\mu$ such that
$\vfi S^2 = \mu \vfi$ (in \cite{VD1}, our $\mu$ is denoted by $\tau$!). It is 
not so difficult to prove
in an algebraic way that $|\mu|=1$. The question remains open if there exists 
an example
of an algebraic quantum group (in our sense) with $\mu \neq 1$.

\medskip

It is clear that $\vfi S$ is a non-zero right invariant linear functional on 
$A$. However, in
general, $\vfi S$ will not be positive. Later, we will use our C$^*$-algebra 
approach to prove 
the existence of a non-zero positive right invariant linear functional on $A$.

Of course, we have similar faithfulness and uniqueness results about right 
invariant linear 
functionals.

\medskip

Another non-trivial property about $\vfi$ is the existence of a unique 
automorphism $\rho$
on $A$ such that $\vfi(a b) = \vfi(b \rho(a))$ for every $a,b \in A$. We call 
this the weak KMS-property of $\vfi$ (In \cite{VD1}, our mapping $\rho$ is 
denoted by $\sigma$!).
This weak KMS-property will be crucial to extend $\vfi$ to a weight on the 
C$^*$-algebra level.
Moreover, we have that $\rho(\rho(a^*)^*) = a$ for every $a \in A$.

As usual there exists a similar object $\rho'$ for the right invariant 
functional $\vfi S$, 
i.\ e.\ $\rho'$ is an automorphism on $A$ such that 
$(\vfi S)(a b) = (\vfi S)(b \rho'(a))$ for every $a,b \in A$.

Using the antipode, we can connect $\rho$ and $\rho'$ via the formula
$S \rho' = \rho S$. Furthermore, we have that $S^2$ commutes with $\rho$ and 
$\rho'$. 
The interplay between $\rho$,$\rho'$ and $\de$ is given by the following 
formulas:
$$\de \rho = (S^2 \od \rho)\de \text{\ \ \ \ and \ \ \ \ } \de \rho' = (\rho' 
\od S^{-2})\de.$$

\medskip

It is also possible to introduce the modular function of our algebraic 
quantum group. This
is an invertible element $\sde$ in $M(A)$ such that $(\vfi \od \io)(\de(a)(1 
\ot b)) = \vfi(a) 
\, \sde  \, b $ for every $a,b \in A$. 

Concerning the right invariant functional, we have that $(\io \od \vfi 
S)(\de(a)(b \ot 
1)) = (\vfi S)(a) \, \sde^{-1} \, b$ for every $a,b \in A$.

This modular function is, like in the classical group case, a one dimensional 
(generally unbounded)
corepresentation of our algebraic quantum group:
$$ \de(\sde) = \sde \od \sde \text{\ \ \ \ \ \ } \vep(\sde)= 1 \text{\ \ \ \ 
\ \ } 
S(\sde) = \sde^{-1} .$$

As in the classical case, we can relate the left invariant functional to our 
right invariant 
functional via the modular function: for every $a \in A$, we have that
$$ \vfi(S(a)) =  \vfi(a \sde) = \mu \, \vfi(\sde a) .$$
Not surprisingly, we have also that $\rho(\sde) = \rho'(\sde) = \mu^{-1} 
\sde$.

Another connection between $\rho$ and $\rho'$ is given by the equality
$\rho'(a) = \sde \rho(a) \sde^{-1}$ for all $a \in A$.

\medskip

We have also a property which says, loosely speaking, that every element of 
$A$ has compact 
support:

Consider $a_1,\ldots\!,a_n \in A$. Then there exists an element $c$ in $A$ 
such that
$c a_i = a_i c = a_i$ for every $i \in \{1,\ldots\!,n\}$.

\medskip

In a last part, we are going to say something about  duality.

We define the subspace $\ah$ of $A'$ as follows:
$$\ah = \{ \vfi a \mid a \in A \} = \{ a \vfi \mid a \in A \}.$$
Like in the theory of Hopf$\,^*$-algebras, we turn $\ah$ into a non-
degenerate 
$^*$-algebra:
\begin{enumerate}
\item For every $\om_1,\om_2 \in \ah$ and $a \in A$, we have that $(\om_1 
\om_2)(a) = 
(\om_1 \od \om_2)(\de(a))$.
\item For every $\om \in \ah$ and $a \in A$, we have that $\om^*(a) = 
\overline{\om(S(a)^*)}$.
\end{enumerate}
We should remark that a little bit of care has to be taken by defining the 
product and 
the $^*$-operation in this way.

Also, a comultiplication $\deh$ can be defined on $\ah$
such that $\deh(\om)(x \ot y) = \om(x y)$ for every $\om \in \ah$ and $x,y 
\in A$.

Again, this has to be made more precise. This can be done by embedding $M(A 
\od A)$ into 
$(A \od A)'$ in the right way but we will not go into this subject. A 
definition of the 
comultiplication $\deh$ without the use of such an embedding can be found in 
definition $4.4$
of \cite{VD1}.

In this way, $(\ah,\deh)$ becomes a Multiplier Hopf$\,^*$-algebra. The counit 
$\hat{\vep}$
and the antipode $\hat{S}$ are such that
\begin{enumerate}
\item For every $\om \in \ah$, we have that $\hat{\vep}(\om) = \om(1)$.
\item For every $\om \in \ah$ and every $a \in A$, we have that 
$\hat{S}(\om)(a) = \om(S(a))$.
\end{enumerate}

For any $a \in A$, we define $\hat{a} = a \vfi \in \ah$. The mapping $A 
\rightarrow \ah :
a \mapsto \hat{a}$ is a bijection, which is in fact nothing else but the 
Fourier transform.

\medskip

Next, we define the linear functional $\psih$ on $\ah$ such that
$\psih(\hat{a}) = \vep(a)$ for every $a \in A$. It is possible to prove that 
$\psih$
is right invariant. 

Also, we have that $\psih(\hat{a}^* \hat{a}) = \vfi(a^* a)$ for every $a \in 
A$. This implies
that $\psih$ is a non-zero positive right invariant linear functional on 
$\ah$.

Later, we will also prove that $A$ possesses a non-zero positive right 
invariant linear 
functional.  In a similar way, this functional will give rise to a non-zero 
positive left
invariant linear functional on $\ah$. This will imply that $(\ah,\deh)$ is 
again an
algebraic quantum group.

\section{The reduced bi-C$^*$-algebra.}

In the rest of this paper we are going to construct gradually a \cst ic 
quantum group in the sense of Masuda, Nakagami \& Woronowicz out of our 
algebraic quantum group $(A,\de)$. In this 
second section, we are first going to construct the \cst\ together with the 
comultiplication.
In fact, we will always be working with the reduced \cst ic quantum group, 
because all
our construnctions are on the GNS-space of our left invariant functional 
$\vfi$ on $A$.
It should be noted however that it is also possible to construct a universal 
\cst ic quantum group out of our algebraic one. This section revolves around 
a multiplicative unitary
like in the case of Kac-algebras and the work of Baaj \& Skandalis 
\cite{B-S}. Therefore, the 
definitions and results appearing in this section will resemble the ones in 
those works.
  
\medskip

For the rest of the paper, we fix a GNS-pair $(H,\la)$ for $\vfi$. Note that 
$\la$ is injective because $\vfi$ is faithfull.

\begin{proposition}   \label{def1.1}
There exists a unique unitary operator $W$ on $H \ot H$ such that
$ W \,(\la \od \la)(\de(b)(a \ot 1)) = \la(a) \ot \la(b) $
for all $a,b \in A$.
\end{proposition}

The proof of this proposition is straightforward.
The left invariance of $\vfi$ implies that
$$ \langle  (\la \od \la)(\de(b)(a \ot 1)) , (\la \od \la)
(\de(d)(c \ot 1))  \rangle  =  \langle  \la(a) \ot \la(b) ,
\la(c) \ot \la(d)  \rangle $$
for all $a,b,c,d \in A$. Remembering that $\de(A)(A \od 1) = A \od A$,
the result follows easily.

\medskip

There is yet another way of expressing $W$:

\begin{proposition}  \label{prop1.1}
For all $a,b \in A$, we have that  
$ W (\la(a) \ot \la(b)) =  
(\la \od \la)(\,(S^{-1} \od \io)(\de(b))(a \ot 1)\,).$
\end{proposition}
\begin{demo}
There exist $p_1,\ldots\! ,p_n \, ,  \,q_1,\ldots\! ,q_n \in A$
such that $\sum_{i=1}^n \de(q_i)(p_i \ot 1)  
= a \ot b$. Using lemma 5.4 of \cite{VD6}, we see that
$$\sum_{i=1}^n S(p_i) \ot q_i = (S(a) \ot 1) \de(b).$$
Applying $S^{-1} \od \io$ to this equation and remembering that $S^{-1}$ is
antimultiplicative, we see that
$$\sum_{i=1}^n p_i \ot q_i = (S^{-1} \od \io)(\de(b))\,(a \ot 1).$$
Therefore,
\begin{eqnarray*}
& & W (\la(a) \ot \la(b))
=  \sum_{i=1}^n W\,(\la \od \la)(\de(q_i)(p_i \ot 1)) \\
& & \spat = \sum_{i=1}^n \la(p_i) \ot \la(q_i)
= (\la \od \la)(\,(S^{-1} \od \io)(\de(b))(a \ot 1)\,).
\end{eqnarray*}
\end{demo}

We would like to represent $A$ by left multiplication as operators on $H$. In 
general,
this would result in unbounded operators. However, in the case of our left 
invariant functional $\vfi$, we get bounded operators. In order to prove 
this, we will need the following short lemma. 

\begin{lemma}   \label{lem1.2}
For every $a,b,c \in A$, we have that
$$ (\io \ot \om_{\la(a),\la(b)} )(W) \, \la(c)
= \la(\,(\io \ot \vfi)(\de(b^*)(1 \ot a)) \, c).$$
\end{lemma}
\begin{demo}
For all $d \in A$, we have that
\begin{eqnarray*}
& &  \langle  (\io \ot \om_{\la(a),\la(b)})(W) \, \la(c) , \la(d)  \rangle
=   \langle  W (\la(c) \ot \la(a)) , \la(d) \ot \la(b) \rangle  \\
& & \spat =   \langle  \la(c) \ot \la(a) , W^* (\la(d) \ot \la(b))  \rangle  
=  \langle  \la(c) \ot \la(a) , (\la \od \la)( \de(b)(d \ot 1) )  \rangle  \\
& & \spat = (\vfi \od \vfi)(\,(d^* \ot 1) \de(b^*) (c \ot a)) 
= \vfi( d^* (\io \od \vfi)(\de(b^*)(1 \ot a) ) \, c ) \\
& & \spat = \langle  \la( \, (\io \ot \vfi)(\de(b^*)(1 \ot a) ) \,c ), \la(d) 
 \rangle .
\end{eqnarray*}
In this chain of equalities, we did not use anything but the definition of 
$W$.
\end{demo}

\begin{lemma}
Let $x \in A$, then there exist a unique bounded operator $T$ on $H$
such that $T \la(c) = \la(x c)$ for all $c \in A$.
\end{lemma}
By the previous lemma, this lemma is true for every element of the form
$x = (\io \ot \vfi)( \de(b^*)(1 \ot a) )$ with $a,b \in A$. Because every
element in $A$ is a sum of such elements, our lemma is true for all $x \in A$.
This allows us to introduce the GNS-representation $\pi$ of $\vfi$:

\begin{definition}
We define the mapping $\pi$ from $A$ into $\op$ such that $\pi(a) \la(b)
= \la(a b)$ for all $a,b \in A$. Then $\pi$ is an injective $^*$-homomorphism
such that $\pi(A) H$ is dense in $H$.
\end{definition}
The last statement of this definition follows easily from the fact that
$A^2=A$. Furthermore, using lemma \ref{lem1.2}, we have for all $a,b \in A$:
\begin{equation}
\pi(\,(\io \ot \vfi)( \de(b^*)(1 \ot a)) \,)
= (\io \ot \om_{\la(a),\la(b)})(W)      \label{eq1.1}
\end{equation}
Now, we are in a position to define our reduced \cst \ $\ar$ associated to 
$(A,\de)$.

\begin{definition}
We define $\ar$ as the closure of $\pi(A)$ in $\op$. So, $\ar$
is a non-degenerate sub-\cst \ of $B(H)$.
\end{definition}

Equation \ref{eq1.1} implies that 
\begin{equation}
\ar = \text{\ \ closure of \ \ }
\{ (\io \ot \om)(W) \mid \om \in \cop^* \}  \text{\ \ in \  }  \op \ , 
\label{eq1.2} 
\end{equation}  
which is something familiar to us.

\medskip

We will define our comultiplication on the $C^*$-algebra level. It will
be denoted by the same symbol as the comultiplication on the $^*$-algebra 
level
but it will always be clear from the context which comultiplication is under
consideration.

\begin{definition}
We define the mapping $\de$ from $\ar$ into $B(H \ot H)$ such that
$\de(x) = W^* (1 \ot x) W$ for all $x \in \ar$. Then $\de$ is an injective
$^*$-homorphism.
\end{definition}

The next lemma guarantees that the comultiplication on $\ar$ is
an extension of the comultiplication on $A$.

\begin{lemma}
For any $a \in A$ and $x \in A \od A$, we have that $(\pi \od \pi)(x) \,
\de(\pi(a)) = (\pi \od \pi)(x \de(a))$ and $\de(\pi(a)) \, (\pi \od \pi)(x)
= (\pi \od \pi)(\de(a) x).$
\end{lemma}
\begin{demo}
Choose $b,c \in A$, then
\begin{eqnarray*}
& & (\pi \od \pi)(x)\, W^* \,(1 \ot \pi(a))\,(\la(b) \ot \la(c))
= (\pi \od \pi)(x)\, W^* \,(\la(b) \ot \la(ac)) \\
& &  \spat =(\pi \od \pi)(x)\, (\la \od \la)(\de(ac)(b \ot 1)\,)
= (\la \od \la)(x \de(ac) (b \ot 1)\,) \\
& & \spat  = (\la \od \la)(x  \de(a) \de(c) (b \ot 1)\,)
= (\pi \od \pi)(x  \de(a))\,(\la \ot \la)(\de(c)(b \ot 1)\,) \\
& & \spat = (\pi \od \pi)(x  \de(a)) \, W^* \, (\la(b) \ot \la(c)).
\end{eqnarray*}
Therefore,
$$ (\pi \od \pi)(x)\,W^*\,(1 \ot \pi(a)) = (\pi \od \pi)(x  \de(a)) \, W^*,$$
so
$$
(\pi \od \pi)(x) \de(\pi(a)) = (\pi \od \pi)(x) W^* (1 \ot \pi(a)) W
= (\pi \od \pi)(x \de(a)).
$$
The other equality can be proven in a similar way.
\end{demo}

Using this lemma, it is easy to infer formulas like:
$$ \de(\pi(a)) \, (1 \ot \pi(b)) = (\pi \od \pi)(\de(a)(1 \ot b)\,)$$
for all $a,b \in A$. Hence, we can conclude the following.

\begin{lemma}
We have that $\de$ is a non-degenerate $^*$-homomorphism from $\ar$ into
$M(\ar \ot \ar)$ such that $\de(\ar)(\ar \ot 1)$ and $\de(\ar)(1 \ot \ar)$
are dense subsets of $\ar \ot \ar$.
\end{lemma}

The coassociativity on the $^*$-algebra level is transferred to the
coassocitivity on the $C^*$-algebra level.

\begin{lemma}
We have that $\de$ is coassociative on the \cst \ level.
\end{lemma}
\begin{demo}
Choose $a \in A$.
For any $b,c \in A$, we have that
\begin{eqnarray*}
& &  (\de \ot \io) \de( \pi(a) ) \, ( \pi(b) \ot 1 \ot \pi(c) ) \\
& & \spat\spat =  (\de \ot \io) (\de(\pi(a)) (1 \ot \pi(c))\,) \, (\pi(b) \ot 
1 \ot 1) \\
& & \spat\spat = (\de \ot \io) (\,(\pi \od \pi)(\de(a)(1 \ot c))\,) \,(\pi(b) 
\ot 1 \ot 1)
\\
& & \spat\spat = (\pi \od \pi \od \pi)( \, (\de \od \io)(\de(a)(1 \ot c)) \, 
(b \ot 1 \ot
1)\,)\\
& & \spat\spat = (\pi \od \pi \od \pi)( \, (\io \od \de)(\de(a)(b \ot 1)) \, 
(1 \ot 1 \ot
c)\,) \\
& & \spat\spat = (\io \ot \de) (\,(\pi \od \pi)(\de(a)(b \ot 1))\,) \, (1 \ot 
1 \ot
\pi(c)) \\
& & \spat\spat = (\io \ot \de) (\de(\pi(a))(\pi(b) \ot 1)\,)\,(1 \ot 1 \ot 
\pi(c)) \\
& & \spat\spat = (\io \ot \de) \de(\pi(a)) \, (\pi(b) \ot 1 \ot \pi(c)).
\end{eqnarray*}
The nondegeneracy of $\ar$ implies that $(\de \ot \io) \de( \pi(a) )
= (\io \ot \de) \de(\pi(a))$. Because $\pi(A)$ is dense in $\ar$, we must have
that $(\de \ot \io) \de = (\io \ot \de) \de$.
\end{demo}

In the next theorem, we gather all these results.

\begin{theorem}
We have that $A_r$ is a non-degenerate sub-\cst\ of $\op$ and $\de$ is a non-
degenerate
injective $^*$-homomorphism from $A_r$ to $M(A_r \ot A_r)$ such that:
\begin{enumerate}
\item $(\de \ot \io)\de = (\io \ot \de)\de$
\item The vectorspaces $\de(A_r)(A_r \ot 1)$ and $\de(A_r)(1 \ot A_r)$ are 
dense subsets
      of $A_r \ot A_r$.
\end{enumerate}
\end{theorem}

In fact, this theorem says that $(A_r,\de)$ satisfies the first part of the 
definition
of Masuda, Nakagami \& Woronowicz.

The coassociativity on the $^*$-algebra level implies also easily the
next proposition.

\begin{proposition} \label{prop1.3}
We have that $W$ is multiplicative: $W_{12} W_{13} W_{23} = W_{23} W_{12}$.
\end{proposition}

\medskip \medskip

In a second part of this section , we will also represent the dual Multiplier 
Hopf$\,^*$-algebra $\ah$ on $H$. This will be done in a similar way as before 
and therefore the proofs will be left out.

Remember from section 1 that we have a non-zero positive right invariant 
linear functional  $\psih$ on $\ah$. Moreover, $\psih(\hat{b}\,^* \hat{a}) = 
\vfi(b^* a)$ for all $a,b
\in A$ \ (a). 

Next, we define the linear map $\lah$ from $\ah$ into $H$ such that
$\lah(\hat{a})= \la(a)$ for every $a \in A$. By using equality (a), we see 
that $(H,\lah)$ is a GNS-pair for $\psih$.

\medskip

In a first lemma, we want to prove an expression for $W$ in terms of our new 
GNS-pair.

In the proof of this lemma, we will need the following formula:
\begin{equation}
(\io \od \vfi)(\,(1 \ot a)\de(b)\,)
= S(\,(\io \od \vfi)(\de(a)(1 \ot b))\,)     \label{eq1.3}
\end{equation}
for all $a,b \in A$. A proof of this result can be found in proposition
3.11 of \cite{VD1}. It is in fact nothing else but an algebraic form of the 
strong left
invariance in the definition of Masuda, Nakagami \& Woronowicz.

\begin{lemma}
For every $\om_1,\om_2 \in \ah$, we have that
$W\,(\lah(\om_1) \ot \lah(\om_2)\,)
= (\lah \od \lah)(\deh(\om_1)(1 \ot \om_2)\,)$
\end{lemma}
\begin{demo}
There exist $a,b \in A$ such that $\om_1=\hat{a}$ and $\om_2=\hat{b}$.

Choose $x,y \in A$, then we have that
\begin{eqnarray*}
& & \hspace{-4.5ex} \left( (S^{-1} \od \io)(\de(b)) \, (a \ot 1) \right)\hoed 
(x \ot y)
=    (\vfi \od \vfi)(\,(x \ot y)\,(S^{-1} \od \io)(\de(b))\, (a \ot 1)\,) \\
& & \hspace{-4.5ex} \spat \vfi(x\, (\io \od \vfi)(\,(1 \ot y)(S^{-1} \od 
\io)(\de(b))\,) \,a) 
=  \vfi(x \, S^{-1}(\,(\io \od \vfi)(\,(1 \ot y)\de(b))\,) \,a ) \ ,
\end{eqnarray*}
so, using equation \ref{eq1.3}, we see that
\begin{eqnarray*}
& & \left( (S^{-1} \od \io)(\de(b))\, (a \ot 1) \right)\hoed (x \ot y)
= \vfi(x\, (\io \od \vfi)(\de(y)(1 \ot b))\,a) \\
&  & = (\vfi \od \vfi)(\,(x \ot 1)\de(y)(a \ot b)\,)
= (\hat{a} \od \hat{b})(\,(x \ot 1)\de(y)\,)  \ ,
\end{eqnarray*}
which, by definition 4.4 of \cite{VD1}, equals
$$((\deh(\om_1)(1 \ot \om_2))(x \ot y).$$
It follows that
$$\left( (S^{-1} \od \io)(\de(b)) \, (a \ot 1) \right)\hoed
= \deh(\om_1)(1 \ot \om_2).$$
Using proposition \ref{prop1.1}, we see that
\begin{eqnarray*}
& & W(\lah(\om_1) \ot \lah(\om_2)) = W (\la(a) \ot \la(b))
= (\la \od \la)(\,(S^{-1} \od \io)(\de(b)) \, (a \ot 1) \,) \\
& & \spat =  (\lah \od \lah)( \,\bigl(\,(S^{-1} \od \io)(\de(b)) \, (a \ot 
1)\,
\bigr)\hoed\ )
= (\lah \od \lah)(\,\deh(\om_1)(1 \ot \om_2)\,).
\end{eqnarray*}
\end{demo}

With this expression, we can do the same things for the dual $\ah$ as
we did for $A$ itself.

\begin{definition}
We define the mapping $\pih$ from $\ah$ into $\op$ such that $\pih(\om)
\lah(\th) = \lah(\om \th)$ for all $\om,\th \in \ah$. Then $\pih$ is
an injective $^*$-homomorphism such that $\pih(\ah) H$ is dense in $H$.
\end{definition}

Moreover, for every $\th,\et \in \ah$, we have that
\begin{equation}
(\om_{\lah(\th),\lah(\et)} \ot \io)(W) = \pih(\,(\psih \od \io)(\,(\et^* \ot
1)\deh(\th))\,).    \label{eq1.4}
\end{equation}

\begin{definition}
We define $\arh$ as the closure of $\pih(\ah)$ in $\op$. So, $\arh$
is a non-degenerate sub-\cst\ of $B(H)$.
\end{definition}

Equation \ref{eq1.4} implies that
$$\arh = \text{\ \ closure of \ \ } \{ (\om \ot \io)(W) \mid \om \in \cop^* \}
\text{\ \ in\ } \op,$$ which is again something familiar.

\begin{definition}
We define the mapping $\deh$ from $\arh$ into $B(H \ot H)$ such that
$\deh(x) = W (x \ot 1) W^*$ for all $x \in \ar$. Then $\deh$ is an injective
$^*$-homorphism.
\end{definition}

\begin{lemma}
For any $\om \in \ah$ and $\th \in \ah \od \ah$, we have that $(\pih \od
\pih)(\th)
\, \deh(\pih(\om)) = (\pih \od \pih)(\th \, \deh(\om))$ and
$\deh(\pih(\om)) \, (\pih \od \pih)(\th) = (\pih \od \pih)(\deh(\om) \, \th 
).$
\end{lemma}

\begin{theorem}
We have that $\ah_r$ is a non-degenerate sub-\cst\ of $\op$ and $\deh$ is a 
non-degenerate
injective $^*$-homomorphism from $\ah_r$ to $M(\ah_r \ot \ah_r)$ such that:
\begin{enumerate}
\item $(\deh \ot \io)\deh = (\io \ot \deh)\deh$
\item The vectorspaces $\deh(\ah_r)(\ah_r \ot 1)$ and $\deh(\ah_r)(1 \ot 
\ah_r)$ are dense                               subsets of $\ah_r \ot \ah_r$.
\end{enumerate}
\end{theorem}

It is also interesting to prove another formula about $\pih$.

\begin{lemma}   \label{prop1.2}
For every $\om \in \ah$, $x \in A$, we have that $\pih(\om) \la(x) =
\la(\,(\sh^{-1}(\om) \od \io) \de(x)\,)$
\end{lemma}
\begin{demo}
Choose $a \in A$. Then
$$(\om \hat{x})(a)= (\om (x\vfi))(a) = (\om \od \vfi)(\de(a)(1 \ot x))
= \om(\,(\io \od \vfi)(\de(a) (1 \ot x))\,)$$
By using equation \ref{eq1.3}, we see that
\begin{eqnarray*}
(\om \hat{x})(a) & = & \omega(S^{-1}(\,(\io \od \vfi)( (1 \ot
a)\de(x))\,)\,) 
= \sh^{-1}(\om) \, (\,(\io \od \vfi)(\,(1 \ot a) \de(x)\,)\,) \\
& = & \vfi( \, (\sh^{-1}(\om) \od \io)((1 \ot a)\de(x))\,) 
= \vfi(a \, (\sh^{-1}(\om) \od \io)\de(x)\,) \\
& = &\bigl( (\sh^{-1}(\om) \od \io) \de(x) \bigr) \hoed (a).
\end{eqnarray*}
So, $\om \, \hat{x} = \bigl( (\sh^{-1}(\om) \od \io) \de(x) \bigr) \hoed $.
Hence,
$$
\pih(\om) \la(x) = \pih(\om) \lah(\hat{x})
= \lah(\om \hat{x})
= \la(\,(\sh^{-1}(\om) \od \io) \de(x)\,)
$$
\end{demo}

\begin{lemma}
For every $a,b \in A$ we have that
$(\om_{\la(a),\la(b)} \ot \io)(W) = \pih(a \vfi b^*)$.
\end{lemma}
\begin{demo}
Choose $y \in A$, then
\begin{eqnarray*}
 \langle  (\om_{\la(a),\la(b)} \ot \io)(W)\,\lah(\hat{x}) , \la(y)  \rangle
&=&  \langle  (\om_{\la(a),\la(b)} \ot \io)(W)\,\la(x) , \la(y)  \rangle  \\
&=&  \langle  W\,(\la(a) \ot \la(x)\,) , \la(b) \ot \la(y)  \rangle  \\
\end{eqnarray*}
Therefore, proposition \ref{prop1.1} implies that 
$$\langle  (\om_{\la(a),\la(b)} \ot \io)(W)\,\lah(\hat{x}) , \la(y)  \rangle
= \langle  (\la \od \la)(\,(S^{-1} \od \io)(\de(x)) \, (a \ot 1) \,) ,
\la(b) \ot \la(y)  \rangle. $$
By using this equality, we see that
\begin{eqnarray*}
& & \hspace{-3ex} \langle  (\om_{\la(a),\la(b)} \ot \io)(W)\,\lah(\hat{x}) , 
\la(y)  \rangle  
= (\vfi \od \vfi)(\,(b^* \ot y^*)(S^{-1} \od \io)(\de(x)) (a \ot 1)\,) \\
& & \hspace{-3ex} \spat = (a \vfi b^* \od \vfi)(\,(S^{-1} \od \io)((1 \ot 
y^*)\de(x))\,)
= ( \sh^{-1}(a \vfi b^*) \od \vfi)(\,(1 \ot y^*)\de(x) \,) \\
& & \hspace{-3ex} \spat = \vfi(\,(\sh^{-1}(a \vfi b^*) \od \io)((1 \ot 
y^*)\de(x))\,)
= \vfi(y^* \,(\sh^{-1}(a \vfi b^*) \od \io)\de(x) \,) \\
& & \hspace{-3ex} \spat =  \langle \la(\,(\sh^{-1}(a \vfi b^*) \od 
\io)\de(x)\,),\la(y) \rangle
=  \langle \lah(\,(a \vfi b^*)\,\hat{x}\,),\la(y) \rangle   \ ,
\end{eqnarray*}
where we used the previous lemma in the last step.
\end{demo}

\medskip

In the last part of this section, we want to prove that $W$ belongs to
$M(A_r \ot \ah_r)$.

Remembering that $\ah$ is a subset of the $A'$, we can identify $A \od \ah$
with a subspace of the vector space of linear mappings from $A$ into $A$
in the natural way.

\begin{lemma}        \label{lem1.1}
Let $\om \in \ah$ and $a \in A$. Define the linear mapping $F$ from $A$ into
$A$ such that $F(x) = (\io \od \om)(\de(x)(a \ot 1)\,)$.
Then $F$ belongs to $A \od \ah$ and $W\,(\pi \od \pih)(a \ot \om)
= (\pi \od \pih)(F)$.
\end{lemma}
\begin{demo}
There exist $b$ in $A$ such that $\om = b \vfi$. Moreover, there exist
$p_1,\ldots\!,p_n \, , \, q_1,\ldots\!,q_n \in A$ such that
$a \ot b = \sum_{i=1}^n \de(p_i)(q_i \ot 1).$

Now, we have for every $x \in A$ that
$$ F(x) = (\io \od \vfi)(\de(x)(a \ot b)\,) 
= \sum_{i=1}^n (\io \od \vfi)(\de(x p_i)(q_i \ot 1)\,)
=  \sum_{i=1}^n q_i \, \vfi(x p_i).
$$
Therefore, $F = \sum_{i=1}^n q_i \ot p_i \vfi$, which clearly belongs to
$A \od \ah$.

Because $a \ot b = \sum_{i=1}^n \de(p_i)(q_i \ot 1)$ and $\flip(S \od
S)\de = \de S$, we have that
$$ S(b) \ot S(a) = \sum_{i=1}^n (1 \ot S(q_i))\de(S(p_i)). \text{\ \ \ \ \ \ 
\ (a)}$$
Choose $c,d \in A$. Using first proposition \ref{prop1.2} and then
proposition \ref{prop1.1}, we see that
\begin{eqnarray*}
& & W\,(\pi \od \pih)(a \ot \om)\,(\la(c) \ot \la(d))
= W\,(\la(ac) \ot \la(\,(\sh^{-1}(\om) \od \io)\de(d)\,) \\
& &\spat =  (\la \od \la)(\,(S^{-1} \od \io)\bigl( \de(\,(\sh^{-1}(\om) \od
\io)\de(d)\,) \bigr) \,\, (ac \ot 1)\,) \\
& & \spat = (\la \od \la)(\,(S^{-1} \od \io)\bigl(\,(S(ac) \ot 1)
\de(\,(\sh^{-1}(\om) \od \io)\de(d)\,)\,\bigr)\,).   \text{\ \ \ \ \ \ \ \ \ 
(b)}
\end{eqnarray*}
Furthermore,
\begin{eqnarray*}
& & (S(ac) \ot 1)\de(\,(\sh^{-1}(\om) \od \io)\de(d)\,)
= (S(ac) \ot 1) \de(\,(\om S^{-1} \od \io)\de(d)\,) \\
& & \spat = (S(ac) \ot 1) \de(\,(\vfi S^{-1} \od \io)((S(b) \ot 1)\de(d))\,)\\
& & \spat = (\vfi S^{-1} \od \io \od \io)\bigl(\,(1 \ot S(ac) \ot 1)
(\io \od \de)(\,(S(b) \ot 1)\de(d)\,)\,\bigr) \\
& & \spat = (\vfi S^{-1} \od \io \od \io)\bigl(\,(S(b) \ot  S(c)S(a) \ot 1)
(\io \od \de)\de(d) \,\bigr) \\
& & \spat = (\vfi S^{-1} \od \io \od \io)\bigl(\,(S(b) \ot S(c)S(a) \ot 1)
(\de \od \io)\de(d) \,\bigr).
\end{eqnarray*}
Refering to (a), we see that
\begin{eqnarray*}
& & (S(ac) \ot 1)\de(\,(\sh^{-1}(\om) \od \io)\de(d)\,) \\
& & \spat = \sum_{i=1}^n (\vfi S^{-1} \od \io \od \io)\bigl(\,(1 \ot 
S(c)S(q_i)
\ot 1) (\de(S(p_i)\,) \ot 1)\,(\de \od \io)\de(d)\,\bigr) \\
& & \spat = \sum_{i=1}^n (\vfi S^{-1} \od \io \od \io)\bigl(\,(1 \ot S(q_i c)
\ot 1) (\de \od \io)((S(p_i) \ot 1)\de(d))\,\bigr) \\
& & \spat = \sum_{i=1}^n S(q_i c) \ot (\vfi S^{-1} \od \io)(\,(S(p_i) \ot
1)\de(d)\,) \text{\ \ \ \ \ \ \ \ \ \ \ \ \ \ \ \ (c)} \\
& & \spat = \sum_{i=1}^n S(q_i c) \ot (\sh^{-1}(p_i \vfi) \od \io)(\de(d)),
\end{eqnarray*}
where in step (c), we used the right invariance of $\vfi S^{-1}$.

Combining this result with equation (b), we get that
\begin{eqnarray*}
& & W\,(\pi \od \pih)(a \ot \om)\,(\la(c) \ot \la(d))
= \sum_{i=1}^n \la(q_i c) \ot \la(\,(\sh^{-1}(p_i \vfi) \od \io)\de(d)\,) \\
& & \spat = \sum_{i=1}^n \pi(q_i) \la(c) \ot  \pih(p_i \vfi) \la(d)
= (\pi \od \pih)(F) (\la(c) \ot \la(d)).
\end{eqnarray*}
Hence, $W\,(\pi \od \pih)(a \ot \om) = (\pi \od \pih)(F)$.
\end{demo}

\begin{proposition}
We have that $W\,(\pi(A) \od \pih(\ah)) = \pi(A) \od \pih(\ah)$.
\end{proposition}
\begin{demo}
From the previous lemma, it follows that $W\,(\pi(A) \od \pih(\ah))
\subseteq \pi(A) \od \pih(\ah)$.

Now, we prove the other inclusion.
Choose $a,b \in A$. Then there exist $p_1,\ldots\!,p_n \, , q_1,\ldots\!,q_n 
\in A$ such that $\sum_{i=1}^n p_i \ot q_i = \de(b)(a \ot 1)$.

For every $x \in A$, we have that
\begin{eqnarray*}
& & \sum_{i=1}^n (\io \od q_i \vfi)(\de(x)(p_i \od 1))
= \sum_{i=1}^n (\io \od \vfi)(\de(x)(p_i \od q_i)\,) \\
& & \spat = \sum_{i=1}^n (\io \od \vfi)(\de(xb)(a \ot 1)\,)
= a \, \vfi(xb) = (a \ot b \vfi)(x).
\end{eqnarray*}
By using \ref{lem1.1}, we see that
$W\,(\pi \od \pih)(\sum_{i=1}^n p_i \ot q_i \vfi) = \pi(a) \ot \pih(b \vfi)$.
\end{demo}

\begin{proposition}
We have that $W$ belongs to $M(A_r \ot \ah_r)$.
\end{proposition}
\begin{demo}
From the previous lemma, we can conclude that $W(A_r \ot \ah_r)
= A_r \ot \ah_r$. If we multiply with $W^*$ to the left,
we get $A_r \ot \ah_r = W^*\,(A_r \ot \ah_r)$.
Taking the adjoint of this equation, gives $A_r \ot \ah_r = (A_r \ot \ah_r)
W$. Hence, $W$ belongs to $M(A_r \ot \ah_r)$.
\end{demo}

\section{The modular group of the left Haar weight.} \label{sec2}

In this section, we want to introduce a norm continuous one-parameter group 
$\si$
on $A_r$ which, in a later stage, will play the role of the modular
group of the left Haar weight. We will define this one-parameter group with
the use of the theory of left Hilbert algebras. However, the one-parameter 
group
which we get in this way has initially two shortcomings:
\begin{itemize}
\item $\si$ leaves only the von Neumann algebra $A_r''$ invariant, where we 
want $\si$ 
      to leave the \cst\ $A_r$ invariant.
\item $\si$ is strongly continuous, we would like to have that $\si$ is norm 
continuous.
\end{itemize}
The main part of this section consists of showing that $\si$ does not have 
this shortcomings.
For this purpose,we will have to depend on the quantum group structure.
In a last part of this section, we prove that $\pi(A)$ consists of elements 
which are
analytic with respect to $\si$ and we show that $\rho$ is the restriction of 
$\si_{-i}$ to $A$.

\medskip

First, we introduce a natural left Hilbert algebra $\cu$.

\begin{definition}
We define $\cu = \la(A)$, so $\cu$ is a dense subspace of $H$.
We will make $\cu$ into a $^*$-algebra in the following way.
\begin{enumerate}
\item For all $a,b \in A$, we have that $\la(a) \la(b) = \la(ab)$.
\item For all $a \in A$, we have that $\la(a)^* = \la(a^*)$.
\end{enumerate}
\end{definition}

\begin{proposition}
We have that $\cu$ is a left Hilbert algebra on $H$.
\end{proposition}
\begin{demo}
\begin{enumerate}
\item Because $A^2=A$, $\cu^2 = \cu$. So, $\cu^2$ is dense in $H$.
\item Choose $a \in A$. It follows immediately that $\la(a)v = \pi(a)v$
for all $v \in \cu$. Hence, the mapping $\cu \rightarrow \cu:v \mapsto
\la(a)v$ is bounded.
\item For every $a \in A \, , \, v_1,v_2 \in \cu$, we have that
\begin{eqnarray*}
 \langle \la(a) v_1,v_2 \rangle  &=&  \langle \pi(a) v_1,v_2 \rangle   =  
\langle v_1,\pi(a^*)v_2 \rangle  \\
&=&  \langle v_1,\la(a^*)v_2 \rangle  =  \langle v_1,\la(a)^*v_2 \rangle .
\end{eqnarray*}
\item For every $a,b \in A$, we have that
\begin{eqnarray*}
 \langle \la(a)^*,\la(b) \rangle
&=&  \langle \la(a^*),\la(b) \rangle  = \vfi(b^* a^*) \\
&=& \vfi(a^* \rho(b^*)) =  \langle \la(\rho(b^*)),\la(a) \rangle .
\end{eqnarray*}
Because $\la(A)$ is dense in $H$, it follows easily that the map
$\cu \rightarrow \cu:v \mapsto v^*$ is closable.
\end{enumerate}
\end{demo}

We will fix some terminology , cf.\ left Hilbert algebras (see e.\ g.\ 
\cite{Stra}).

It is clear that $L_{\la(a)} = \pi(a)$ for all $a \in A$.

We define the mapping $T$ as the closed antilinear map from within H into H
such that $\cu$ is a core for $T$ and $T v = v^*$ for every $v \in \cu$.
From item 4 of the preceding proof, it follows that $\cu$ is a subset of
$D(T^*)$ and that $T^* \la(a) = \la(\rho(a^*)\,)$ for all $a \in A$.
Furthermore, we define $\nab = T^*T$.

Then $\cu$ is a subset of $D(\nab)$ and $\nab \la(a) = \la(\rho(a))$ for all
$a \in A$. So, $\cu$ is also a subset of $D(\nab^{-1})$ and
$\nab^{-1} \la(a) = \la(\rho^{-1}(a))$ for all $a \in A$.

We also define $J$ as the anti-unitary operator arising from the polar
decomposition of $T$. Hence, $T = J \nah$.

\medskip

Later in this section, we will define our one-parameter group $\si$ in the 
same 
way as it is done in the theory of left Hilbert algebras. However, we first 
have to do 
a little bit of work to show that this one-parameter group does not have the 
shortcomings
which were mentioned above. This will be done in the next part.

\begin{lemma}
Let $a \in A$, then
\begin{enumerate}
\item $\vfi(a^* \sde a)$ is positive
\item $\vfi(a^* \sde a) = 0$ if and only if $a=0$.
\end{enumerate}
\end{lemma}
\begin{demo}
\begin{enumerate}
\item Choose $a \in A$.
There exist $b \in A$ such that $\vfi(b^* b) = 1$, then  
$\sde = (\vfi \od \io)\de(b^* b)$.
So,
\begin{eqnarray*}
& & \vfi(a^* \sde a)
= \vfi(a^* (\vfi \od \io)(\de(b^* b)) \,  a)
= (\vfi \od \vfi)(\,(1 \ot a^*)\de(b^* b)(1 \ot a)\,) \\
& & \spat =  \langle  (\la \od \la)(\de(b)(1 \ot a)\,) , (\la \od 
\la)(\de(b)(1 \ot a)\,)  \rangle
\end{eqnarray*}
Hence, $\vfi(a^* \sde a)$ is positive.
\item Assume that $a \in A$ and that $\vfi(a^* \sde a)=0$.
By using the Schwarz equality for positive sesquilinear forms,
it follows for every $b \in A$ that
$$\vfi(b^* \sde a) \leq \vfi(b^* \sde b) \, \vfi(a^* \sde a) \ ,$$
therefore $\vfi(b^* \sde a)=0$.
Because $\vfi$ is faithful, $\sde a= 0$. The invertibility of $\sde$
implies that $a=0$.
\end{enumerate}
\end{demo}

This lemma allows us to define a Hilbert space $\hde$ together
with an injective linear map $\lade$ from $A$ to $\hde$
such that
\begin{enumerate}
\item $\lade$ has dense range in $\hde$
\item $ \langle \lade(a),\lade(b) \rangle
= \vfi(b^* \sde a)$ for all $a,b \in A$.
\end{enumerate}
We will use this new Hilbert space to construct a useful operator.

\begin{lemma}
Let $a,b \in A$, then
$ \langle \lade(S(a^*)),\lade(b) \rangle  =  \langle \la(\sde^{-1} S(b)^* 
\sde) , \la(a) \rangle $.
\end{lemma}
\begin{demo}
We have that
$$  \langle \lade(S(a^*)),\lade(b) \rangle
= \vfi(b^* \sde S(a^*)) = \vfi(S(a^* \sde^{-1} S(b)^*)\,) $$
From section 1, we know that $\vfi S = \sde \, \vfi$,
hence
$$  \langle \lade(S(a^*)),\lade(b) \rangle  = \vfi(a^* \sde^{-1} S(b)^* \sde)
=  \langle \la(\sde^{-1} S(b)^* \sde) , \la(a)  \rangle .$$
\end{demo}

Because $\lade(A)$ is dense in $\hde$, this lemma justifies the following
definition.

\begin{definition}     \label{def2.2}
We define the closed antilinear  mapping  $E$ from within $H$ into $\hde$
such that $\la(A)$ is a core for $E$ and $E \la(a) = \lade(S(a^*))$
for all $a \in A$.

It is clear that $E$ is a densely defined injective operator which has a
dense range.
\end{definition}

The previous lemma, implies easily that
$ \langle Ev,\lade(a) \rangle  =  \langle \la(\sde^{-1} S(a)^* \sde),v 
\rangle $ for all $a \in A$ and $v \in
D(E)$.
In fact, you can use this equality to prove that $E$ is injective.
Also, we can infer from this lemma that $\lade(A)$ is a subset of $D(E^*)$
and that \begin{equation}
E^* \lade(a) = \la(\sde^{-1} S(a)^* \sde)  \label{eq2.1} \end{equation}
for all $a \in A$.

The following operator plays a crucial role in the rest of this section.

\begin{definition}
We define $P = E^* E$. Then $P$ is an injective positive operator in $H$.
\end{definition}

From definition \ref{def2.2} and equation \ref{eq2.1}, it follows that 
$\la(A)$
is a subset of $D(P)$ and that $P \la(a) = \la(\sde^{-1} S^{-2}(a) \sde)$
for every $a \in A$.

\medskip

Now, we want to prove some commutation relations involving $W$, $P$
and $\nab$. The technique that we will use, will be helpful to us
in several situations. First, we prove the following frequently used lemma.

\begin{lemma}  \label{lem2.2}
Consider Hilbert spaces $K,L,H_1,H_2$,\  a unitary operator $U$ from $K$
to $H_2$ and a unitary operator $V$ from $H_1$ to $L$,
$F$ a closed linear operator from within $K$ into $H_1$,
$G$ a closed linear operator from within $H_2$ into $L$.

Suppose there exists a core $C$ for $F$ such that $U(C)$ is a core for
$G$ and such that $V(F(v)) = G(U(v))$  for every $v \in C$.
Then we have that $V F = G U$.
\end{lemma}
\begin{demo}
\begin{itemize}
\item Choose $v \in D(F)$. Then there exists a sequence $(v_n)_{n=1}^\infty$ 
in
$C$ such that $(v_n)_{n=1}^\infty \rightarrow v$
and \newline $(F(v_n)\,)_{n=1}^\infty \rightarrow F(v)$.

For every $n \in \N$, we have that $U(v_n)$ belongs to $D(G)$ and $G(U(v_n))
= V(F(v_n))$. So, $(U(v_n)\,)_{n=1}^\infty$ $\rightarrow U(v)$
and $(G(U(v_n))\,)_{n=1}^\infty \rightarrow V(F(v))$.

Because $G$ is closed, $U(v)$ must belong to the domain of $G$
and $G(U(v)) = V(F(v))$.

\item Choose $v \in K$ such that $U(v)$ belongs to $D(G)$.
Then there exists a sequence $(v_n)_{n=1}^\infty$ in $C$
such that $(U(v_n)\,)_{n=1}^\infty \rightarrow  U(v)$
and $(G(U(v_n))\,)_{n=1}^\infty \rightarrow G(U(v))$.

For every $n \in \N$, we have that $V(F(v_n)) = G(U(v_n))$.
So, $(v_n)_{n=1}^\infty \rightarrow v$ and $(F(v_n)\,)_{n=1}^\infty
\rightarrow V^*(G(U(v)))$. Because $F$ is closed, we must have
that $v$ belongs to the domain of $F$ and $F(v) = V^*(G(U(v)))$,
hence $V(F(v)) = G(U(v)).$
\end{itemize}
\end{demo}

\begin{lemma}
Let $a,b,c,d \in A$, then
\begin{eqnarray*}
& &  \langle (\lade \od \lade)(\flip(\de(b))(a \ot 1)\,) ,
(\lade \od \lade)(\flip(\de(d))(c \ot 1)\,) \rangle   \\
& & \spat \spat =  \langle (\lade(a) \ot \lade(b) , \lade(c) \ot \lade(d) 
\rangle .
\end{eqnarray*}
\end{lemma}
\begin{demo}
We have that
\begin{eqnarray*}
& &  \langle (\lade \od \lade)(\flip(\de(b))(a \ot 1)\,) ,
(\lade \od \lade)(\flip(\de(d))(c \ot 1)\,) \rangle  \\
& & \spat = (\vfi \od \vfi)(\,(c^* \ot 1)\flip(\de(d^*))(\sde \ot 
\sde)(\flip(\de(b))
(a \ot 1)\,) \\
& & \spat = (\vfi \od \vfi)(\,(1 \ot c^*)\de(d^*)(\sde \ot \sde)\de(b)(1 \ot 
a)\,).
\end{eqnarray*}
If we use first that $\de(\sde) = \sde \ot \sde$ and then the fact that
$(\vfi \od \io)\de(x) = \vfi(x) \sde$ for all $x \in A$, we get that
\begin{eqnarray*}
& &  \langle (\lade \od \lade)(\flip(\de(b))(a \ot 1)\,) ,
(\lade \od \lade)(\flip(\de(d))(c \ot 1)\,) \rangle  \\
& & \spat = (\vfi \od \vfi)(\,(1 \ot c^*)\de(d^* \sde b)(1 \ot a)\,)
= \vfi(c^* \sde a) \, \vfi(d^* \sde b) \\
& & \spat =  \langle \lade(a) \ot \lade(b) , \lade(c) \ot \lade(d) \rangle .
\end{eqnarray*}
\end{demo}

Now, we are able to prove a first commutation relation.

\begin{proposition}   \label{prop2.2}
We have that $(P \ot P)\,W = W\,(P \ot P)$.
\end{proposition}
\begin{demo}
From the previous lemma, it follows that there exists a unique unitary 
operator
$V$ on $\hde \ot \hde$ such that
$$V\,(\lade(a) \ot \lade(b)) = (\lade \od \lade)(\flip(\de(b))(a \ot
1)\,).$$
We know that $W^* (\la(A) \od \la(A)) = \la(A) \od \la(A)$.
Hence, $\la(A) \od \la(A)$ and $W^*(\la(A) \od \la(A))$ are cores for
$E \ot E$.

Choose $a,b \in A$. Then
\begin{eqnarray*}
& & (E \ot E)W^* \, (\la(a) \ot \la(b))
= (E \ot E)\,(\la \od \la)(\de(b)(a \ot 1)\,) \\
& & \spat = (\lade \od \lade)(\,(S \od S)(\,\bigl(\de(b)(a \ot 
1)\,\bigr)^*)\,)  \\
& & \spat = (\lade \od \lade)(\,(S \od S)(\de(b^*))(S(a^*) \ot 1) \, ) \\
& & \spat = (\lade \od \lade)(\flip(\de(S(b^*)))(S(a^*) \ot 1) \, ) \\
& & \spat = V \, (\lade \od \lade)( S(a^*) \ot S(b^*) \, ) 
= V(E \ot E) \,(\la(a) \ot \la(b)).
\end{eqnarray*}
Using lemma \ref{lem2.2}, we get that $(E \ot E)W^* = V(E \ot E)$.

By taking the adjoint of this equation, we infer that $W (E^* \ot E^*)
= (E^* \ot E^*) V^*$. By using this two equaltities, we conclude that
\begin{eqnarray*}
 W(P \ot P) &=& W (E^* \ot E^*) (E \ot E)
= (E^* \ot E^*) V^* (E \ot E) \\
 &=& (E^* \ot E^*) (E \ot E) W
= (P \ot P) W.
\end{eqnarray*}
\end{demo}

In order to prove another commutation relation, we will first need
a little result in the algebraic case.

\begin{lemma} \label{lem2.1}
We have that $(\rho^{-1} \od \rho')\de = \de S^{-2}$.
\end{lemma}
\begin{demo}
Choose $a \in A$.
From section 1, we know that
$ \rho'(b) = \sde \rho(b) \sde^{-1}  $
for all $b \in A$.
It follows easily, using the fact that $\rho(\sde) = \frac{1}{\mu} \sde$,
that
$\rho^{-1}(b) = \sde \rho'^{-1}(b) \sde^{-1} $
for all $b \in A$  \ (a).

We use these two equations to conclude that
$$ (\rho^{-1} \od \rho')\de(a) = (\sde \ot \sde)(\rho'^{-1} \od \rho)(\de(a))
(\sde \ot \sde).    \  \ \ (b)$$
From section 1, we know that
$$
(\io \od \rho)\de = (S^{-2} \od \io)\de\rho \text{\ \ \ \ and\ \ \ \ }
(\rho'^{-1} \od \io)\de = (\io \od S^{-2})\de\rho'^{-1}.
$$
We also have that $S^{-2}$ commutes with $\rho$ and $\rho'$.
Using all this information, equation (b) implies
that
$$ (\rho^{-1} \od \rho')\de(a)
=(\sde \ot \sde)(S^{-2} \od S^{-2})(\de(\rho(\rho'^{-1}(a)))\,)(\sde \ot 
\sde).$$
Equation (a) implies that $\rho(\rho'^{-1}(a)) = \sde^{-1} a \sde$.
So, remembering that $\de(\sde) = \sde \ot \sde$ and $S^{-2}(\sde)
= \sde$, we see that
$$ (\rho^{-1} \od \rho')\de(a) = (S^{-2} \od S^{-2})\de(a), $$
which equals $\de(S^{-2}(a)).$
\end{demo}

\begin{lemma}
Consider $a,b,c,d \in A$, then
\begin{eqnarray*}
& &  \langle (\la \od \lade)(\,(\io \od S)((a \ot 1)\de(b))\,) ,
(\la \od \lade)(\,(\io \od S)((c \ot 1)\de(d))\,) \rangle  \\
& & \spat\spat\spat =  \langle \la(a) \ot \la(b) , \la(c) \ot \la(d) \rangle .
\end{eqnarray*}
\end{lemma}
\begin{demo}
We have that
\begin{eqnarray*}
& &  \langle (\la \od \lade)(\,(\io \od S)((a \ot 1)\de(b))\,) ,
(\la \od \lade)(\,(\io \od S)((c \ot 1)\de(d))\,) \rangle  \\
& & \spat = (\vfi \od \vfi)(\,(\io \od S)((c \ot 1)\de(d))^* (1 \ot \sde)
(\io \od S)((a \ot 1)\de(b))\,) \\
& & \spat = (\vfi \od \vfi)(\,(\io \od S^{-1})(\de(d^*)(c^* \ot 1))(1 \ot 
\sde)
(\io \od S)((a \ot 1)\de(b))\,) \\
& & \spat = (\vfi \od \vfi)(\,(\io \od S^{-1})(\de(d^*)) (c^* a \ot \sde)
(\io \od S)(\de(b)) \, ).
\end{eqnarray*}
From this equation, we get by the definition of $\rho$ that
\begin{eqnarray*}
& &  \langle (\la \od \lade)(\,(\io \od S)((a \ot 1)\de(b))\,) ,
(\la \od \lade)(\,(\io \od S)((c \ot 1)\de(d))\,) \rangle   \\
& & \spat = (\vfi \od \vfi)(\,(\rho^{-1} \od \rho^{-1})\bigl( (\io \od 
S)\de(b) \bigr)
(\io \od S^{-1})(\de(d^*)) (c^* a \ot \sde) \,).
\end{eqnarray*}
From section 1, we know that $\rho^{-1} S =S \rho'$. Using this
fact and the previous lemma, we get that
\begin{eqnarray*}
& & \langle (\la \od \lade)(\,(\io \od S)((a \ot 1)\de(b))\,) ,
(\la \od \lade)(\,(\io \od S)((c \ot 1)\de(d))\,) \rangle  \\
& & \spat = (\vfi \od \vfi)(\,(\io \od S)\bigl( (\rho^{-1} \od \rho')\de(b) 
\bigr)
(\io \od S^{-1})(\de(d^*)) (c^* a \ot \sde) \,) \\
& & \spat = (\vfi \od \vfi)(\,(S^{-2} \od S^{-1})(\de(b))
(\io \od S^{-1})(\de(d^*)) (c^* a \sde^{-1} \ot 1)(\sde \ot \sde)\,).
\end{eqnarray*}
So, remembering that $\vfi S  = \sde \, \vfi$, this equation implies
that
\begin{eqnarray*}
& &  \hspace{-4.5ex} \langle (\la \od \lade)(\,(\io \od S)((a \ot 
1)\de(b))\,) ,
(\la \od \lade)(\,(\io \od S)((c \ot 1)\de(d))\,) \rangle  \\
& & \hspace{-4.5ex}\spat = (\vfi \od \vfi)(\,(S(c^* a \sde^{-1}) \ot 1)(S \od 
\io)(\de(d^*))
(S^{-1} \od \io)(\de(b))\,) \\
& & \hspace{-4.5ex}\spat = (\vfi \od \vfi)(\,(\sde S(a) S^{-1}(c)^* \ot 
1)(S^{-1} \od
\io)(\de(d))^* (S^{-1} \od \io)(\de(b)) \,) \\
& & \hspace{-4.5ex}\spat = (\vfi \od \vfi)(\,(S^{-1}(c)^* \ot 1)(S^{-1} \od 
\io)(\de(d))^*
(S^{-1} \od \io)(\de(b)) (\rho(\sde S(a)) \ot 1)\,) \\
& & \hspace{-4.5ex}\spat = (\vfi \od \vfi)([(S^{-1} \od \io)(\de(d)) (S^{-
1}(c) \ot 1)]^*
[(S^{-1} \od \io)(\de(b))(\rho(\sde S(a)) \ot 1)]\,).
\end{eqnarray*}
Because $W$ is unitary and because of proposition \ref{prop1.1}, this last
term equals
$$ (\vfi \od \vfi)(\,(S^{-1}(c) \ot d)^* (\rho(\sde S(a)) \ot b)\,) $$
Hence,
\begin{eqnarray*}
& &  \langle (\la \od \lade)(\,(\io \od S)((a \ot 1)\de(b))\,) ,
(\la \od \lade)(\,(\io \od S)((c \ot 1)\de(d))\,) \rangle  \\
& & \spat = \vfi(S^{-1}(c)^* \rho(\sde S(a))\,) \, \vfi(d^* b)
= \vfi( \sde S(a) S(c^*)) \, \vfi(d^* b) \\
& & \spat = \vfi( S(c^* a \sde^{-1}) \, ) \, \vfi(d^* b).
\end{eqnarray*}
Because $\vfi S = \sde \, \vfi$ , we have that
\begin{eqnarray*}
& &  \langle (\la \od \lade)(\,(\io \od S)((a \ot 1)\de(b))\,) ,
(\la \od \lade)(\,(\io \od S)((c \ot 1)\de(d))\,) \rangle   \\
& & \spat = \vfi(c^* a) \, \vfi(d^* b)
=  \langle \la(a) \ot \la(b) , \la(c) \ot \la(d) \rangle .
\end{eqnarray*}
\end{demo}

We are now able to prove the second commutation relation.

\begin{proposition}
We have that $(\nab \ot \nab) W = W (\nab \ot P)$.
\end{proposition}
\begin{demo}
Because of the previous lemma, there exists a unique unitary operator from
$H \ot H$ to $H \ot \hde$ such that
$$V\,(\la(a) \ot \la(b)) = (\la \od \lade)(\,(\io \od S)((a \ot
1)\de(b))\,).$$
Again, $\la(A) \od \la(A)$ is a core for $T \ot T$ and
$W^*(\la(A) \od \la(A))$ is a core for $T \ot E$.

Choose $a,b \in A$. Then
\begin{eqnarray*}
& & (T \ot E)W^* \, (\la(a) \ot \la(b))  = (T \ot E) \, (\la \od \la)(\de(b)(a
\ot 1)\,) \\
& & \spat = (\la \od \lade)(\,(\io \od S)([\de(b)(a \ot 1)]^*) \,)
= (\la \od \lade)(\,(\io \od S)((a^* \ot 1)\de(b^*)) \,) \\
& & \spat =  V (\la \od \la)(a^* \ot b^*)
= V(T \ot T) \, (\la(a) \ot \la(b)).
\end{eqnarray*}
Using once again lemma \ref{lem2.2}, we see that
$(T \ot E)W^* = V(T \ot T)$.
Just like in proposition \ref{prop2.2}, we deduce that
$(\nab \ot \nab) W = W (\nab \ot P)$.
\end{demo}

How can we use this commutation relation?

Consider $\om \in \cop^*$ and $s \in \R$.
We use the preceding proposition to infer
that $$(\nab^{is} \ot 1) W (\nab^{-is} \ot 1) = (1 \ot \nab^{-is})W(1 \ot
P^{is}).$$ Applying $\io \ot \om$ to this equality gives us
$$\nab^{is} \, (\io \ot \om)(W) \, \nab^{-is} = (\io \ot P^{is} \om
\nab^{-is})(W). $$ Hence,
\begin{enumerate}
\item For all $t \in \R$, we have that $\nab^{it} \, (\io \ot \om)(W) \,
\nab^{-it}$ belongs to $A_r$.
\item The function $\R \rightarrow A_r : t \mapsto \nab^{it}\,(\io \ot 
\om)(W)   \, \nab^{-it}$ is norm continuous.
\end{enumerate}
Using these results and formula \ref{eq1.2} of section 2, we are able
to conclude for any $x \in A_r$,
\begin{enumerate}
\item For all $t \in \R$, we have that $\nab^{it} x \nab^{-it}$
      belongs to $A_r$.
\item The function $\R \rightarrow A_r : t \mapsto \nab^{it} x
      \nab^{-it}$ is norm continuous.
\end{enumerate}

This discussion justifies the following definition.

\begin{definition}   \label{def2.3}
We define the norm-continuous one-parameter group $\si$ on $A_r$
such \newline that $\si_t(x) = \nab^{it} x \nab^{-it}$ for all $t \in \R$ and
$x \in A_r$.
\end{definition}

Remember from the preceding discussion that
$$\si_t(\,(\io \ot \om)(W)\,)
= (\io \ot P^{it} \om \nab^{-it})(W)$$
for every $t \in \R$ and $\om \in \cop^*$.

\medskip

By proposition \ref{prop2.2}, we can justify in a similar way as before the
following definition.

\begin{definition}    \label{def2.1}
We define the norm-continuous one-parameter group $K$ on $A_r$
such \newline that
$K_t(x) = P^{it} x P^{-it}$
for all $t \in \R$ and $x \in A_r$.
\end{definition}

Just as before, we have that
$$ K_t(\,(\io \ot \om)(W)\,) = (\io \ot P^{it} \om P^{-it})(W) $$
for every $t \in \R$ and $\om \in \cop^*$.

\medskip

The proof of the following proposition is now rather straightforward.

\begin{proposition}   \label{prop2.1}
For every $t \in \R$, we have that
$(\si_t \ot K_t)\de = \de \si_t$.
\end{proposition}
\begin{demo}
Choose $x \in A_r$.
Then
\begin{eqnarray*}
\de(\si_t(x)) & = & W^* (1 \ot \si_t(x)) W
= W^* (1 \ot \nab^{it} x \nab^{-it}) W \\
& = & W^* (\nab^{it} \ot \nab^{it})(1 \ot x)(\nab^{-it} \ot \nab^{-it}) W \\
& = & (\nab^{it} \ot P^{it}) W^* (1 \ot x) W (\nab^{-it} \ot P^{-it}) \\
& = & (\nab^{it} \ot P^{it}) \de(x) (\nab^{-it} \ot P^{-it})
=  (\si_t \ot K_t)\de(x).
\end{eqnarray*}
\end{demo}

\medskip\medskip

In a next step, we want to prove that for every element $a \in A$,
we have that $\pi(a)$ belongs to the domain of ${\cal D}(\si_{-i})$ and
$\si_{-i}(\pi(a)) = \pi(\rho(a))$.
Because $\si$ is implemented by $\nab$, we have to prove
that $\pi(a) \nab^{-1} \subseteq \nab^{-1} \pi(\rho(a))$.

We know that
\begin{eqnarray*}
& & \pi(a) \nab^{-1} \la(b) = \pi(a) \la(\rho^{-1}(b))
= \la(a\,\rho^{-1}(b)) \\
& & \spat = \la(\rho^{-1}(\rho(a) b))
= \nab^{-1} \la(\rho(a) b) = \nab^{-1} \pi(\rho(a)) \la(b).
\end{eqnarray*}

This gives an indication that the above claim is true, but we can not
conclude from this that
$\pi(a) \nab^{-1} \subseteq \nab^{-1} \pi(\rho(a))$
because $\la(A)$ does not have to be a core for $\nab$.
Therefore, we have to prove it in another way.
First, we will prove a small result in left Hilbert algebra theory.

\begin{lemma}
Consider $v,w \in \cu''\cap D(\nab)$ such that $\nab(v),\nab(w)$ belong
to $\cu''$. Then $v w$ belongs to $D(\nab)$ and $\nab(v w) = \nab(v) \nab(w)$.
\end{lemma}
\begin{demo}
Because $\nab = T^* T$, \ $v,w$ belong to $D(T)$, \ $T(v),T(w)$ belong to
$D(T^*)$ and $T^*(T(v)) = \nab(v)$ and $T^*(T(w)) = \nab(w)$.

The fact that $T(w)$ belongs to $D(T^*)$ and $T^*(T(w))= \nab(w)$ implies
that $T(w)$, $\nab(w)$ are right closable and $R_{\nab(w)} \subseteq
R_{T(w)}^*$.
Remark also that $\cu''$ is a subset of $D(R_{\nab(w)}) , D(R_{T(w)})$.

Choose $u \in \cu''$. Then
\begin{eqnarray*}
& &  \langle T(u),T(vw) \rangle  =  \langle u^*,(vw)^* \rangle  =  \langle 
u^*,w^* v^* \rangle  \\
& &  \spat = \langle w u^*,v^* \rangle  =  \langle (u w^*)^*,v^* \rangle  =  
\langle T(u w^*),T(v) \rangle  \\
& &  \spat = \langle T^*(T(v)\,),u w^*) \rangle  =  \langle \nab(v),u w^* 
\rangle  =  \langle \nab(v),R_{T(w)}(u) \rangle  \\
& &  \spat = \langle R_{T(w)}^*(\nab(v)\,),u \rangle  =  \langle 
R_{\nab(w)}(\nab(v)\,),u \rangle
=  \langle \nab(v) \nab(w),u \rangle . 
\end{eqnarray*}
Because $\cu''$ is a core for $T$, we have for all $u \in D(T)$
that $ \langle T(u),T(v w) \rangle  =    \langle  \nab(v) \nab(w),u \rangle $.
Hence, $T(vw)$ belongs to $D(T^*)$ and $T^*(T(v w)) = \nab(v) \nab(w)$.
So $vw$ belongs to $D(\nab)$ and $\nab(vw) = \nab(v) \nab(w)$.
\end{demo}

We can even have weaker assumptions:

\begin{lemma}
Consider $v \in \cu''\cap D(\nab)$ such that $\nab(v)$ belongs to $\cu''$
and consider $w \in D(\nab)$. Then $v w$ belongs to $D(\nab)$
and $\nab(v w) = \nab(v) \nab(w)$.
\end{lemma}
\begin{demo}
Put $C=\{ u \in \cu''\cap D(\nab) \mid \nab(u) \in \cu'' \}$,
because the Tomita algebra is a core for $\nab$, $C$ must be a core for 
$\nab$.
Hence, there exists a sequence $(u_n)_{n=1}^\infty$ in $C$ such that
$(u_n)_{n=1}^\infty \rightarrow w$ and $(\nab(u_n)\,)_{n=1}^\infty$.

We know by the previous lemma that $v u_n \in D(\nab)$ and
$\nab(v u_n) = \nab(v) \nab(u_n)$ for all $n \in \N$.
Because $v$ and $\nab(v)$ are left bounded, we get that
$(v u_n)_{n=1}^\infty \rightarrow v u$ and
$(\nab(v u_n)\,)_{n=1}^\infty \rightarrow \nab(v) \nab(u)$.
The closedness of $\nab$ implies that $v u$ belongs to $D(\nab)$ and
$\nab(v u) = \nab(v) \nab(u)$.
\end{demo}

The next lemma is an easy consequence of the previous one.

\begin{lemma}
Consider $v \in \cu''\cap D(\nab^{-1})$ such that $\nab^{-1}(v)$ belongs to
$\cu''$   and consider $w \in D(\nab^{-1})$. Then $v w$ belongs to
$D(\nab^{-1})$ and $\nab^{-1}(v w) =  \nab^{-1}(v) \nab^{-1}(w)$.
\end{lemma}

We formulate these two lemmas in another way which is more useful to us.

\begin{lemma}

\begin{enumerate}
\item Consider $v \in \cu''\cap D(\nab)$ such that $\nab(v) \in \cu''$.
      \newline Then $L_{\nab(v)} \nab \subseteq \nab L_v$.
\item Consider $v \in \cu''\cap D(\nab^{-1})$ such that $\nab^{-1}(v) \in
      \cu''$. \newline Then $L_{\nab^{-1}(v)} \nab^{-1} \subseteq \nab^{-     
 1} L_v$.
\end{enumerate}
\end{lemma}

Now, we are in a position to use these results in our case.

\begin{proposition}
Consider $a \in A$. Then $\pi(a) \, \nab \subseteq \nab \, \pi(\rho^{-1}(a))$
and  $\pi(a) \, \nab^{-1} \subseteq \nab^{-1} \, \pi(\rho(a)\,)$.
\end{proposition}
\begin{demo}
We know that $\la(\rho^{-1}(a))$ belongs to $\cu''\cap D(\nab)$
and $\nab \la(\rho^{-1}(a)) = \la(a)$, which also belongs to $\cu''$.
Hence, by the previous proposition,
$$ \pi(a) \, \nab = L_{\la(a)} \nab = L_{\nab \la(\rho^{-1}(a))} \, \nab
\subseteq \nab \,L_{\la(\rho^{-1}(a))} = \nab \, \pi(\rho^{-1}(a)).$$

The other equality is proven in a similar way.
\end{demo}

\begin{corollary}
Consider $a \in A$ and $n \in \Z$. Then $\pi(a) \nab^n \, \subseteq \nab^n \,
\pi(\rho^{-n}(a))$.
\end{corollary}

Remembering that $\si$ is implemented by $\nab$ and using the previous
corollary, we get the following

\begin{proposition}  \label{prop2.3}
Consider $a \in A$. Then $\pi(a)$ is an analytic element for $\si$
and $\si_{n i}(\pi(a)) = \pi(\rho^{-n}(a))$ for all $n \in \Z$.
\end{proposition}

\section{A construction procedure connected with the polar decomposition of 
the antipode.}

In this section, we wil introduce a fairly general construction procedure 
which in later
sections provides a tool for making a polar decomposition of the antipode. We 
will use this procedure in different cases to get different implementations 
of the scaling group and anti-unitary antipode. The techniques used in this 
section resemble the ones we used for the construction of the modular group 
$\si$.

\medskip

We consider a positive linear functional $\eta$ on $A$ such that there
exist an invertible element $x \in M(A)$ such that
$\eta(a) = \vfi(S(a) x)$ for all $a \in A$.
Furthermore, we take an invertible element $y \in M(A)$.
Let $(K,\ga)$ be a GNS-pair of $\eta$.

\begin{lemma}
Consider $a,b \in A$. Then $ \langle \ga(S(a)^* y),\ga(b) \rangle  =
 \langle \la(S(b^*) \, x \, \rho(S(y))\,),\la(a) \rangle $.
\end{lemma}
\begin{demo}
We have that
\begin{eqnarray*}
& &  \langle \ga(S(a)^* y),\ga(b) \rangle  = \eta(b^* S(a)^* y) = \eta(b^* 
S^{-1}(a^*) y) \\
& & \spat =\eta(S^{-1}(S(y) a^* S(b^*))\,) =
\vfi( S(y) a^* S(b^*) x )\\
& & \spat = \vfi(a^* S(b^*) \, x \, \rho(S(y))\,) 
=  \langle \la(S(b^*) \, x \, \rho(S(y))\,),\la(a) \rangle .
\end{eqnarray*}
\end{demo}

Because, by definition, $\ga(A)$ is dense in $K$, the following definition
is justified.

\begin{definition}
We define $G$ as the closed antilinear map from within $H$ into $K$ such
that $\la(A)$ is a core for $G$ and such that $G \la(a) = \ga(S(a)^* y)$
for every $a \in A$.

It is clear that $G$ is a densely defined injective operator which has dense
range.
\end{definition}

Again, it follows easily that
$$  \langle G v , \ga(b) \rangle  =  \langle \la(S(b^*) \, x \, \rho(S(y))\,) 
, v \rangle  $$
for every $v \in D(G)$.
So $\ga(A)$ is a subset of $D(G^*)$ and $G^* \ga(b) =
\la(S(b^*) \, x \, \rho(S(y))\,)$ for all $b \in A$.

\begin{corollary} \label{prop3.1}
We have that $\ga(A)$ is a core for $G^{-1}$ and
that $G^{-1} \ga(a) =  \la(S(a)^* S(y^{-1})^*)$ for all $a \in A$.
\end{corollary}
\begin{demo}
Because $\la(A)$ is a core for $G$, we have that $G \la(A)$ is a core
for $G^{-1}$. It is clear that $G \la(A) = \ga(A)$, so $\ga(A)$ is a core
for $G^{-1}$. The formula for $G^{-1}$ is easy to prove.
\end{demo}

Now, we want to make a polar decomposition of $G$.
We define $\cj$ as the anti-unitary operator from $H$ to $K$ arising from
the polar decomposition of $G$.
Furthermore, we define $\cp = G^* G$, so $\cp$ is a positive injective 
operator
in $H$ such that $G = \cj \cp^\frac{1}{2}$.

Again, we want to prove some commutation relations. We will
use the same kind of techniques as in the previous section.

\medskip

The left invariance of $\vfi$ allows us to define a unitary operator
on $K \ot H$ such that $V\,(\ga \od \la)(\de(b)(a \ot 1)\,)
= \ga(a) \ot \la(b)$ for all $a,b \in A$.

\begin{proposition}   \label{prop3.2}
We have that $(\cj \ot J)W = V^* (\cj \ot J)$
and $ (\cp \ot \nab)  W =    W (\cp \ot \nab)$.
\end{proposition}
\begin{demo}
For all $a,b \in A$, we have that
\begin{eqnarray*}
& &  V^* (G \ot T) \, (\la(a) \ot \la(b))
= V^* (\ga(S(a)^* y) \ot \la(b^*)) \\
& & \spat = (\ga \od \la)(\de(b^*)(S(a)^* y \ot 1)\,)
= (\ga \od \la)(\,[S(a) \ot 1)\de(b)]^* (y \ot 1) \,) \\
& & \spat = (\ga \od \la)(\,(S \od \io)\bigl(\,(S^{-1} \od \io)(\de(b))\,(a 
\ot 1)
\,\bigr)^* (y \ot 1)\,) \\
& & \spat = (G \ot T)(\la \od \la)(\,(S^{-1} \od \io)(\de(b)) \, (a \ot 1)\,)
= (G \ot T) W \, (\la(a) \ot \la(b)),
\end{eqnarray*}
where, in the last equality, we used proposition \ref{prop1.1}.
Using lemma \ref{lem2.2} once again, we see that
$V^* (G \ot T) = (G \ot T) W$.
As in the proof of proposition \ref{prop2.2}, we get that
$(\cp \ot \nab) W = W (\cp \ot \nab)$

Hence, $W^* (\cp^\frac{1}{2} \ot \nab^\frac{1}{2}) W
= \cp^\frac{1}{2} \ot \nab^\frac{1}{2} $.
Therefore, we see that
\begin{eqnarray*}
& & V\,(\cj \ot J)\,W \,(\cp^\frac{1}{2} \ot \nab^\frac{1}{2})
= V\,(\cj \ot J)\,(\cp^\frac{1}{2} \ot \nab^\frac{1}{2}) \,W \\
& & \spat = V\,(G \ot T)\,W = G \ot T
= (\cj \ot J)\,(\cp^\frac{1}{2} \ot \nab^\frac{1}{2}).
\end{eqnarray*}
Because $\cp^\frac{1}{2} \ot \nab^\frac{1}{2}$ has dense range, it follows
that $V\,(\cj \ot J)\,W = \cj \ot J$.
\end{demo}

The proof of the following lemma is similar to the proof of lemma 
\ref{lem1.2}.

\begin{lemma}   \label{lem3.1}
Consider $a,b,c \in A$.
Then $(\io \ot \om_{\la(a),\la(b)})(V)\,\ga(c) =  
 \ga(\,(\io \od \vfi)(\de(b^*)(1 \ot a))\,c)$.
\end{lemma}

Then, in the same way as in the first section, we get for every $a$ in $A$
the existence of a unique bounded operator $T$ on $K$
such that $T \, \ga(c) = \ga(a c)$ for all $c \in A$.

\begin{definition}
We define the mapping $\th$ from $A$ into $B(K)$
such that $\th(a) \ga(c) = \ga(a c)$ for all $a,c \in A$.
Then $\th$ is a $^*$-homomorphism from $A$ to $B(K)$ such
that $\th(A) K$ is dense in $K$
\end{definition}

We will need the following easy consequence of lemma \ref{lem3.1}.
\begin{equation}
(\io \ot \om_{\la(a),\la(b)})(V^*)
= \th(\,(\io \od \vfi)((1 \ot b^*)\de(a))\,)  \label{eq3.1}  
\end{equation}
for every $a,b \in A$.

\medskip

Now, it is time for some new commutation relations.

\begin{proposition}    \label{prop3.3}
Consider $a \in A$. Then
\begin{enumerate}
\item $\th(a) \, G \subseteq G \,\pi(S(a)^*) $
\item $\pi(a) \, G^* \subseteq G^* \, \th(S(a^*)) $
\item $\pi(a) \, G^{-1} \subseteq G^{-1} \, \th(S(a)^*) $
\item $\th(a) \, (G^{-1})^* \subseteq (G^{-1})^* \pi(S(a^*)) $
\end{enumerate}
\end{proposition}
\begin{demo}
\begin{enumerate}
\item Consider $b \in A$.
      Then we have for all $c \in A$ that
      \begin{eqnarray*}
      & & \th(b) G \, \la(c) = \th(b) \ga(S(c)^* y) = \ga(b S(c)^* y) \\
      & & \spat = \ga(S(S(b)^* c)^* y) = G \la(S(b)^* c) = G \pi(S(b)^*) \, 
\la(c).
      \end{eqnarray*}
      Because $\th(b)$ and $\pi(S(b)^*)$ are bounded and $\la(A)$ is a core
      for $G$, it follows easily that $\th(b)\,G \subseteq G \pi(S(b)^*)$.
\item Consider $b \in A$.
      Apply 1) with $a$ equal to $S^{-1}(b)$, then
      $\th(S^{-1}(b)) \, G \subseteq G \, \pi(b^*)$.
      Taking the adjoint of this equation gives
      $$(G \, \pi(b^*))^*
      \subseteq (\th(S^{-1}(b)) \, G)^* . \ \ \ \text{(a)}  $$

      We know that
      $$ (G \, \pi(b^*)\,)^* \supseteq \pi(b^*)^* \, G^*
      = \pi(b) G^* .$$
      The boundedness of $\th(S^{-1}(b))$ implies
      that $$ (\th(S^{-1}(b)) \, G)^* = G^* \, \th(S^{-1}(b))^*
      = G^* \, \th(S(b^*)) .$$
      Combining this two relations with relation (a) gives us
      $\pi(b) \, G^* \subseteq G^* \th(S(b^*))$.
\item This is proven in a similar way as 1), using corollary \ref{prop3.1}.
\item This is proven in a similar way as 2), using 3).
\end{enumerate}
\end{demo}

\begin{proposition}
Consider $a \in A$. Then
\begin{enumerate}
\item $\pi(a) \, \cp \subseteq \cp \, \pi(S^{-2}(a))$
\item $\pi(a) \, \cp^{-1} \subseteq \cp^{-1} \pi(S^2(a)).$
\end{enumerate}
\end{proposition}
\begin{demo}
\begin{enumerate}
\item Using the previous proposition, we see that
      $$
      \pi(a)\,\cp = \pi(a) \, G^* G \subseteq G^* \, \th(S(a^*)) \, G \\
      \subseteq G^* G \, \pi(S(S(a^*))^*) = \cp \pi(S^{-2}(a)).
      $$
\item This is proven in a similar way, using the fact that $\cp^{-1}
      = (G^{-1})(G^{-1})^*$.
\end{enumerate}
\end{demo}

\begin{corollary}
Consider $a \in A$ and $n \in \Z$. Then $\pi(a) \cp^n \subseteq \cp^n
\pi(S^{-2n}(a))$.
\end{corollary}

\medskip

The following results must be seen in the light of the manageability of our 
multiplicative unitary $W$, as defined by Woronowicz in \cite{Wor3}.

\begin{proposition}  \label{prop3.4}
Consider $u_1 \in D(\cp^\frac{1}{2})$, $u_2 \in D(\cp^{-\frac{1}{2}})$ and
$v_1,v_2 \in H$. Then
$$  \langle W^* (\cp^\frac{1}{2} u_1 \ot J v_1) , \cp^{-\frac{1}{2}} u_2 \ot 
J v_2 \rangle
=  \langle W (u_1 \ot v_2), u_2 \ot v_1 \rangle .$$
\end{proposition}
\begin{demo}
Choose $a_1,a_2 \in A$.

Because $G = \cj \cp^\frac{1}{2}$ and $(G^{-1})^* = \cj \, \cp^{-
\frac{1}{2}}$,
we have that $u_1 \in D(G)$, $u_2 \in D((G^{-1})^*)$
and $\cp^\frac{1}{2} u_1 = \cj^*(G u_1)$, $\cp^{-\frac{1}{2}} u_2
= \cj^*((G^{-1})^* u_2)$.

Using equation \ref{eq3.1} and proposition \ref{prop3.3}, we see that
\begin{eqnarray*}
& & (\io \ot \om_{\la(a_2),\la(a_1)})(V^*) \, (G^{-1})^*
= \th(\,(\io \od \vfi)((1 \ot a_1^*)\de(a_2))\,) \, (G^{-1})^* \\
& & \spat  \subseteq (G^{-1})^* \, \pi( S(\,(\io \od \vfi)((1 \ot
a_1^*)\de(a_2))^*\,)\,) \\
& & \spat = (G^{-1})^* \, \pi( S(\,(\io \od \vfi)(\de(a_2^*)(1 \ot 
a_1))\,)\,) \\
& & \spat = (G^{-1})^* \, \pi(\,(\io \od \vfi)((1 \ot a_2^*)\de(a_1))\,) .
\end{eqnarray*}
where, in the last equality, we used equation \ref{eq1.3} of section 2.

So, using equation \ref{eq1.1}, we get that
$$ (\io \ot \om_{\la(a_2),\la(a_1)})(V^*) \, (G^{-1})^*
\subseteq (G^{-1})^* \, (\io \ot \om_{\la(a_1),\la(a_2)})(W^*) .$$
Hence, $(\io \ot \om_{\la(a_1),\la(a_2)})(W^*) \, u_2$ belongs
to $D(\,(G^{-1})^*)$ and
$$ (G^{-1})^*(\,(\io \ot \om_{\la(a_1),\la(a_2)})(W^*) \, u_2 )
= (\io \ot \om_{\la(a_2),\la(a_1)})(V^*) \, (G^{-1})^* u_2  .
\ \ \ \text{(a)} $$
Furthermore,
\begin{eqnarray*}
& &  \langle W^*  (\cp^\frac{1}{2} u_1 \ot J \la(a_1)) ,
\cp^{-\frac{1}{2}} u_2 \ot J \la(a_2)  \rangle  \\
& & \spat =  \langle W^*  (\cj^* G u_1 \ot J \la(a_1)) , \cj^* (G^{-1})^* u_2 
\ot J
\la(a_2)  \rangle   \\
& & \spat  =  \langle W^* (\cj^* \ot J)(G u_1 \ot \la(a_1)) , (\cj^* \ot J)
((G^{-1})^* u_2 \ot \la(a_2))  \rangle  \\
& & \spat =  \langle (\cj^* \ot J) V  (G u_1 \ot \la(a_1)) , (\cj^* \ot J)
((G^{-1})^* u_2 \ot \la(a_2)) \rangle   \text{\ \ \ \ \ \ (b)} \\
& & \spat =  \langle  (G^{-1})^* u_2 \ot \la(a_2) , V  (G u_1 \ot \la(a_1))  
\rangle  \\
& & \spat =  \langle V^* ((G^{-1})^* u_2 \ot \la(a_2)) , G u_1 \ot \la(a_1)  
\rangle  \\
& & \spat =  \langle (\io \ot \om_{\la(a_2),\la(a_1)})(V^*)\,(G^{-1})^* u_2 , 
G u_1  \rangle .
\end{eqnarray*}
where, in step (b), we used proposition \ref{prop3.2}.

Combining this result with equation (a), we get that
\begin{eqnarray*}
& &  \langle W^* (\cp^\frac{1}{2} u_1 \ot J \la(a_1)) ,
\cp^{-\frac{1}{2}} u_2 \ot J \la(a_2)  \rangle  \\
& & \spat =  \langle  (G^{-1})^*(\,(\io \ot \om_{\la(a_1),\la(a_2)})(W^*) \, 
u_2 ) ,
G u_1 \rangle  \\
& & \spat =  \langle  u_1 , (\io \ot \om_{\la(a_1),\la(a_2)})(W^*)\, u_2  
\rangle  \\
& & \spat =  \langle  (\io \ot \om_{\la(a_2),\la(a_1)})(W) \, u_1, u_2  
\rangle  \\
& & \spat =  \langle  W (u_1 \ot \la(a_2)) , u_2 \ot \la(a_1)  \rangle  .
\end{eqnarray*}
Now, the result follows easily.
\end{demo}

\begin{lemma}
Suppose that $S(y^*) \, x \, \rho(S(y)) \in \C 1$. Then we have that
\begin{eqnarray*}
& &  \langle (\ga \od \ga)( \flip(\de(b)(1 \ot a))(y \ot y) \,) ,
(\ga \od \ga)( \flip(\de(d)(1 \ot c))(y \ot y) \,)  \rangle  \\
& & \spat =  \langle (\ga \od \ga)(\,(a \ot b)(y \ot y)\,) ,
(\ga \od \ga)(\,(c \ot d)(y \ot y)\,) \rangle .
\end{eqnarray*}
\end{lemma}
\begin{demo}
By supposition, there exists a complex number $r$ such that 
$S(y^*) \, x \, \rho(S(y)) = r 1$. We have that
\begin{eqnarray*}
& &  \langle (\ga \od \ga)( \flip(\de(b)(1 \ot a))(y \ot y) \,) ,
(\ga \od \ga)( \flip(\de(d)(1 \ot c))(y \ot y) \,)  \rangle  \\
& & \spat = (\eta \od \eta)([\flip(\de(d)(1 \ot c))(y \ot y)]^*
[\flip(\de(b)(1 \ot a))(y \ot y)]) \\
& & \spat = (\eta \od \eta)(\,(y^* \ot y^*) \flip((1 \ot c^*)\de(d^*))
\flip(\de(b)(1 \ot a))(y \ot y)\,) \\
& & \spat = (\eta \od \eta)(\,(y^* \ot (c y)^*)\de(d^* b)(y \ot ay)\,).
\end{eqnarray*}
Using our connection between $\vfi$ and $\eta$, we get that
\begin{eqnarray*}
& &  \hspace{-4ex} \langle (\ga \od \ga)( \flip(\de(b)(1 \ot a))(y \ot y) \,) 
,
(\ga \od \ga)( \flip(\de(d)(1 \ot c))(y \ot y) \,)  \rangle  \\
& & \hspace{-4ex}\spat = (\vfi \od \eta)(\,(S \od \io)\bigl((y^* \ot (c 
y)^*)\de(d^* b)
(y \ot ay) \bigr) (x \ot 1)\,) \\
& & \hspace{-4ex}\spat = (\vfi \od \eta)(\,(S(y) \ot 1)\,(S \od \io)\bigl((1 
\ot (cy)^*)
\de(d^*b) (1 \ot ay)\bigr) \, (S(y^*) \ot 1)(x \ot 1)\,) \\
& & \hspace{-4ex}\spat = (\vfi \od \eta)(\,(S \od \io)\bigl((1 \ot (cy)^*) 
\de(d^*b)
(1 \ot ay)\bigr) \, (S(y^*) \, x \, \rho(S(y)) \ot 1)\,) \\
& & \hspace{-4ex}\spat = r \, (\vfi \od \eta)(\,(S \od \io)\bigl((1 \ot 
(cy)^*) \de(d^*b)
(1 \ot ay)\bigr) \,) \\
& & \hspace{-4ex}\spat = r \, (\, \vfi S \od \eta)(\,(1 \ot (cy)^*) 
\de(d^*b)(1 \ot ay)\,).
\end{eqnarray*}
Hence, the right invariance of $\vfi S$ implies that
\begin{eqnarray*}
& &  \langle (\ga \od \ga)( \flip(\de(b)(1 \ot a))(y \ot y) \,) ,
(\ga \od \ga)( \flip(\de(d)(1 \ot c))(y \ot y) \,)  \rangle  \\
& & \spat = r \, \vfi(S(d^* b)) \, \eta((cy)^* (ay))
= \vfi(S(d^* b) S(y^*) \, x \, \rho(S(y))\,) \, \eta((cy)^* (ay)) \\
& & \spat = \vfi(S(y) S(d^* b) S(y^*) x ) \, \eta((cy)^* (ay))
= \vfi(S(y^* d^* b y) x ) \, \eta((cy)^* (ay))  \\
& & \spat = \eta((d y)^* (b y)) \, \eta((cy)^* (ay)) \\
& & \spat =  \langle (\ga \od \ga)(\,(a \ot b)(y \ot y)\,) ,
(\ga \od \ga)(\,(c \ot d)(y \ot y)\,) \rangle .
\end{eqnarray*}
\end{demo}

\begin{proposition} \label{prop3.5}
If $S(y^*) \, x \, \rho(S(y)\,) \in \C 1$, then we have
that $W(\cp \ot \cp) =  (\cp \ot \cp) W$.
\end{proposition}
\begin{demo}
The previous lemma implies the existence of a unitary operator $U$ on $K \ot 
K$.
such that
$$ U\,(\ga \od \ga)(\flip(\de(b)(1 \ot a))(y \ot y)\,)
= (\ga \od \ga)(\,(a \ot b)(y \ot y)\,) $$
for all $a,b \in A$.

Choose $c,d \in A$, then
\begin{eqnarray*}
& & (G \ot G)W \, (\la \od \la)(\de(d)(c \ot 1)\,)
= (G \ot G)\,(\la \od \la)(c \ot d)  \\
& & \spat = (\ga \od \ga)(\,(S(c)^* \ot S(d)^*)(y \ot y)\,) \\
& & \spat = U (\ga \od \ga)(\flip\bigl(\de(S(d)^*)(1 \ot S(c)^*)\bigr)(y \ot 
y)\,) \\
& & \spat = U (\ga \od \ga)([(S(c) \ot 1) \flip(\de(S(d)))]^* (y \ot y) \,) \\
& & \spat = U (\ga \od \ga)([(S(c) \ot 1)(S \od S)(\de(d))]^* (y \ot y) \,) \\
& & \spat = U (\ga \od \ga)(\,(S \od S)(\de(d)(c \ot 1))^* (y \ot y) \,) \\
& & \spat = U(G \ot G)\,(\la \od \la)(\de(d)(c \ot 1)\,).
\end{eqnarray*}
Refering to lemma \ref{lem2.2} once again, we see
that $(G \ot G)W = U(G \ot G)$.
As before, we can conclude that $(\cp \ot \cp) W = W (\cp \ot \cp)$.
\end{demo}

In a later stage, we will apply our construction procedure in a case where 
the assumption of
this proposition is satisfied. Then we will have by definition that $W$ is 
amenable.

\section{The polar decomposition of the antipode.}

In this section, we are going to make a polar decomposition of the antipode
$S$.  We will use the construction procedure of section 4 for this.
In this case,  we take $\eta=\vfi$ and $y=1$.
For all $a \in A$, we have that $\frac{1}{\mu} \vfi(S(a) \sde) = \vfi(a)$.
Hence, $x$ is equal to $\frac{1}{\mu} \, \sde$.
Furthermore, we take in this case $K=H$ and $\ga=\la$. So, $\th=\pi$.
We put $I = \cj$ and $M = \cp$.
First, we give a summary of the results of section 4 in this particular case.

\medskip

We have that $G$ is the closed antilinear map from within $H$ into $H$
such that $\la(A)$ is a core for $G$ and $G \la(a) = \la(S(a)^*)$
for every $a \in A$. It is clear that $G$ is involutive in this case.

Therefore,  $I$ is an involutive anti-unitary transformation on $H$,  \
$M=G^* G$ is an injective positive operator in $H$
such that $G = I M^\frac{1}{2} = M^{-\frac{1}{2}} I$.
For $t \in \R$ we have that $M^{it} =  I M^{it} I$ and
$M^t = I M^{-t} I$.

We have the following commutations:
$$ \pi(a)\,G \subseteq G\,\pi(S(a)^*)   \text{\ \ \ \ and  \ \ \ \ }
\pi(a)\,G^* \subseteq G^*\,\pi(S(a^*)) $$
for all $a \in A$.
Moreover, we have for all $a \in A$ and $n \in \Z$ that
\begin{equation}
\pi(a)\,M^n \subseteq M^n\, \pi(S^{-2n}(a)) . \label{eq4.1}  
\end{equation}

In this case, $V$ is equal to $W$, so
$$ (M \ot \nab)W = W(M \ot \nab) \text{\ \ \ \ and \ \ \ \ }
(I \ot J)W = W^* (I \ot J)  . $$

Using these two commutation relations, the following definitions are justified
in the same way as definition \ref{def2.3}.

\begin{definition}
We define the norm-continuous one-parameter group $\tau$ on $A_r$ such
\newline that $\tau_t(x) = M^{it} x M^{-it}$ for every $x \in A_r$ and $t \in 
\R$.
\end{definition}

Again, it follows that
\begin{equation}
\tau_t(\,(\io \ot \om)(W)\,) = (\io \ot \nab^{it} \om \nab^{-it})(W)  
\label{eq4.2}
\end{equation}
for all $t \in \R$ and $\om \in \cop^*$.

\begin{definition}
We define the involutive $^*$-anti automorphism $R$ on $A_r$
such that $R(x) = I x^* I$ for every $x \in A_r$.
\end{definition}

Similarly, we have that
$$R(\,(\io \ot \om)(W)\,) = (\io \ot \om(J\,.^*J)\,)(W) $$
for every $\om \in \cop^*$.

\begin{proposition}
For every $t \in \R$, we have that $\tau_t R = R \tau_t$.
\end{proposition}
\begin{demo}
Choose $x \in A_r$. We know that $M^{it}$ and $I$ commute.
Hence,
\begin{eqnarray*}
& & \tau_t(R(x)) = M^{it} I x^* I M^{-it} = I M^{it} x^* M^{-it} I \\
& & \spat = I (M^{it} x M^{-it})^* I = R(\tau_t(x)).
\end{eqnarray*}
\end{demo}

\begin{corollary}
For every $z \in \C$, we have that $\tau_z R = R \tau_z$.
\end{corollary}

Using relation \ref{eq4.1} and remembering that $\tau$ is implemented by
$M$, we get the following

\begin{proposition}
Consider $a \in A$. Then $\pi(a)$ is an analytic element of $\tau$
and $\tau_{n i}(\pi(a)) = \pi(S^{-2n}(a))$ for every $n \in \Z$.
\end{proposition}

Finally,  we are able to describe the polar decomposition of the antipode.

\begin{theorem}
Consider $a \in A$, then $\pi(S(a)) = R(\tau_{-\frac{i}{2}}(\pi(a))\,)$.
\end{theorem}
\begin{demo}
We know from the beginning of this section that
$\pi(a)\,G \subseteq G\, \pi(S(a)^*)$,
so $\pi(a)\,M^{-\frac{1}{2}} I \subseteq M^{-\frac{1}{2}} I \, \pi(S(a)^*)$.
Therefore,
$$ \pi(a)\, M^{-\frac{1}{2}}  \subseteq  M^{-\frac{1}{2}} I \pi(S(a))^* I
= M^{-\frac{1}{2}} R(\pi(S(a))\,) .$$
Remembering that $\tau$ is implemented by $M$, we get that
$ \tau_{-\frac{i}{2}}(\pi(a)) = R(\pi(S(a))\,)$.
The result follows because $R$  is involutive.
\end{demo}

The rest of this section will devoted to proving some interesting
commutation relations between $\tau, R$ and $\de$. They are always 
consequences of commutation
relations between $W$,$M$,$\nab$,$I$ and $J$.

\begin{proposition}  \label{prop4.1}
For every $t \in \R$, we have that $\de \si_t = (\tau_t \ot \si_t)\de$.
\end{proposition}
\begin{demo}
Choose $x \in A_r$. Then
\begin{eqnarray*}
& &\de(\si_t(a)) = W^*(1 \ot \si_t(x))W
= W^* (1 \ot \nab^{it} x \nab^{-it}) W \\
& & \spat = W^* (M^{it} \ot \nab^{it})(1 \ot x)(M^{-it} \ot \nab^{-it}) W .
\end{eqnarray*}
We know that $W$ commutes with $M \ot \nab$. Hence,
\begin{eqnarray*}
& &\de(\si_t(a))
= (M^{it} \ot \nab^{it}) W^* (1 \ot x) W (M^{-it} \ot \nab^{-it}) \\
& & \spat = (M^{it} \ot \nab^{it}) \de(x) (M^{-it} \ot \nab^{-it})
= (\tau_t \ot \si_t)\de(x).
\end{eqnarray*}
\end{demo}

Using this result, it is not difficult to prove the following one.

\begin{proposition}
For every $t \in \R$, we have that $\de \tau_t = (\tau_t \ot \tau_t) \de$.
\end{proposition}
\begin{demo}
We have that
\begin{eqnarray*}
& & (\,(\tau_t \ot \tau_t)\de \ot \si_t)\de
= (\tau_t \ot \tau_t \ot \si_t)(\de \ot \io)\de \\
& & \spat = (\tau_t \ot \tau_t \ot \si_t)(\io \ot \de)\de
= (\tau_t \ot \de \si_t) \de \\
& & \spat = (\io \ot \de)\de \si_t
= (\de \ot \io)\de \si_t 
= (\de \tau_t \ot \si_t)\de.
\end{eqnarray*}
Therefore,
$$ (\,(\tau_t \ot \tau_t)\de \ot \io)\de
= (\de \tau_t \ot \io)\de  . $$
From this and the fact that $A_r = [\,(\io \ot \om)\de(x) \mid \om \in A_r^*,x
\in A_r]$, it follows easily that $(\tau_t \ot \tau_t)\de = \de \tau_t$.
\end{demo}

A last commutation property concerns $R$ and $\de$.

\begin{proposition}
We have that $\flip(R \ot R)\de = \de R$.
\end{proposition}
\begin{demo}
Choose $\om \in \cop^*$. Then
\begin{eqnarray*}
& & (\flip(R \ot R)\de)(\,(\io \ot \om)(W)\,) \\
& & \spat = \flip(\,(I \ot I) \de(\,(\io \ot \om)(W)\,)^* (I \ot I)\,) \\
& & \spat = \flip(\,(I \ot I) [W^* (1 \ot (\io \ot \om)(W))\,W]^* (I \ot I)\, 
) \\
& & \spat = \flip(\,(I \ot I) W^* (1 \ot (\io \ot \overline{\om})(W^*)) W (I 
\ot I) \,)
\\
& & \spat = \flip(\,(I \ot I)(\io \ot \io \ot \overline{\om})(W_{12}^* 
W_{23}^* W_{12})
(I \ot I)\,) \\
& & \spat = \flip(\,(I \ot I)(\io \ot \io \ot \overline{\om})(\,(W_{13} 
W_{23})^*)
(I \ot I)\,).
\end{eqnarray*}
where, in the last equality, we used the pentagon equation (proposition
\ref{prop1.3}).
Using this equality, it is not so difficult to check that
\begin{eqnarray*}
& & (\flip(R \ot R)\de)(\,(\io \ot \om)(W)\,) \\
& & \spat =\flip(\,(\io \ot \io \ot \om)\bigl( (\,(I\,.^*I) \ot (I\,.^*I) \ot
\io)(W_{13} W_{23})\bigr) \,)\\
& & \spat =\flip(\,(\io \ot \io \ot \om)\bigl((\,(I\,.^*I) \ot \io)(W)_{13} \,
(I\,.^*I) \ot \io)(W)_{23}\bigr)\,).
\end{eqnarray*}
Because $(I \ot J)\,W = W^* (I \ot J)$, we have that
$$(\,(I\,.^*I) \ot \io)(W) = (\io \ot (J\,.^*J)\,)(W)  .$$
Hence,
\begin{eqnarray*}
& & (\flip(R \ot R)\de)(\,(\io \ot \om)(W)\,) \\
& & \spat = \flip(\,(\io \ot \io \ot \om)\bigl(\,(\io \ot 
(J\,.^*J)\,)(W)_{13} \,
(\io \ot (J\,.^*J)\,)(W)_{23}\bigr)\,) \\
& & \spat = \flip(\,(\io \ot \io \ot \om)\bigl(\,(\io \ot \io \ot 
(J\,.^*J)\,)(W_{23}
W_{13})\,\bigr)\,) \\
& & \spat = (\io \ot \io \ot \om)\bigl(\,(\io \ot \io \ot (J\,.^*J)\,)(W_{13}
W_{23})\,\bigr) \\
& & \spat = (\io \ot \io \ot \om)\bigl(\,(\io \ot \io \ot 
(J\,.^*J)\,)(W_{12}^*
W_{23} W_{12})\,\bigr) \\
& & \spat = W^* (\io \ot \io \ot \om)\bigl(\,(\io \ot \io \ot 
(J\,.^*J)\,)(W_{23})\,\bigr)
W \\
& & \spat = W^* \bigl(1 \ot (\io \ot \om(J\,.^*J)\,)(W)\,\bigr) W \\
& & \spat = W^* \bigl(1 \ot R(\,(\io \ot \om)(W)\,)\,\bigr) W
= \de\bigl(R(\,(\io \ot \om)(W)\,)\,\bigr).
\end{eqnarray*}
The proposition follows by equation \ref{eq1.2} of section 2.
\end{demo}
In a later section, we will look at some other implementations of $\tau$ and 
$R$.

\section{The left Haar weight.}

In this section, we want to extend our positive Haar functional on the
$^*$-algebra $A$ to a weight on the \cst\ $A_r$.
We will show that this weight satisfies the KMS-condition and that it is
left invariant in some sense.

\medskip

First, we are going to define the left Haar weight on our \cst\ $A_r$.
In order to make our notation not too heavy, we will use for this weight
the same symbol $\vfi$ as for our left Haar functional on our $^*$-algebra 
$A$.

The context will make clear which meaning $\vfi$ has. The term \lq linear 
functional 
$\vfi$\rq \ always means the left invariant funcional $\vfi$ on the $^*$-
algebra $A$ whereas
the term \lq weight $\vfi$\rq \ always means the left invariant weight $\vfi$ 
on the 
C$^*$-algebra $A_r$. Also, when $\vfi$ works on elements of $A$, $\vfi$ 
will always be the linear fuctional $\vfi$.  When $\vfi$ works on elements of 
$A_r$, $\vfi$ will always be the weight $\vfi$. 

Later in this paper, we also have other symbols with two meanings. In every 
case, it handles
about a symbol which denotes an object connected to the $^*$-algebra $A$ and 
an object connected to the C$^*$-algebra $A_r$ which are each other 
equivalent. In these cases, we will use the same conventions as above.

\begin{definition}
Denote $\tilde{\vfi}$ as the semifinite faithful normal weight on
the von Neumann algebra $A_r''$ associated to the left Hilbert algebra $\cu$.
Then we define $\vfi$ as the restriction of $\tilde{\vfi}$ to $A_r^+$,
hence $\vfi$ is a faithful lower semi-continous weight on $A_r$.
\end{definition}

The next proposition guarantees that our $\vfi$ on the $C^*$-algebra
level is an extension of our $\vfi$ on the $^*$-algebra level.

\begin{proposition}
We have that $\pi(A)$ is a subset of $\Mfi$ and $\vfi(\pi(a)) = \vfi(a)$
for all $a \in A$.
\end{proposition}
\begin{demo}
In this proof, we will use the theory of left Hilbert algebras.

Choose $b \in A$. Then $L_{\la(b)} = \pi(b)$. According to the theory of left
Hilbert algebras, $L_{\la(b)}^* L_{\la(b)}$ belongs to ${\cal
M}_{\tilde{\vfi}}^+$
and $\tilde{\vfi}(L_{\la(b)}^* L_{\la(b)}) =  \langle \la(b),\la(b) \rangle  
= \vfi(b^* b)$.
Because $L_{\la(b)}^* L_{\la(b)} = \pi(b)^* \pi(b) = \pi(b^* b)$,
and $\vfi$ is the restriction of $\tilde{\vfi}$ to $A_r^+$, it follows
that $\pi(b^* b)$ belongs to $\Mfi$ and $\vfi(\pi(b^* b)) = \vfi(b^* b)$.

The proposition follows from polarization and the fact that $A = A^* A$.
\end{demo}

From this proposition, we can conclude dat $\vfi$ is densely defined.

\medskip

From now on, we will use the following GNS-triple of $\vfi$:
\begin{enumerate}
\item $H_\vfi = H$
\item We define the linear map $\lafi$ from $\Nfi$ into $H$ as follows.
      Let $a \in \Nfi$, then there exists a unique $v \in H$ such
      that $v$ is left bounded with respect to $\cu$ and we define
      $\lafi(a) = v$.
\item $\pi_\vfi =$ the identity map on $A_r$
\end{enumerate}
Remark that $\lafi(\pi(a)) = \la(a)$ for all $a \in A$, this guarantees in
fact that $\lafi$ has dense range in $H$.

\medskip

The following result follows immediately from the theory of left Hilbert
algebras.

\begin{proposition}
We have that $\vfi$ is invariant under $\si$. Moreover,
$\lafi(\si_t(x)) = \nab^{it} \lafi(x)$ for all $t \in \R$.
\end{proposition}

\medskip

Now, we are going to prove that the weight $\vfi$ is determined by its
values on $\pi(A)$. More precisely, we will prove that $\la(A)$ is a core
for $\lafi$.

\medskip

We define $\la_0$ as the closure of the mapping $\pi(A) \rightarrow H : x
\mapsto \lafi(x)$
and $A_0$ will be the domain of $\la_0$. It is clear that $\la_0$ is a
restriction of $\lafi$.
Furthermore, it is easy to check that $A_0$ is a left ideal of $A_r$.

\begin{lemma}
Consider $x,y \in \pi(A)$ and $\om \in A_r^*$. Then $(y \om \ot \io)\de(x)$
belongs to $\pi(A)$ and \newline
$\vfi(\,(y \om \ot \io)\de(x)\,) = \om(y) \, \vfi(x)$.
\end{lemma}
\begin{demo}
Choose $a,b \in A$.
Then
\begin{eqnarray*}
& & (\pi(b) \om \ot \io)\de(\pi(a)) = (\om \ot \io)(\de(\pi(a))(\pi(b) \ot
1)\,) \\
& & \spat = (\om \ot \io)(\,(\pi \od \pi)(\de(a)(b \ot 1))\,)
= \pi(\,(\,\om \circ \pi \od \io)(\de(a)(b \ot 1))\,).
\end{eqnarray*}
So, $(\pi(b) \om \ot \io)\de(\pi(a))$ belongs to $\pi(A)$ and the left
invariance of $\vfi$ on the $^*$-algebra level implies that
$$ \vfi(\,(\pi(b) \om \ot \io)\de(\pi(a))\,)
= \vfi(a)\, (\om \circ \pi)(b) = \vfi(\pi(a)) \, \om(\pi(b)).$$
\end{demo}

\begin{lemma}
Consider $x \in A_0$ and $\om \in (A_r)_+^*$. Then $(\om \ot \io)\de(x^* x)$
belongs to $\Mfi^+$ and $\vfi(\,(\om \ot \io)\de(x^* x)\,) \leq
\|\om\|\,\vfi(x^* x)$.
\end{lemma}
\begin{demo}
It is a well known result that there exist $\th \in (A_r)^*_+$ and $y \in A_r$
such that $\om = y \th y^*$ (see \cite{Tay}).
Now, there exist sequences $(x_n)_{n=1}^\infty , (y_n)_{n=1}^\infty$ in
$\pi(A)$ such that $(y_n)_{n=1}^\infty \rightarrow y$,
$(x_n)_{n=1}^\infty \rightarrow x$ and $(\lafi(x_n)\,)_{n=1}^\infty 
\rightarrow
\lafi(x)$.

It is clear that $(\,(y_n \th y_n^* \ot \io)\de(x_n^* x_n)\,)_{n=1}^\infty
\rightarrow (\om \ot \io)\de(x^* x)$.

The lower semicontinuity of $\vfi$ implies that
$$ \vfi(\,(\om \ot \io)\de(x^* x)\,)
\leq \liminf \bigl(\vfi(\,(y_n \th y_n^* \ot \io)\de(x_n^* x_n)\,)\,
\bigr)_{n=1}^\infty .$$
By the previous lemma, the right hand side of this equality equals
\begin{eqnarray*}
& & \liminf (\th(y_n^* y_n) \vfi(x_n^* x_n)\,)_{n=1}^\infty
= \liminf ( \|y_n \th y_n^*\|  \langle \lafi(x_n),\lafi(x_n) \rangle 
)_{n=1}^\infty \\
& & \spat = \|\om\|  \langle \lafi(x),\lafi(x) \rangle  = \|\om\| \, \vfi(x^* 
x).
\end{eqnarray*}
The lemma follows.
\end{demo}

\begin{lemma}  \label{equi1}
Consider $x \in A_0$ and $\om \in A_r^*$. Then
$(\om \ot \io)\de(x)$ belongs to $\Nfi$ and $\|\lafi(\,(\om \ot 
\io)\de(x)\,)\|
\leq \|\om\|\,\|\lafi(x)\|$.
\end{lemma}
\begin{demo}
We know that there exist an element $\th \in (A_r)_+^*$ with
$\|\th\|=\|\om\|$ such that $$ [(\om \ot \io)\de(y)]^* [(\om \ot \io)\de(y)]
\leq \|\om\| (\th \ot \io)\de(y^* y) $$
for all $y \in A_r$ (see proposition 4.6 of \cite{Tak}).
Hence, using the previous lemma, we see that
\begin{eqnarray*}
& & \vfi( [(\om \ot \io)\de(x)]^* [(\om \ot \io)\de(x)] \,)
\leq \|\om\| \, \vfi(\,(\th \ot \io)\de(x^* x)\,)    \\
& & \spat \leq \|\om\| \, \|\th\| \, \vfi(x^* x) = \|\om\|^2 \, \vfi(x^* x) .
\end{eqnarray*}
The lemma follows.
\end{demo}

Finally, we can prove the lemma which will be essential to us:

\begin{lemma}  \label{prop5.3}
Consider $\om \in A_r^*$ and $x \in A_0$. Then $(\om \ot \io)\de(x)$
belongs to $A_0$.
\end{lemma}
\begin{demo}
There exist $\th \in A_r^*$ and $y \in A_r$ such that $\om = y \th$.
Furthermore, there exist sequences $(x_n)_{n=1}^\infty$ and
$(y_n)_{n=1}^\infty$ in $\pi(A)$  such that $(y_n)_{n=1}^\infty \rightarrow 
y$,
$(x_n)_{n=1}^\infty \rightarrow x$ and
$(\lafi(x_n)\,)_{n=1}^\infty \rightarrow \lafi(x)$.

We know already that $(y_n \th \ot \io)\de(x_n)$ belongs to $\pi(A)$ for
all $n \in \N$.
It is also clear that  \newline $(\,(y_n \th \ot \io)\de(x_n)\,)_{n=1}^\infty
\rightarrow (\om \ot \io)\de(x)$.

There exists a positive number M such that $\|y_n \th\| \leq M$
and $\|\lafi(x_n)\| \leq M$ for all $n \in \N$.

So, for all $m,n \in \N$, we have that
\begin{eqnarray*}
& & \| \lafi(\,(y_n \th \ot \io)\de(x_n)\,) -\lafi(\,(y_m \th \ot \io)\de(x_m)
\,)  \| \\
& & \spat \leq \|\lafi(\,(y_n \th \ot \io)\de(x_n)\,)
-\lafi(\,(y_n \th \ot \io)\de(x_m)\,) \| \\
& & \spat + \|\lafi(\,(y_n \th \ot \io)\de(x_m)\,)
-\lafi(\,(y_m \th \ot \io)\de(x_m)\,) \| \\
& & \spat \leq \|y_n \th\| \, \|\lafi(x_n) - \lafi(x_m)\|
+ \|\lafi(x_m)\| \, \| y_n \th - y_m \th \|   \text{\ \ \ (a)} \\
& & \spat \leq M \|\lafi(x_n) - \lafi(x_m)\|
+ \| y_n \th - y_m \th \| M,
\end{eqnarray*}
where in inequality (a), we used the previous lemma.

From this chain of inequalities, we infer that
$\bigl(\lafi(\,(y_n \ot \io)\de(x_n)\,)\,\bigr)_{n=1}^\infty$ is Cauchy and
hence convergent in $H$.
Because $\la_0$ is closed, $(\om \ot \io)\de(x)$ belongs to $A_0$.
\end{demo}

Now, we will prove a result which is essential in the rest of the paper.
First, we need an easy lemma.

\begin{lemma}
Consider elements $a,b,c \in A$. Then
\begin{enumerate}
\item $(\om_{\la(b),\la(c)} \ot \io)\de(\pi(a))
      = \pi(\,(\vfi \od \io)((c^* \ot 1)\de(a)(b \ot 1))\,)$
\item $(\io \ot \om_{\la(b),\la(c)})\de(\pi(a))
      = \pi(\,(\io \od \vfi)((1 \ot c^*)\de(a)(1 \ot b))\,) .$
\end{enumerate}
\end{lemma}
\begin{demo}
We know there exist an element $d \in A$ such that $b=db$.
Hence,
\begin{eqnarray*}
& & \pi(\,(\vfi \od \io)((c^* \ot 1)\de(a)(b \ot 1))\,)
= \pi(\,(\vfi \od \io)((c^* \ot 1)\de(a)(d \ot 1)(b \ot 1))\,) \\
& & \spat = (\om_{\la(b),\la(c)} \ot \io)(\,(\pi \od \pi)(\de(a)(d \ot 1))\,) 
\\
& & \spat = (\om_{\la(b),\la(c)} \ot \io)(\de(\pi(a))(\pi(d) \ot 1)\,) \\
& & \spat = (\om_{\pi(d) \la(b),\la(c)} \ot \io)(\de(\pi(a))\,) \\
& & \spat = (\om_{\la(db),\la(c)} \ot \io)\de(\pi(a)) 
= (\om_{\la(b),\la(c)} \ot \io)\de(\pi(a)).
\end{eqnarray*}
The other equality is proven in a similar way.
\end{demo}

\begin{proposition} \label{prop5.1}
Consider a dense left ideal $N$ in $A_r$ such that for all $a \in N$ and
all $\om \in A_r^*$, we have that $(\io \ot \om)\de(a) \in N$.
Then $\pi(A)$ is a subset of $N$.
\end{proposition}
\begin{demo}
Choose $b \in A$.
Because $\vfi \neq 0$, there exist $c,d \in A$ such that $c \vfi d^* \neq 0$.
We have that $\om_{\la(b),\la(c)} \circ \pi = c \vfi d^*$, therefore
$\om_{\la(b),\la(c)} \neq 0$.  The fact that $N$ is dense in $A_r$ implies
the existence of an element $x \in N$ such that $\om_{\la(b),\la(c)}(x)=1$.

There exist $p_1,\ldots\!,p_n \, , \, q_1,\ldots\!,q_n \, , \, 
r_1,\ldots\!,r_n
\in A$ such  that
$$ (b \ot 1)\de(\rho^{-1}(c) d^*) = \sum_{i=1}^n (1 \ot p_i q_i)\de(r_i) .
\text{\ \ \ (*)} $$
For every $y \in A$, we have that
\begin{eqnarray*}
& & \om_{\la(c),\la(d)} (\pi(y)) \, \pi(b) = \vfi(d^* y c) \, \pi(b)
= \pi(\,\vfi(\rho^{-1}(c) d^* y) \, b) \\
& & \spat = \pi(\,(\io \od \vfi)((b \ot 1) \de(\rho^{-1}(c) d^* y))\,) \\
& & \spat =\sum_{i=1}^n \pi(\,(\io \od \vfi)((1 \ot p_i q_i) \de(r_i y))\,) \\
& & \spat = \sum_{i=1}^n \pi(\,(\io \od \vfi)((1 \ot q_i) \de(r_i y)(1 \ot
\rho(p_i)))\,) \\
& & \spat = \sum_{i=1}^n (\io \ot \om_{\la(\rho(p_i)),\la(q_i^*)})\de(\pi(r_i 
y))\\
& & \spat = \sum_{i=1}^n (\io \ot 
\om_{\la(\rho(p_i)),\la(q_i^*)})\de(\pi(r_i) \pi(
y))\,).
\end{eqnarray*}
Because $\pi(A)$ is dense in $A_r$, we must have that
$$ \om_{\la(c),\la(d)}(z) \, \pi(b)
= \sum_{i=1}^n (\io \ot \om_{\la(\rho(p_i)),\la(q_i^*)})\de(\pi(r_i)\,z) . $$
for all $z \in A_r$.
In particular, we have that
$$ \pi(b) = \om_{\la(c),\la(d)}(x) \, \pi(b)
= \sum_{i=1}^n (\io \ot \om_{\la(\rho(p_i)),\la(q_i^*)})\de(\pi(r_i)\,x). $$
Because $x$ belongs to $N$ and $N$ is a left ideal in $A_r$,
$\pi(r_i) \,x$ will belong to $N$ for all $i \in \{1,\ldots\!,n\}$.
By the above equality, we see that $\pi(b)$ belongs to $N$.
\end{demo}

Of course, also the following proposition is true.

\begin{proposition}  \label{prop5.2}
Consider a dense left ideal $N$ in $A_r$ such that for all $a \in N$ and
all $\om \in A_r^*$, we have that $(\om \ot \io)\de(a) \in N$.
Then $\pi(A)$ is a subset of $N$.
\end{proposition}

\bf Sketch of proof: \rm
The proof starts in the same way as the previous one but step (*)
is replaced by :

There exist $p_1,\ldots\!,p_n \, , \, q_1,\ldots\!,q_n \, , \, r_1,
\ldots\!,r_n \in A$ such that
$$(1 \ot b \sde^{-1}) \de(\rho^{-1}(c) d^*)
= \sum_{i=1}^n (p_i q_i \ot 1)\de(r_i) .$$
Then you prove in a similar way as in the previous proof that
$$ \om_{\la(c),\la(d)}(\pi(y)) \, \pi(b)
= \sum_{i=1}^n (\om_{\la(\rho(p_i)),\la(q_i^*)} \ot \io)\de(\pi(r_i)
\pi(y))$$ for every $y \in A$ but in stead of the left invariance of $\vfi$,
you use the fact that $(\vfi \od \io)\de(a) = \vfi(a) \sde$
for any $a \in A$.
From there one, the proof runs along the same lines as the preceding one.

\medskip

Now, we will give a first application of this result.

\begin{lemma}  \label{lem5.1}
Consider $t \in \R$, then $\si_t(A_0)$ is a subset of $A_0$.
\end{lemma}
\begin{demo}
First, we proof that $\si_t(\pi(A))$ is a subset of $A_0$.

\begin{list}{}{\setlength{\leftmargin}{.4cm}}

\item Because $A_0$ is a dense left ideal in $A_r$, $\si_{-t}(A_0)$ is a dense
left ideal in $A_r$.

Choose $a \in A_0$ and $\om \in A_r^*$. Using proposition \ref{prop4.1}, we 
see
that
\begin{eqnarray*}
& & (\om \ot \io)\de(\si_{-t}(a)) = (\om \ot \io)(\tau_{-t} \ot \si_{-t})
\de(a) \\
& & \spat = \si_{-t}(\,(\om \tau_{-t} \ot \io)\de(a)\,).
\end{eqnarray*}
Lemma \ref{prop5.3} implies that $(\om \tau_{-t} \ot \io)\de(a)$
belongs to $A_0$, so $(\om \ot \io)\de(\si_{-t}(a))$ belongs to $\si_{-
t}(A_0)$.

By proposition \ref{prop5.2}, we have that $\pi(A) \subseteq \si_{-t}(A_0)$,
so $\si_t(\pi(A)\,) \subseteq A_0$.

\end{list}

Choose $x \in A_0$. Then there exists a sequence $(x_n)_{n=1}^\infty$
in $\pi(A)$ such that $(x_n)_{n=1}^\infty \rightarrow x$
and $(\lafi(x_n)\,)_{n=1}^\infty \rightarrow \lafi(x)$.

From the first part of the proof, we know that $\si_t(x_n)$ belongs to $A_0$
for all $n \in \N$.
Furthermore, it is clear that $(\si_t(x_n)\,)_{n=1}^\infty \rightarrow
\si_t(x)$. Because $\si$ is invariant under $\vfi$,
we know that
$$\| \lafi(\si_t(x_m)) - \lafi(\si_t(x_n)) \|
= \| \lafi(x_m) - \lafi(x_n) \| $$
for every $m,n \in \N$.
This implies that $(\lafi(\si_t(x_n))\,)_{n=1}^\infty$ is Cauchy
and hence convergent in $H$. The closedness of $\la_0$ implies that
$\si_t(x)$ belongs to $A_0$.
\end{demo}

We are now able to prove a major result of this paper which says that the 
left invariant
weight is completely determined by its values on $\pi(A)$

\begin{theorem}  \label{thm5.1}
We have that $\pi(A)$ is a core for $\lafi$.
\end{theorem}
\begin{demo}
Because $\pi(A)$ is dense in $A_r$, there exists a bounded net $(e_j)_{j \in 
J}$
in $\pi(A)$ such that $(e_j)_{j \in J}$ converges strictly to 1.

For every $j \in J$, we define
$$u_j = \frac{1}{\sqrt{\pi}} \int \exp(-t^2) \,
\si_t(e_j) \, dt  \in A_r,$$ it is clear that $u_j$ belongs to
$D(\si_\frac{i}{2})$ and $$\si_\frac{i}{2}(u_j)
= \frac{1}{\sqrt{\pi}} \int \exp(-(t - \frac{i}{2})^2) \, \si_t(e_j) \,dt.$$

Furthermore, $(u_j)_{j \in J}$ and $(\si_\frac{i}{2}(u_j)\,)_{j \in J}$
are both bounded nets in $A_r$ which converge strictly to 1, because
$(e_j)_{j \in J}$ converges strictly to 1.

Next, we show that $u_j$ belongs to $A_0$  for all $j \in J$.

\begin{list}{}{\setlength{\leftmargin}{.4cm}}

\item From lemma \ref{lem5.1}, we know that $\si_t(e_j)$ belongs to $A_0$ for 
every
$t \in \R$. Also, $\la_0(\si_t(e_j)) = \nab^{it} \la_0(e_j)$ for every $t \in
\R$.  Therefore, the function $\R \rightarrow H : t \mapsto \exp(-t^2) \, 
\la_0
(\si_t(e_j))$ is integrable. The closedness of $\la_0$ implies
that $u_j$ belongs to $A_0$.

\end{list}

Choose $x \in \Nfi$. Then $(x u_j)_{j \in J} \rightarrow x$.

For all $j \in J$, we have that $x u_j$ belongs to $A_0$ (because $A_0$ is a
left ideal) and $$\la_0(x u_j) = \lafi(x u_j) = J \si_\frac{i}{2}(u_j)^* J
\lafi(x)$$ (here, we used lemma \ref{lem5.2}).
Therefore, $(\la_0(x u_j)\,)_{j \in J} \rightarrow \lafi(x)$.

Because $\la_0$ is closed, we see that $x \in A_0$ and $\la_0(x) = \lafi(x)$.

At the end, we see that $A_0 = \Nfi$, so $\pi(A)$ is a core for $\lafi$.
\end{demo}

This result allows us to prove the left invariance of $\vfi$ on the
\cst\ level.

\begin{theorem}
Consider $x \in \Mfi$, then $\de(x)$ belongs to $\overline{{\cal M}}_{\io \ot
\vfi}$ and $(\io \ot \vfi)\de(x) = \vfi(x) \, 1$.
\end{theorem}
\begin{demo}
Choose $y \in \Nfi$.

Take $z \in A_r$. Then there exist sequences $(a_n)_{n=1}^\infty$ and
$(b_n)_{n=1}^\infty$ in $A$ such that $(\pi(a_n)\,)_{n=1}^\infty \rightarrow
y$, \ $(\pi(b_n)\,)_{n=1}^\infty \rightarrow z$ and
$(\lafi(\pi(a_n)\,)_{n=1}^\infty \rightarrow \lafi(y)$.

Choose $c_1,c_2,d_1,d_2 \in A$. Then
\begin{eqnarray*}
& &  \langle (\io \ot \lafi)(\de(\pi(c_1))(\pi(d_1)\ot 1)\,) ,
(\io \ot \lafi)(\de(\pi(c_2))(\pi(d_2) \ot 1)\,) \rangle      \\
& & \spat =  \langle (\io \ot \lafi)(\,(\pi \od \pi)(\de(c_1)(d_1 \ot 1))\,) ,
(\io \ot \lafi)(\,(\pi \od \pi)(\de(c_2)(d_2 \ot 1))\,) \rangle    \\
& & \spat = (\pi \od \vfi)(\,(d_2^* \ot 1) \de(c_2^* c_1) (d_1 \ot 1)\,)
= \vfi(c_2^* c_1) \pi(d_2^* d_1) \\
\end{eqnarray*}
Therefore,
\begin{eqnarray*}
& &  \langle (\io \ot \lafi)(\de(\pi(c_1))(\pi(d_1) \ot 1)\,) ,
(\io \ot \lafi)(\de(\pi(c_2))(\pi(d_2) \ot 1)\,) \rangle      \\
& & \spat\spat =  \langle \lafi(\pi(c_1)) , \lafi(\pi(c_2)) \rangle  \, 
\pi(d_2)^* \pi(d_1) \ .
\hspace{3.5cm} \text{(a)}
\end{eqnarray*}

It is clear that $(\de(\pi(a_n))(\pi(b_n)\ot 1)\,)_{n=1}^\infty \rightarrow
\de(y)(z \ot 1)$.
Using equation (a), we get for any $m,n \in \N$ that
\begin{eqnarray*}
& & \| (\io \ot \lafi)(\de(\pi(a_m))(\pi(b_m) \ot 1) \,) -
(\io \ot \lafi)(\de(\pi(a_n))(\pi(b_n) \ot 1)\,) \|^2 \\
& & \spat = \|  \langle \lafi(\pi(a_m)) , \lafi(\pi(a_m)) \rangle  \,
\pi(b_m)^* \pi(b_m) 
-   \langle \lafi(\pi(a_m)) , \lafi(\pi(a_n)) \rangle  \,
\pi(b_m)^* \pi(b_n) \\
& & \spat -   \langle \lafi(\pi(a_n)) , \lafi(\pi(a_m)) \rangle  \,
\pi(b_n)^* \pi(b_m) 
+   \langle \lafi(\pi(a_n)) , \lafi(\pi(a_n)) \rangle  \,
\pi(b_n)^* \pi(b_n) \|.
\end{eqnarray*}
This equation implies that
$\bigl(\,(\io \ot \lafi)(\de(\pi(a_n))(\pi(b_n) \ot 1)\,\bigr)_{n=1}^\infty$
is a Cauchy sequence and hence convergent in $H \ot A$.
The closedness of $\io \ot \lafi$ implies that
$\de(y)(z \ot 1)$ belongs to ${\cal N}_{\io \ot \vfi}$
and that \newline
$\bigl(\,(\io \ot \lafi)(\de(\pi(a_n))(\pi(b_n) \ot 
1)\,)\,\bigr)_{n=1}^\infty$
converges to $(\io \ot \lafi)(\de(y)(z \ot 1)\,)$.

By equation (a), we have that
\begin{eqnarray*}
& &  \langle (\io \ot \lafi)(\de(\pi(a_n))(\pi(b_n) \ot 1)\,)   ,
(\io \ot \lafi)(\de(\pi(a_n))(\pi(b_n) \ot 1)\,)   \rangle  \\
& & \spat\spat =   \langle \lafi(\pi(a_n)) , \lafi(\pi(a_n)) \rangle  \,
\pi(b_n)^* \pi(b_n)
\end{eqnarray*}
for every $n \in \N$.
So,
$$ \bigl(\, \langle (\io \ot \lafi)(\de(\pi(a_n))(\pi(b_n) \ot 1)\,)   ,
(\io \ot \lafi)(\de(\pi(a_n))(\pi(b_n) \ot 1)\,)   \rangle   
\bigr)_{n=1}^\infty $$
converges to $ \langle \lafi(a),\lafi(a) \rangle  \, b^* b$.
Hence,
\begin{eqnarray*}
& &  \langle (\io \ot \lafi)(\de(y)(z \ot 1)\,)   ,
(\io \ot \lafi)(\de(y)(z \ot 1)\,)   \rangle  \\
& & \spat\spat =   \langle \lafi(y) , \lafi(y) \rangle  \, z^* z.
\end{eqnarray*}
This implies that
$$ (\io \ot \vfi)(\,(z^* \ot 1)\de(y^* y)(z \ot 1)\,)
= \vfi(y^* y) \, z^* z. $$

For every $\om \in {\cal G}_\vfi$, we have that
$$z^*\, (\io \ot \om)(\de(y^* y)) \, z
= (\io \ot \om)(\,(z^* \ot 1)\de(y^* y)(1 \ot z)\,). \text{\ \ \ (a)}$$
By the previous discussion, we know that
$$\bigl(\,(\io \ot \om)((z^* \ot 1)\de(y^* y)(1 \ot z))\,\bigr)_{\om \in {\cal
G}_\vfi} \rightarrow  \vfi(y^*y) \, z^*z .$$
Therefore, equation (a) implies that
$$\bigl(\,z^*\, (\io \ot \om)(\de(y^* y)) \, z\,\bigr)_{\om \in {\cal
G}_\vfi} \rightarrow  \vfi(y^*y) \, z^*z .$$

Consequently, we get that $\de(y^*y)$ belongs to
$\overline{{\cal M}}_{\io \ot
\vfi}$ and $(\io \od \vfi)\de(y^* y) = \vfi(y^* y) 1$.

The proposition follows by polarization.
\end{demo}

\begin{corollary}
Consider $x \in \Mfi$ and $\om \in A_r^*$. Then $(\om \ot \io)\de(x)$
belongs to $\Mfi$ and $\vfi(\,(\om \ot \io)\de(x)\,) = \vfi(x) \, \om(1)$.
\end{corollary}

It is also possible to prove the strong left invariance proposed in the 
definition of Masuda, Nakagami \& Woronowicz. This is in fact a direct 
consequence of formula \ref{eq1.3}.

\begin{proposition}   \label{prop5.4}
Consider $a,b \in \Nfi$ and $\om \in A_r^*$ such that $\om R
\tau_{-\frac{i}{2}}$
is bounded and call $\th$ the unique element in $A_r^*$ which extends
$\om R \tau_{-\frac{i}{2}}$. Then
$b^* \, (\om \ot \io)\de(a)$
and $(\th \ot \io)(\de(b^*)) \,a$ belong to
$\Mfi$ and
$$\vfi(\,b^* \, (\om \ot \io)\de(a)\,)
= \vfi(\,(\th  \ot \io)(\de(b^*)) \,a \,) \ .$$
\end{proposition}
\begin{demo}
Because of the left invariance of $\vfi$, it follows that $(\om \ot 
\io)\de(a)$
belongs to $\Nfi$ and $(\th \ot \io)(\de(b^*))$
belongs to $\Nfi^*$. This implies that
$b^* \, (\om \ot \io)\de(a)$
and $(\th \ot \io)(\de(b^*)) \,a$ belong to
to $\Mfi$.

Because $\pi(A)$ is a core for $\lafi$, there exist sequences
$(a_n)_{n=1}^\infty$ and $(b_n)_{n=1}^\infty$ in $A$ such
that $(\pi(a_n)\,)_{n=1}^\infty \rightarrow a$, \
$(\pi(b_n)\,)_{n=1}^\infty \rightarrow b$, \
$(\lafi(\pi(a_n))\,)_{n=1}^\infty \rightarrow \lafi(a)$ and
$(\lafi(\pi(b_n))\,)_{n=1}^\infty \rightarrow \lafi(b)$.

The left invariance of $\vfi$ implies that $\bigl(\,\lafi(\,(\om \ot
\io)\de(\pi(a_n))\,)\,\bigr)_{n=1}^\infty \rightarrow \lafi(\,(\om
\ot \io) \de(a)\,)$ and \newline
$\bigl(\,\lafi(\,(\overline{\th} \ot \io)\de(\pi(b_n))\,)\,\bigr)_{n=1}^\infty
\rightarrow \lafi(\,(\overline{\th} \ot \io)\de(b)\,)$.

This implies that
$$\bigl(\, \langle \lafi(\,(\om \ot \io)\de(\pi(a_n))\,) , \lafi(\pi(b_n))
\rangle \bigr)_{n=1}^\infty$$
converges to $\vfi(\,b^* (\om \ot \io)\de(a)\,)$ \ \ (a) \newline
and
$$\bigl(\, \langle \lafi(\pi(a_n)) ,  \lafi(\,(\overline{\th} \ot
\io)\de(\pi(b_n)) \,) \rangle \bigr)_{n=1}^\infty$$
converges to $\vfi(\,(\th \ot \io)(\de(b^*)) \, a\,)$ \ \ (b).

Using equation \ref{eq1.3}, we have for every $n \in \N$ that
\begin{eqnarray*}
& & \hspace{-0.7cm} \langle \lafi(\,(\om \ot \io)\de(\pi(a_n))\,) , 
\lafi(\pi(b_n)) \rangle
= \vfi(\,\pi(b_n^*) (\om \ot \io)\de(\pi(a_n)) \,) \\
& & \hspace{-0.7cm}\spat = \vfi(\,(\om \ot \io)(\,(1 \ot 
\pi(b_n^*))\de(\pi(a_n))\,)\,)
= \vfi(\,(\om \ot \io)(\,(\pi \od \pi)((1 \ot b_n^*)\de(a_n))\,)\,) \\
& & \hspace{-0.7cm}\spat = \om(\pi(\,(\io \od \vfi)((1 \ot b_n^*)\de(a_n))\,))
= \om(\pi(\,S(\,(\io \od \vfi)(\de(b_n^*)(1 \ot a_n))\,)\,)) \\
& & \hspace{-0.7cm}\spat = \om(R(\tau_{-\frac{i}{2}}(\pi(\,(\io \od 
\vfi)(\de(b_n^*)(1 \ot
a_n))\,))))
= \th(\pi(\,(\io \od \vfi)(\de(b_n^*)(1 \ot a_n))\,)) \\
& & \hspace{-0.7cm}\spat = \vfi(\,(\th \ot \io)(\,(\pi \od \pi)(\de(b_n^*)(1 
\ot a_n))\,)\,)
= \vfi(\,(\th \ot \io)(\de(\pi(b_n^*))(1 \ot \pi(a_n)))\,) \\
& & \hspace{-0.7cm}\spat = \vfi(\,(\th \ot \io)(\de(\pi(b_n^*)) \, \pi(a_n)\,)
= \langle \lafi(\pi(a_n)) , \lafi(\,(\overline{\th} \ot \io)\de(\pi(b_n))\,)
\rangle. \end{eqnarray*}
Together with (a) and (b), this gives us that
$$\vfi(\,b^* \, (\om \ot \io)\de(a)\,)
= \vfi(\,(\th  \ot \io)(\de(b^*)) \,a \,) \ .$$
\end{demo}

\section{Invariance properties of bi-C$^*$-isomorphisms and group-like 
elements.}

In this section, we will consider group-like elements and $^*$-isomorphisms 
which commute with the comultiplication.  In particular, their behaviour with 
respect to the 
one-parameter groups $\si$ and $\tau$ , as well as their behaviour with 
respect to the anti-unitary
antipode $R$, is investigated. The propositions involved will mostly handle 
about some
form of relative invariance.  We will use a lot of these results for the 
first time in the next section, where we introduce the modular function of 
our \cst ic quantum group $(A_r,\de)$.

\begin{proposition}  \label{prop6.2}
Consider a $^*$-automorphism $\al$ on $A_r$ such that there exists
a $^*$-automorphism $\be$ on $A_r$ such that $\de \al = (\be \ot \al)\de$.
Then there exists a unique strictly positive number $r$ such that
$\vfi \al = r \vfi$.
\end{proposition}
\begin{demo}
Because $\Nfi$ is a dense left ideal of $A_r$, $\al^{-1}(\Nfi)$ is
a dense left ideal of $A_r$.
Choose $\om \in A_r^*$ and $a \in \al^{-1}(\Nfi)$. Then
\begin{eqnarray*}
& & \al(\,(\om \ot \io)\de(a)\,) = (\om \be^{-1} \ot \io)(\be \ot \al)\de(a) 
\\
& & \spat\spat (\om \be^{-1} \ot \io)\de(\al(a)).
\end{eqnarray*}
Because $\al(a)$ belongs to $\Nfi$ and because of the left invariance of 
$\vfi$,
we have that $(\om \be^{-1} \ot \io)\de(\al(a))$ belongs to $\Nfi$.
Therefore, $(\om \ot \io)\de(a)$ belongs to $\al^{-1}(\Nfi)$.
Proposition \ref{prop5.2} implies that $\pi(A) \subseteq \al^{-1}(\Nfi)$,
hence, $\al(\pi(A)) \subseteq \Nfi$.
Now, using the fact that $A^* A=A$, it is not difficult to see that
$\al(\pi(A)) \subseteq \Mfi$.

This allows us to define the positive linear functional $\phi$ on $A$
such that $\phi(a) = \vfi(\al(\pi(a))\,)$ for every $a \in A$.

Choose $a \in A$.
For any $b,c \in A$, we have that
\begin{eqnarray*}
& &\hspace{-.6cm} \vfi(c^* (\io \od \phi)(\de(a)) \, b)
=   \phi(\,(\vfi \od \io)((c^* \ot 1)\de(a)(b \ot 1))\,) \\
& & \hspace{-1cm}\spat = \vfi(\al(\pi(\,(\vfi \od \io)((c^* \ot 1)\de(a)(b 
\ot 1))\,)\,)\,) 
= \vfi(\al(\,(\om_{\la(b),\la(c)} \ot \io)(\de(\pi(a))\,)\,) \\
& & \hspace{-1cm}\spat = \vfi(\,(\om_{\la(b),\la(c)} \be^{-1} \ot 
\io)\bigl((\be \ot
\al)\de(\pi(a))\bigr)\,) 
=   \vfi(\,(\om_{\la(b),\la(c)} \be^{-1} \ot \io)\de(\al(\pi(a)))\,) \\
& & \hspace{-1cm}\spat = (\om_{\la(b),\la(c)} \be^{-1})(1)\, 
\vfi(\al(\pi(a))) 
= \om_{\la(b),\la(c)}(1) \, \phi(a) = \vfi(c^* b) \, \phi(a)
\end{eqnarray*}
where we used  the left invariance of the weight $\vfi$ in the third last
equality.
The faithfulness of $\vfi$ implies that $(\io \od \phi)\de(a) = \phi(a) 1$.

By unicity of left invariant functionals on the $^*$-algebra level, we get 
the existence of a strictly positive
number $r$ such that $\phi = r \vfi$. So, $\vfi(\al(x)) = r \vfi(x)$ for every
$x \in \pi(A)$.

So, we have for any $y \in \pi(A)$, that $\al(y)$ belongs to
$\Nfi$ and $\| \lafi(\al(y)) \| = \sqrt{r} \| \lafi(y) \|$.

Remembering that $\pi(A)$ is a core for $\lafi$, we get that
$\al(y)$ belongs to $\Nfi$ and $\| \lafi(\al(y)) \| = \sqrt{r} \| \lafi(y)
\|$ for every $y \in \Nfi$.

\medskip

Using the fact that $(\be^{-1} \ot \al^{-1})\de = \de \al^{-1}$, we
infer completely analogously that $\al^{-1}(\Nfi) \subseteq \Nfi$,
so $\Nfi \subseteq \al(\Nfi)$.

Combining these two results, we conclude that $\al(\Nfi) = \Nfi$
and $\| \lafi(\al(y)) \| = \sqrt{r} \| \lafi(y) \|$ for all $y \in \Nfi$.
The result follows.
\end{demo}

\begin{proposition}
Consider a $^*$-automorphism $\al$ on $A_r$ such that there exists
a $^*$-automorphism $\be$ on $A_r$ such that $\de \al = (\be \ot \al)\de$.
Then $\si_t \al = \al \si_t$ and $\tau_t \be = \be \tau_t$ for every
$t \in \R$.
\end{proposition}
\begin{demo}
By the previous proposition, we know that there exists a strictly positive
number $r$ such that $\vfi \al = r \vfi$. Because $\vfi$ is KMS with respect 
to
$\si$, we have that $\si_t \al = \al \si_t$ for every $t \in \R$.

Choose $s \in \R$.
Then
\begin{eqnarray*}
& & (\be \tau_s \ot \si_s)\de = (\be \ot \io)\de\si_s 
= (\io \ot \al^{-1})\de\al\si_s \\
& & \spat = (\io \ot \al^{-1})\de\si_s\al 
= (\io \ot \al^{-1})(\tau_s \ot \si_s)\de\al \\
& & \spat =(\tau_s \ot \si_s)(\io \ot \al^{-1})(\be \ot \al)\de 
= (\tau_s \be \ot \si_s)\de.
\end{eqnarray*}
Therefore, $(\be \tau_s \ot \io)\de = (\tau_s \be \ot \io)\de$, which,
as before, implies that $\be \tau_s = \tau_s \be$.
\end{demo}

\begin{corollary} \label{prop6.4}
There exist a unique strictly positive number $\nu$
such that $\vfi \tau_t = \nu^t \vfi$ for all $t \in \R$.
For all $s,t \in \R$, we have that $\tau_s \si_t = \si_t \tau_s$.
\end{corollary}

We see that $\vfi$ is relatively invariant with respect to the 
one-parameter group $\tau$. In the definition of Masuda, Nakagami \&
Woronowicz, it is assumed that $\tau$ is invariant with respect to 
$\tau$. Untill now, we did not find any examples where we have a real
relative invariance, nor could we prove that $\vfi$ is invariant with
respect to $\tau$. Therefore, the question remains open whether
the definition of Masuda, Nakagami \& Woronowicz should be adapted in this 
respect.

\begin{proposition}
Consider a $^*$-automorphism $\al$ on $A_r$ such that there exists
a $^*$-automorphism $\be$ on $A_r$ such that $\de \al = (\be \ot \al)\de$.
Then $R \be = \be R$.
\end{proposition}
\begin{demo}
We know from the previous results that there exists a strictly positive number
$r$ such that $\vfi \al = r \, \vfi$. Moreover, we have for every $t \in \R$
that $\be \tau_t = \tau_t \be$, which implies that $\tau_{-\frac{i}{2}}
 \be =  \be  \tau_{-\frac{i}{2}} $.

Choose $a \in \Nfi$ and $\om \in A_r^*$. Take $\th \in A_r^*$.
For every $n \in \N$, we define $\th_n \in A_r^*$ such that
$$\th_n(x) = \frac{n}{\sqrt{\pi}} \int \exp(-n^2 t^2) \, \th(\tau_t(x)) \, dt 
$$
for every $x \in A_r$.
It is clear that $(\th_n)_{n=1}^\infty$ converges weakly$^*$ to $\th$.

\medskip

Fix $m \in \N$ and define $\eta \in A_r^*$ such that
$$\eta(x) = \frac{m}{\sqrt{\pi}} \int \exp(-m^2 (t-\frac{i}{2})^2) \,
\th(\tau_t(x)) \, dt $$
for every $x \in A_r$.
It is not so difficult to see that $\eta(\tau_{-\frac{i}{2}}(x)) = \th_m(x)$
for every $x \in {\cal D}(\tau_{-\frac{i}{2}})$. This implies that
$\eta \tau_{-\frac{i}{2}}$ is bounded and $\overline{\eta
\tau_{-\frac{i}{2}}} = \th_m$ \ \ (a).

Choose $b \in \Nfi$. Then we see that
\begin{eqnarray*}
& & \vfi(\,(\th_m R \be \ot \io)(\de(a^*)) \, b )
= \vfi(\al^{-1}(\,(\th_m R \ot \io)(\de(\al(a)^*))\,) \, b) \\
& & \spat = \vfi(\al^{-1}(\,(\th_m R \ot \io)(\de(\al(a)^*)) \, \al(b)\,))
= r^{-1} \, \vfi(\,(\th_m R \ot \io)(\de(\al(a)^*)) \, \al(b)\,).
\end{eqnarray*}
We have that $\eta \, R \, \tau_{-\frac{i}{2}} = \eta \, \tau_{-\frac{i}{2}}
\, R$. Together with (a), this implies that $\eta \, R \tau_{-\frac{i}{2}}$ is
bounded and $\overline{\eta R \tau_{-\frac{i}{2}}} = \th_m R$.
Hence, using the foregoing chain of equalities and proposition \ref{prop5.4}
we get that
\begin{eqnarray*}
& & \vfi(\,(\th_m R \be \ot \io)(\de(a^*)) \, b )
= r^{-1} \, \vfi(\al(a)^*  (\eta \ot \io)\de(\al(b))\,) \\
& & \spat = r^{-1} \, \vfi(\al(a)^*  \al((\eta \be \ot \io)\de(b))\,)
= r^{-1} \, \vfi(\al(a^* (\eta \be \ot \io)\de(b)\,)\,) \\
& & \spat = \vfi(a^* (\eta \be \ot \io)\de(b)\,).
\end{eqnarray*}
Furthermore, $\eta \be R \tau_{-\frac{i}{2}} = \eta \tau_{-\frac{i}{2}} \be 
R$.
So, by using (a) once more, we see that $\eta \be R \tau_{-\frac{i}{2}}$
is bounded and $\overline{\eta \be R \tau_{-\frac{i}{2}}}
= \th_m \be R$. Refering to proposition \ref{prop5.4} once again, we get that
$$ \vfi(a^* (\eta \be \ot \io)\de(b)\,)
= \vfi(\,(\th_m \be R \ot \io)(\de(a^*)) \, b) \ ,$$
which implies that
$$ \vfi(\,(\th_m R \be \ot \io)(\de(a^*)) \, b )
= \vfi(\,(\th_m \be R \ot \io)(\de(a^*)) \, b) \ .$$

Because $\vfi$ is faithful, we can conclude from  this that
$$(\th_m R \be \ot \io)(\de(a^*)) = (\th_m \be R \ot \io)(\de(a^*)).$$
So, applying $\om$ to this equation, gives
$$ \th_m \bigl(\,(\be R)((\io \ot \om)\de(a^*))\,\bigr)
=  \th_m \bigl(\,(R \be)((\io \ot \om)\de(a^*))\,\bigr).$$

\medskip

Because $(\th_n)_{n=1}^\infty$ converges weakly$^*$ to $\th$, we get that
$$ \th \bigl(\,(\be R)((\io \ot \om)\de(a^*))\,\bigr)
=  \th \bigl(\,(R \be)((\io \ot \om)\de(a^*))\,\bigr).$$

Consequently,
$(\be R)((\io \ot \om)\de(a^*))
=  (R \be)((\io \ot \om)\de(a^*))$.
The density conditions imply that $\be R = R \be$.
\end{demo}

\medskip

We would like to mention the following application of these results:

Consider a locally compact group $G$ together with a norm continuous one 
parameter group $\al$
on $A_r$ such that $(\al_t \ot \al_t)\de = \de \al_t$ for every $t \in \R$. 
Then we get from
the preceding propositions that the weight $\vfi$ is relatively invariant 
under $\al$ and that
$\al$ commutes with $\si$,$\tau$ and $R$. These kind of properties allow us 
to make
a left Haar weight with corresponding modular group, a scaling group and a 
anti-unitary 
antipode on the level of the crossed product $G \times_\al A_r$, turning it 
into a 
C$^*$-algebraic quantum group according to Masuda, Nakagami \& Woronowicz.

\medskip

Now, we turn our attention to group-like elements.

\begin{proposition}      \label{prop6.3}
Consider a unitary element $u \in M(A_r)$ such that there exists an element $v
\in M(A_r)$ such that $\de(u)=v \ot u$.
Then there exist a unique strictly positive number $\lambda$ such
that $\si_t(u) = \lambda^{it} u$ for every $t \in \R$.
Moreover, $\tau_t(v) = v$ for all $t \in \R$.
\end{proposition}
\begin{demo}
We first proof the assertion about $\si$.
\begin{enumerate}
\item Define the $^*$-automorphisms $\al,\be$ on $A_r$
      such that $\al(a) = u^* a u$ and $\be(a)=v^* a v$ for all $a \in A$.
      The formula $\de(u) = v \ot u$ implies easily that $(\be \ot \al)\de
      = \de \al$. By proposition \ref{prop6.2} there exists a strictly 
positive
      number $\mu$ such that $\vfi \al = \mu \vfi$.
      Consequently, we have for every $x \in \Mfi$ that $u^* x u$ belongs
      to $\Mfi$ and $\vfi(u^* x u) = \mu \vfi(x)$.

      Hence, we have for all $y \in \Nfi$ that $(y u)^* (y u) = u^* y^* y u 
\in
      \Mfi^+$ so that $y u$ belongs to $\Nfi$.  Thus, $\Nfi u \subseteq \Nfi$.
\item Analogously, working with $u^*$ in stead of $u$, one proves
      that $\Nfi u^* \subseteq \Nfi$, so $u \Nfi^* \subseteq \Nfi^*$.
\end{enumerate}
From 1) and 2), it follows immediately that $u \Mfi \subseteq \Mfi$
and $\Mfi u \subseteq \Mfi$.

Choose $x,y \in \Nfi$.
From 2), it follows that $u y^*$ belongs to $\Nfi^*$, so $u y^* x$ belongs
to $\Mfi$. Using 1), we get that
$$ \vfi(y^* x u) = \vfi(u^* (u y^* x) u) = \mu \vfi(u y^* x) .$$
Using lemma \ref{lemA1}, we
conclude that $\si_t(u) = \mu^{-it} u $ for all $t \in \R$.
Putting $\lambda = \frac{1}{\mu}$, we see that $\si_t(u) = \lambda^{it} u$
for all $t \in \R$.

Now, we will prove the assertion about $\tau$.
Choose $t \in \R$. Then
\begin{eqnarray*}
& & \tau_t(v) \ot \si_t(u) = (\tau_t \ot \si_t)\de(u) = \de(\si_t(u)) \\
& & \spat = \lambda^{it} \de(u) = \lambda^{it} \, v \ot u =  v \ot \si_t(u).
\end{eqnarray*}
Therefore, $\tau_t(v) = v$.
\end{demo}

In the following proposition, we prove that every one dimensional unitary
corepresentation of $A_r$ is of an algebraic nature. We first prove a little
lemma.

\begin{lemma}                \label{lem6.1}
Consider $a,b,c \in A$ and $x \in M(A_r)$, then
$(\io \ot \om_{\la(b),\la(c)})(\de(x)) \, \pi(a)$ belongs to $\pi(A)$.
\end{lemma}
\begin{demo}
We know that there exist $p_1,\ldots\!,p_n \, , \, q_1,\ldots\!,q_n \, , \,
r_1,\ldots\!,r_n \in A$ such that
$$a \ot b \rho(c^*) = \sum_{i=1}^n \de(p_i \rho(q_i^*)) (r_i \ot 1).$$
For every $y \in A$, we have that
\begin{eqnarray*}
& & (\io \ot \om_{\la(b),\la(c)})(\de(\pi(y))(\pi(a) \ot 1)\,)
= (\io \ot \om_{\la(b),\la(c)})(\,(\pi \od \pi)(\de(y)(a \ot 1))\,) \\
& & \spat = \pi(\,(\io \od \vfi)((1 \ot c^*)\de(y)(a \ot b))\,)
= \pi(\,(\io \od \vfi)(\de(y)(a \ot b \rho(c^*))\,) \\
& & \spat = \sum_{i=1}^n \pi(\,(\io \od \vfi)(\de(y p_i \rho(q_i^*))(r_i \ot 
1))\,)
= \sum_{i=1}^n \pi(\vfi( y p_i \rho(q_i^*)) \, r_i) \\
& & \spat = \sum_{i=1}^n  \vfi(q_i^* y p_i) \, \pi(r_i)
= \sum_{i=1}^n \om_{\la(p_i),\la(q_i)}(\pi(y)) \pi(r_i).
\end{eqnarray*}
Because $\pi(A)$ is dense in $A_r$ and because of strict continuity arguments
we can replace $\pi(y)$ in this equation by any element in $M(A_r)$. Hence,
$$ (\io \ot \om_{\la(b),\la(c)})(\de(x)(\pi(a) \ot 1)\,)
= \sum_{i=1}^n \om_{\la(p_i),\la(q_i)}(x) \, \pi(r_i).$$
Consequently,
$(\io \ot \om_{\la(b),\la(c)})(\de(x)) \pi(a)$ belongs to $\pi(A)$.
\end{demo}

\begin{proposition}  \label{prop6.1}
Consider a unitary $u$ in $M(A_r)$ such that $\de(u) = u \ot u$.
Then there exists a unique $x \in M(A)$ such that
$u \, \pi(a) = \pi(x a)$ and $\pi(a) \, u = \pi(a  x)$ for every $a \in A$.
Moreover, $\de(x) = x \ot x$.
\end{proposition}
\begin{demo}
Choose $a \in A$.
Now, there exist $b,c \in A$ such that $ \langle u \, \la(b) , \la(c) \rangle 
 = 1$
Therefore,
\begin{eqnarray*}
& & u\, \pi(a) = u \, \pi(a) \, \om_{\la(b),\la(c)}(u)
= (\io \ot \om_{\la(b),\la(c)})(u \ot u)\, \pi(a) \\
& & \spat\spat = (\io \ot \om_{\la(b),\la(c)})(\de(u)) \, \pi(a).
\end{eqnarray*}
which, by the previous lemma, implies that $u \, \pi(a)$ belongs to $\pi(A)$.

Using the fact that $\de(u^*) = u^* \ot u^*$, we conclude also
that $u^* \pi(A) \subseteq \pi(A)$, hence $\pi(A) \, u \subseteq \pi(A)$.

Because $u \, \pi(A)$ and $\pi(A) \, u$ are subsets of $\pi(A)$, it follows
easily that there exist $x \in M(A)$ such that
$u \, \pi(a) = \pi(x a)$ and $\pi(a) \, u = \pi(a  x)$ for every $a \in A$.
It is also not so difficult to prove that $x$ is unitary and that
$\de(x) = x \ot x$.
\end{demo}

\begin{proposition}   \label{prop6.5}
Consider a unitary element $u \in M(A_r)$ such that  $\de(u) = u \ot u$.
Then $R(u) = u^*$.
\end{proposition}
\begin{demo}
The previous proposition ensures the existence of a uniqe unitary  element $x
\in M(A)$ such that $\pi(a)\,u = \pi(a x)$ and $u \, \pi(a) = \pi(x a)$ for 
all
$a \in A$. Moreover, $\de(x) = x \ot x$. This last equation implies that
$S(x) = x^*$.

Choose $a \in A$.
For any $t \in \R$, we have that
$$ \tau_t(\pi(S(a) x^*)) = \tau_t(\pi(S(a))\, u^* )
= \tau_t(\pi(S(a))\,)\, u^*  ,$$
where we used proposition  \ref{prop6.3}.
This implies that
$$\tau_\frac{i}{2}(\pi(S(a) x^*)) = \tau_\frac{i}{2}(\pi(S(a))\,) \, u^*
= R(\pi(a)) \, u^*.$$
Furthermore,
\begin{eqnarray*}
& & R(\pi(a))\,R(u) = R(u \pi(a)) = R(\pi(xa)) =
\tau_\frac{i}{2}(\pi(S(xa))\,) \\
& & \spat = \tau_\frac{i}{2}(\pi(S(a)S(x))\,) = \tau_\frac{i}{2}(\pi(S(a) 
x^*))
=  R(\pi(a)) \, u^*.
\end{eqnarray*}
So, we can conclude that $R(u) = u^*$.
\end{demo}

Now, we want to use these results about unitary elements to prove
results about strictly positive elements affiliated with the \cst \, $A_r$. 
They
will apply to the modular function of our quantum group.

\begin{proposition}
Consider a strictly positive element $\al \,\et \,  A_r$ such that there 
exists a
strictly positive element $\be \,\et \,  A_r$ such that $\de(\al) = \be \ot 
\al$.
Then there exists a unique strict positive number $\lambda$
such that $\si_t(\al) = \lambda^t \, \al$ for every $t \in \R$.
Moreover, $\tau_t(\be) = \be$ for every $t \in \R$.
\end{proposition}
\begin{demo}
First, we will prove the assertion about $\si$.

Choose $s \in \R$. Then $\de(\al^{si}) = \be^{si} \ot \al^{si}$.
By proposition \ref{prop6.3}, there exists a unique strictly positive
number $\lambda_s$ such that $\si_t(\al^{si}) = (\lambda_s)^{it} \al^{si}$
for every $t \in \R$.

\begin{list}{}{\setlength{\leftmargin}{0.4cm}}

\item Choose $s_1,s_2 \in \R$.

For every $t \in \R$, we have that
\begin{eqnarray*}
& & (\lambda_{s_1 + s_2})^{it} \al^{(s_1 + s_2)i}
= \si_t(\al^{(s_1+s_2)i})
= \si_t(\al^{s_1 i}) \si_t(\al^{s_2 i}) \\
& & \spat = (\,(\lambda_{s_1})^{it} \al^{s_1 i})(\,(\lambda_{s_2})^{it} 
\al^{s_2 i})
= (\lambda_{s_1} \lambda_{s_2})^{it} \al^{(s_1 + s_2)i} \ ,
\end{eqnarray*}
hence, $(\lambda_{s_1 + s_2})^{it} = (\lambda_{s_1} \lambda_{s_2})^{it}$.

This implies that $\lambda_{s_1 + s_2} = \lambda_{s_1} \lambda_{s_2}$.
\end{list}

We put $\lambda = \lambda_{1} \in \R_0^+$.

From the previous result, it follows easily that $\lambda_q = \lambda^q$ for
every $q \in \Q$.

Fix $t \in \R$.

For any $s \in \R$, we have that $(\lambda_s)^{it}  1 = \si_t(\al^{si})
\al^{-si}$. Therefore, because the mapping $\R \rightarrow M(A_r):s \mapsto
\si_t(\al^{si}) \al^{-si}$ is strictly continuous, the mapping
$\R \rightarrow \C : s \mapsto (\lambda_s)^{it}$ is continuous.

We know also that $(\lambda_q)^{it} = \lambda^{qit}$ for $q \in \Q$,
thus $(\lambda_s)^{it} =\lambda^{ist}$ for every $s \in \R$.

Hence,
$$
\si_t(\al)^{si} = \si_t(\al^{si}) = (\lambda_s)^{it} \al^{si}
= \lambda^{ist} \al^{is} =  (\lambda^t \al)^{is}
$$
for all $s \in \R$. Consequently, $\si_t(\al) = \lambda^t \al$.

Next, we quickly prove the assertion concerning $\tau$.

Choose $t \in \R$

We have for any $s \in \R$ that $\de(\al^{si}) = \be^{si} \ot \al^{si}$,
so proposition \ref{prop6.3} implies that $\tau_t(\be^{si}) = \be^{si}$,
hence $\tau_t(\be)^{si}= \be^{si}$.
Therefore, we see that $\tau_t(\be) = \be$.
\end{demo}

\begin{proposition}
Consider a stictly positive element $\al \,\et \,  A_r$
such that $\de(\al) = \al \ot \al$. Then $R(\al) = \al^{-1}$.
\end{proposition}
\begin{demo}
Choose $t \in R$. Then $\de(\al^{it}) = \al^{it} \ot \al^{it}$.
Proposition \ref{prop6.5} implies that $R(\al^{it}) = \al^{-it}$,
hence $R(\al)^{it} = (\al^{-1})^{it}$.
This implies that $R(\al) = \al^{-1}$.
\end{demo}

\medskip

\begin{lemma} \label{lem6.2}
Consider $x \in M(A_r)$ such that $\de(x) = x \ot 1$. Then
$x$ belongs to $\C 1$.
\end{lemma}
\begin{demo}
By proposition \ref{prop2.1}, we have that
$$\de(\si_t(x)) = (\si_t \ot K_t)\de(x) = (\si_t \ot K_t)(x \ot 1)
= \si_t(x) \ot 1   \text{\ \ \ (a)}$$
for every $t \in \R$.
Fix $n \in \N$ and define $x_n \in M(A_r)$ such that
$$x_n a = \frac{n}{\sqrt{\pi}} \int \exp(-n^2 t^2) \, \si_t(x) a \, dt $$
for all $a \in A_r$. Then, equation (a) implies that $\de(x_n) = x_n \ot 1$.

Define the function $F$ from $\C$ to $M(A_r)$ such that
$$F(z) a = \frac{n}{\sqrt{\pi}} \int \exp(-n^2(t-z)^2) \, \si_t(x) a \, dt$$
for all $z \in \C$ and $a \in A_r$.
It is clear that $F$ is strictly analytic and hence analytic on
$\C$. Moreover, we have that $F(u) = \si_u(x_n)$ for every $u \in \R$.

Using a standard technique, it is not difficult to prove, using lemma
\ref{prop5.2}, that $\Nfi \, x_n \subseteq \Nfi$ and that
$\lafi(a x_n) = J F(\frac{i}{2})^* J \lafi(a)$ for every
$a \in \Nfi$. So, we get that $\Mfi \, x_n \subseteq \Mfi$.

Now, we are in a position to use a technique which we borrowed from 
Woronowicz.

There exists an element $d \in \Mfi$ such that $\vfi(d)=1$. By the previous
result, we know that $d x_n$ belongs also to $\Mfi$.
Moreover, $\de(d)(x_n \ot 1) = \de(d x_n)$.

Choose $\om \in A_r^*$ and apply $\om \ot \io$ to the previous equation,
this gives $(x_n \om \ot \io)\de(d) = (\om \ot \io)\de(d x_n) $ \ \ (b).

The left invariance of $\vfi$ implies that $(x_n \om \ot \io)\de(d)$ belongs
to $\Mfi$ and
$$ \vfi(\,(x_n \om \ot \io)\de(d) \,) = (x_n \om)(1) \, \vfi(d) = \om(x_n) .$$
At the same time, we get that $(\om \ot \io)\de(d x_n)$ belongs to $\Mfi$
and
$$\vfi(\,(\om \ot \io)\de(d x_n)\,) = \om(1) \, \vfi(d x_n).$$
Combining these two equalities, using equation (b), gives us that
$\om(x_n) = \om( \vfi(d x_n) 1)$.

Therefore, we have that $x_n = \vfi(d x_n) 1$.

It is clear that $(x_n)_{n=1}^\infty$ converges strictly to
$x$. Using the Cauchy criterium, this implies easily that there exists a 
complex
number $c$ such that $(\vfi(d x_n)\,)_{n=1}^\infty \rightarrow c$.
Then, $x$ must be equal to $c 1$.
\end{demo}

Using the fact that $\flip(R \ot R)\de = \de R$, it is straightforward to 
arrive at the following conclusion.

\begin{lemma}
Consider $x \in M(A_r)$ such that $\de(x) = 1 \ot x$. Then
$x$ belongs to $\C 1$.
\end{lemma}

\medskip

It is clear that all the things we prove about $A_r$, we can also prove about
$\ah_r$. For instance, we can extend the right Haar functional $\psih$ on 
$\ah$
to a weight on $\ah_r$ in the same way as we did for $\vfi$. We will go a 
little
bit further into this because this gives a nice implementation of the polar 
decomposition
of the antipode.

\medskip

As can be expected, we will introduce again a left Hilbert algebra:

\begin{definition}
We define $\cuh = \lah(\ah)$, so $\cuh$ is a dense subspace of $H$.
We will make $\cuh$ into a $^*$-algebra in the following way.
\begin{enumerate}
\item For all $a,b \in \ah$, we have that $\lah(a) \lah(b) = \lah(ab)$.
\item For all $a \in A$, we have that $\lah(a)^\circ = \lah(a^*)$.
\end{enumerate}
Just as before, we have that $\cuh$ is a left Hilbert algebra on $H$.
\end{definition}

Again, we will introduce some terminology in connection with this left 
Hilbert algebra.

We will have that $\hat{L}_{\lah(a)} = \pih(a)$ for all $a \in \ah$.

Define the mapping $\Th$ as the closed antilinear map from within H into H
such that $\cuh$ is a core for $\Th$ and $\Th v = v^\circ $ for every $v \in 
\cu$.
Furthermore, we define $\nabh = \Th^* \Th$.

We also define $\Jh$ as the anti-unitary operator arising from the polar
decomposition of $\Th$. Hence, $\Th = \Jh \nabh^\frac{1}{2}$.

\medskip

We are going to interprete these map $\Th$ in another way which will have 
some implications
on the level of the antipode.

Remember that, by definition, $\la(a)=\lah(\hat{a})$ for every $a \in A$. So, 
$\cuh$
equals $\la(A)$ as a subspace of $H$. Therefore $\la(A)$ is a core for $T$.

\begin{lemma}
We have that $\hat{a}\,^* = (S(a)^* \sde)\hoed$ for every $a \in A$.
\end{lemma}
\begin{demo}
Choose $x \in A$. Then
$$ \hat{a}\,^*(x) = \overline{\hat{a}(S(x)^*)} = \overline{\vfi(S(x)^* a)}
= \vfi(a^* S(x)) = \vfi(S(x S(a)^*)\,) \ .$$
Because $\vfi S = \sde \vfi$, we get that 
$$ \hat{a}\,^*(x) = \vfi(x S(a)^* \sde) = (S(a)^* \sde)\hoed\,(x) \ .$$
\end{demo}

This lemma implies easily that $\Th \la(a) = \la(S(a)^* \sde)$ for every $a 
\in A$.
A short investigation learns that we are in a situation where the construction
of section 4 applies.  Consider the construction procedure of section 4 with 
$\eta=\vfi$,
$x=\frac{1}{\mu} \, \sde$ and $y=\sde$. Moreover, take $K=H$ and $\ga=\la$.

If we go through this construction procedure, we see that $G=T$. This
implies that $\cp=\nabh$ and $\cj=\Jh$.

In this case $V=W$, so proposition \ref{prop3.2} of section 4 implies
the following result.

\begin{proposition}   
We have that $(\Jh \ot J)W = W^* (\Jh \ot J)$
and $ (\nabh \ot \nab)  W = W (\nabh \ot \nab)$.
\end{proposition}

A similar result can be found in \cite{MasNak}. Using this proposition, we 
can prove easily the
following one.

\begin{proposition} \label{prop6.6}
\begin{enumerate}
\item For every $t \in \R$ and $x \in A_r$, we have that $\tau_t(x) = 
\nabh^{it} \, x \, \nabh^{-it}$.
\item For every $x \in A_r$, we have that $R(x)= \Jh \, x^* \, \Jh$.
\end{enumerate}
\end{proposition}
\begin{demo}
Fix $t \in \R$. Choose $\om \in \cop^*$. We know by the above commutation 
relation that
$$(\nabh^{it} \ot 1)\,W\,(\nabh^{-it} \ot 1) = (1 \ot \nab^{-it}) \, W \, (1 
\ot \nab^{it}) \ .$$
Aplying $\io \ot \om$ to this equation, gives us that
$$ \nabh^{it} \, (\io \ot \om)(W) \, \nabh^{-it}
= (\io \ot \nab^{it} \om \nab^{-it})(W) \ .$$
Therefore, equation \ref{eq4.2} implies that 
$$\nabh^{it} \, (\io \ot \om)(W) \, \nabh^{-it} = \tau_t(\,(\io \ot 
\om)(W)\,) \ . $$
Hence, it follows that $\tau_t(x) = \nabh^{it} \, x \, \nabh^{-it}$ for every 
$x \in A_r$.

The assertion about $R$ is proven in a similar way.
\end{demo}

This implementation of $R$, together with lemma \ref{lem6.2} allows us to 
prove another
result which can be found in \cite{MasNak}.

\begin{proposition}
We have that $M(A_r) \cap M(\ah_r) = \C 1$.
\end{proposition}
\begin{demo}
The following argument is due to \cite{MasNak}, proposition 3.11. 

Choose $x \in M(A_r) \cap M(\ah_r)$. Because $x$ belongs to $M(\ah_r)$, we 
have that
$x$ belongs to ${\cal L}(\cuh)$. By the theory of left Hilbert algebras, we 
know
that $\Jh x^* \Jh$ belongs to ${\cal L}(\cuh)'$. Therefore $R(x)$ belongs to
${\cal L}(\cuh)'$, so we get that $R(x)$ commutes with $\ah_r$.

Remembering that $W$ belongs to $M(A_r \ot \ah_r)$, this implies that
$W^*(1 \ot R(x))W = 1 \ot R(x) \ .$

We know also that $R(x)$ belongs to $M(A_r)$ and 
$$\de(R(x)) = W^*(1 \ot R(x))W = 1 \ot R(x) \ .$$
Lemma \ref{lem6.2} implies that $R(x)$ belongs to $\C 1$, so $x$ belongs
to $\C 1$.
\end{demo}

\section{The modular function.}

We already have a modular function $\sde$ on the $^*$-algebra level. In this 
section,
we are going to introduce the modular function of our $C^*$-algebraic quantum
group. For the notations used in this part, we refer to section \ref{sec2}.

\medskip

It is easy to see that
$ \langle \lade(a),\lade(b) \rangle  =  \langle \la(a),\la(\sde b) \rangle $ 
for all $a,b \in A$.
As before, this result justifies the following definition

\begin{definition}
We define the closed linear operator $L$ from within $H$ into $\hde$
such that $\la(A)$ is a core for $L$ and $L \la(a) = \lade(a)$
for every $a \in A$.

Then $L$ is a densely defined injective operator with dense range.
\end{definition}

As before, it is easy to check that
$$ \langle L v , \lade(a) \rangle  =  \langle v , \la(\sde a)  \rangle $$
for every $a \in A$. It follows that $\lade(A)$ is a subset of $D(L^*)$
and $L^* \lade(a) = \la(\sde a)$ for all $a \in A$.

\medskip

Next, we will give the definition of the modular function of the our quantum
group $A_r$. Again, we will use the same symbol $\sde$ as for the modular 
function on the
$^*$-algebra level. By looking on which elements $\sde$ acts, it should be 
clear which $\sde$ is meant (elements in $A$ vs.\ elements in $A_r$).

\begin{definition}
We define $\sde = L^*L$, so $\sde$ is an injective positive operator in $H$.
\end{definition}

It is clear that $\la(A)$ is a subset of $D(\sde)$ and that
$\sde \la(a) = \la(\sde a)$.
Later, we will prove that $\la(A)$ is a core for $\sde$.
In the next part, we want to prove that $\sde$ is affiliated with $A_r$.

\begin{lemma}
Consider $a,b \in A$, then
\begin{eqnarray*}
& &  \langle  (\lade \od \lade)(\de(b)(a \ot 1)\,) , (\lade \od 
\lade)(\de(d)(c \ot
1)\,)  \rangle  \\
& & \spat =  \langle \la(a) \ot \lade(b) , \la(c) \ot \lade(d) \rangle .
\end{eqnarray*}
\end{lemma}
\begin{demo}
We have that
\begin{eqnarray*}
& &  \langle  (\lade \od \lade)(\de(b)(a \ot 1)\,) , (\lade \od 
\lade)(\de(d)(c \ot
1)\,)  \rangle  \\
& & \spat = (\vfi \od \vfi)(\,(c^* \ot 1) \de(d^*) (\sde \ot \sde) \de(b)(a 
\ot
1)\,) \\
& & \spat = (\vfi \od \vfi)(\,(c^* \ot 1) \de(d^*) \de(\sde) \de(b)(a \ot 
1)\,) \\
& & \spat = (\vfi \od \vfi)(\,(c^* \ot 1) \de(d^* \sde b)(a \ot 1)\,).
\end{eqnarray*}
So, by the left invariance of $\vfi$, we find that
\begin{eqnarray*}
& &  \langle  (\lade \od \lade)(\de(b)(a \ot 1)\,) , (\lade \od 
\lade)(\de(d)(c \ot
1)\,)  \rangle  \\
& & \spat = \vfi(c^* a) \, \vfi(d^* \sde b)
=   \langle \la(a) \ot \lade(b) , \la(c) \ot \lade(d) \rangle .
\end{eqnarray*}
\end{demo}

\begin{lemma}
We have that $(1 \ot \sde)W = W(\sde \ot \sde)$.
\end{lemma}
\begin{demo}
By the previous lemma, there exists a unitary operator $U$ from  $H \ot \hde$
to $\hde \ot \hde$ such that
$$U (\la(a) \ot \lade(b)) = (\lade \ot \lade)(\de(b)(a \ot 1)\,)$$
for every $a,b \in A$. We easily infer that
\begin{eqnarray*}
& & (L \ot L)W^* \, (\la(a) \ot \la(b))
= (L \ot L) (\la \od \la)(\de(b)(a \ot 1)\,) \\
& & \spat = (\lade \od \lade)(\de(b)(a \ot 1)\,)
= U (\la(a) \ot \lade(b)) \\
& & \spat = U (1 \ot L) \, (\la(a) \ot \la(b)) \\
\end{eqnarray*}
for every $a,b \in A$.
Using, lemma \ref{lem2.2} once more, we get that
$ (L \ot L) W = W (1 \ot L) $.
As before, this implies that $W(\sde \ot \sde) = (1 \ot \sde)W$.
\end{demo}

We are now in a position to use a technique from Woronowicz.

\begin{proposition}
We have that $\sde$ is a strictly positive element affiliated with $A_r$
(in the \cst\ sense).
\end{proposition}
\begin{demo}
Choose $t \in \R$ and $\om \in \cop$.
We know that $(1 \ot \sde^{it}) W (1 \ot \sde^{-it}) = W (\sde^{it} \ot 1)$.
Applying $\io \ot \om$ to this equation gives that $(\io \ot \om)(W) \,
\sde^{it} =  (\io \ot \sde^{-it} \om \sde^{it})(W)$.

Consequently, we have for every $\om \in \cop$ that
\begin{enumerate}
\item For every $t \in \R$,  $(\io \ot \om)(W)\,\sde^{it}$ belongs to $A_r$.
\item The mapping $\R \rightarrow A_r: t \mapsto (\io \ot \om)(W)\,\sde^{it}$
      is continuous.
\end{enumerate}
It follows easily that every $a \in A_r$ satisfies:
\begin{enumerate}
\item For every $t \in \R$,  $a\,\sde^{it}$ belongs to $A_r$.
\item The mapping $\R \rightarrow A_r: t \mapsto a\,\sde^{it}$
      is continuous.
\end{enumerate}
The conclusion follows from these two results.
\end{demo}

Sometimes, we will look at $\sde$ as an operator in $H$. In other cases,
we will look at $\sde$ as an element affiliated with $A_r$. It will be clear
from the context which viewpoint is under consideration.

\begin{proposition}
We have that $\de(\sde) = \sde \ot \sde$.
\end{proposition}
\begin{demo}
Choose $t \in \R$.
We know that $(1 \ot \sde^{it}) W = W (\sde^{it} \ot \sde^{it})$.
Hence,
$$
\de(\sde)^{it} = \de(\sde^{it}) = W^* (1 \ot \sde^{it}) W
\sde^{it} \ot \sde^{it} =(\sde \ot \sde)^{it}.
$$
Therefore, $\de(\sde) = \sde \ot \sde$.
\end{demo}

Because of this equation and the results of the previous section, we have
already some nice results about $\sde$.
\begin{enumerate}
\item For every $t \in \R$, we have that $\sde^{it} \, \pi(A) \subseteq
\pi(A)$ and $\pi(A) \, \sde^{it} \subseteq \pi(A)$, so $\sde^{it}$
      belongs to the algebraic multiplier algebra of $\pi(A)$.
\item There exist a unique strictly positive number $\gamma$ such that
      $\si_t(\sde) = \gamma^t \, \sde$.  Later, we will prove that 
$\gamma=\nu^{-1}$.
\item For every $t \in \R$, we have that $\tau_t(\sde) = \sde$.
\item We have that $R(\sde) = \sde^{-1}$.
\end{enumerate}

\begin{lemma}
Consider $a \in A$ and $n \in \Z$. Then $\la(a)$ belongs to $D(\sde^n)$
and $\sde^n \la(a) = \la(\sde^n a)$.
\end{lemma}

This proposition is trivially true for $n=0$. We know already that it is true
for  $n=1$. This implies easily the result for $n=-1$. Induction will give us
the general result.

\begin{lemma}
For every $a \in A$ and $z \in \C$, we have that $\la(a)$ belongs to
$D(\sde^z)$.
\end{lemma}

In the next propositions, we will look at $\sde$ as an element affiliated with
$A_r$.

\begin{proposition}  \label{prop7.1}
Consider $a \in A$ and $n \in \Z$. Then $\pi(a)$ belongs to ${\cal D}(\sde^n)$
and $\sde^n \pi(a) = \pi(\sde^n a)$.
\end{proposition}
\begin{demo}
Choose $x \in A$. Then $\pi(a) \la(x) = \la(a x) \in D(\sde^n)$
and
$$\sde^n (\pi(a) \la(x)) = \la(\sde^n a x ) = \pi(\sde^n a) \la(x) .$$
Using the closedness of $\sde^n$, it is not so difficult to prove that
$\sde^n \pi(a) = \pi(\sde^n a)$ as operators on $H$.

Therefore, $\pi(a)$ belongs to ${\cal D}(\sde^n)$ and $\sde^n \pi(a) = 
\pi(\sde^n
a)$.
\end{demo}

\begin{corollary}
For every $z \in \C$ and $a \in A$, we have that $\pi(a)$ belongs to
${\cal D}(\sde^z)$.
\end{corollary}

\begin{lemma}
Consider an element $\al \, \et \, A_r$ and elements $a \in {\cal D}(\al)$, 
$b \in
A_r$. Then $\de(a)(b \ot 1)$ belongs to ${\cal D}(\de(\al))$
and $\de(\al)(\de(a)(b \ot 1)) = \de(\al(a))(b \ot 1)$.
\end{lemma}
\begin{demo}
Take an approximate unit $(e_k)_{k \in K}$ for $A_r$.
We know by the general theory of affiliated elements that we have for any $k 
\in
K$ that $\de(a)(b \ot e_k)$ belongs to ${\cal D}(\de(\al))$ and
$\de(\al)(\de(a)(b \ot e_k)) = \de(\al(a))(b \ot e_k)$.

Because $\de(a)(b \ot 1)$ and $\de(\al(a))(b \ot 1)$ belong to $A_r \ot A_r$,
we see that
$$(\de(a)(b \ot e_k)\,)_{k \in K} \rightarrow \de(a)(b \ot 1)  $$ and 
$$(\de(\al)(\de(a)(b \ot e_k))\,)_{k \in K} \rightarrow \de(\al(a))(b
\ot 1). $$
The closedness of $\de(\al)$ implies that
$\de(a)(b \ot 1)$ belongs to ${\cal D}(\de(\al))$
and  $\de(\al)(\de(a)(b \ot 1)) = \de(\al(a))(b \ot 1)$.
\end{demo}

We will now show that $\sde^z$ is of an algebraic nature.
The proof is a slight adaptation of the proof of proposition \ref{prop6.1}.

\begin{proposition}
Consider $z \in \C$, then $\sde^z \, \pi(A) \subseteq \pi(A)$.
\end{proposition}
\begin{demo}
Choose $a \in A$
Take also a non zero element $b \in A$, then $\pi(b)$ belongs to ${\cal 
D}(\sde^z)$
and $\sde^z \pi(b) \neq 0$ (because $\sde^z$ is injective).
Hence, there exist $c,d \in A$ such that $ \langle (\sde^z \pi(b)) \, \la(c) ,
\la(d) \rangle \, =1$. Moreover,
there exist $p_1,\ldots\!,p_n \, , \, q_1,\ldots\!,q_n \in A$ such that
$$a \ot b = \sum_{i=1}^n \de(p_i)(q_i \ot 1).$$
We know that $\de(\sde^z) = \sde^z \ot \sde^z$. So, by the previous lemma,
\begin{eqnarray*}
& & \sde^z \pi(a) \ot \sde^z \pi(b) = \de(\sde^z)(\pi(a) \ot \pi(b)) \\
& & \spat = \sum_{i=1}^n \de(\sde^z) (\de(\pi(p_i))(\pi(q_i) \ot 1)\,) \\
& & \spat = \sum_{i=1}^n \de(\sde^z  \pi(p_i))(\pi(q_i) \ot 1).
\end{eqnarray*}
Therefore, we get that
\begin{eqnarray*}
& & \sde^z \pi(a) = (\sde^z \pi(a)) \, \om_{\la(c),\la(d)}(\sde^z \pi(b)) \\
& & \spat = (\io \ot \om_{\la(c),\la(d)})(\sde^z \pi(a) \ot \sde^z \pi(b)) \\
& & \spat = \sum_{i=1}^n (\io \ot \om_{\la(c),\la(d)})(\de(\sde^z 
\pi(p_i))(\pi(q_i)
\ot 1)\,).
\end{eqnarray*}
By lemma \ref{lem6.1}, we have that $\sde^z \pi(a)$ belongs to $\pi(A)$.
\end{demo}

This proposition allows us to define any power of $\sde$ on the algebraic 
level.
Proposition \ref{prop7.1} guarantees that this definition is consistent with
the usual definition of integer powers of $\sde$.

\begin{definition}
Consider $z \in \C$, then there exist a unique $\sde^z$ in $M(A)$ such
that $\pi(\sde^z a) = \sde^z \pi(a)$ for every $a \in A$.
\end{definition}

It is not difficult to prove the usual rules for exponentiation on the
$^*$-algebra level: \begin{enumerate}
\item For any $z \in \C$, we have that $(\sde^z)^* = \sde^{\overline{z}}$.
\item For any $y,z \in \C$, we have that $\sde^y \, \sde^z = \sde^{y+z}$.
\item For any $t \in \R$, we have that $\sde^{it}$ is a unitary in $M(A)$.
\item For any $t \in \R$, we have that $\sde^t$ is positive in the sense
      that $\sde^t = \sde^\frac{t}{2} \, \sde^\frac{t}{2}$, where
      $\sde^\frac{t}{2}$ is self adjoint.
\end{enumerate}

In the next proposition, we show that $\sde$ on the $C^*$-algebra level is
determined by its values on $\pi(A)$.

\begin{proposition}
Consider $z \in \C$, then $\pi(A)$ is a core for $\sde^z$.
\end{proposition}
\begin{demo}
Define $\sde_z$ as the closure of the mapping $\pi(A) \rightarrow A_r :
x \mapsto \sde^z x$. Then $\sde^z$ an extension of $\sde_z$.

\medskip

Choose $n \in \N$ and $x \in \pi(A)$.  We define
$$x_n = \frac{n}{\sqrt{\pi}} \int \exp(-n^2 t^2) \, \sde^{it} x \, dt.$$
For every $t \in \R$, we have that $\exp(-n^2 t^2) \, \sde^{it} x$ belongs to
$\pi(A)$ (which is a  subset of $\sde_z$)
and
$$\sde_z(\exp(-n^2 t^2) \, \sde^{it} x) = \sde^z(\exp(-n^2 t^2) \, \sde^{it} 
x)
= \exp(-n^2 t^2) \, \sde^{it} (\sde^z x) .$$
Consequently, the function $\R \rightarrow A_r : t \mapsto
\sde_z(\exp(-n^2 t^2) \, \sde^{it} x) $ is integrable.
The closedness of $\sde_z$ implies that $x_n$ belongs to the domain of 
$\sde_z$.

\medskip

Because $\pi(A)$ is dense in $A_r$, we have that
$ \langle x_n \mid x \in \pi(A) \text{ and } n \in \N \rangle $ is a core for 
$\sde^z$.
It follows that $\sde_z = \sde^z$, so the proposition is proven.
\end{demo}

We know already that there exists a strictly positive number $\gamma$ such
that $\si_t(\sde) = \gamma^t \, \sde$. In proposition
\ref{prop6.4}, we introduced a strictly positive number $\nu$  such that $\vfi
\tau_t = \nu^t \vfi$ for every $t \in \R$.
We will prove in the next part that $\gamma = \nu^{-1}$.

\begin{lemma}
Consider $a \in \Nfi$ and $b \in {\cal D}(\sde^\frac{1}{2})$.
Then $\de(a)(1 \ot b)$ belongs to ${\cal N}_{\vfi \ot \io}$
and
$$ (\vfi \ot \io)(\,(1 \ot b^*)\de(a^* a)(1 \ot b)\,)
= \vfi(a^* a) \, (\sde^\frac{1}{2} b)^* (\sde^\frac{1}{2} b). $$
\end{lemma}
\begin{demo}
There exist sequences $(a_n)_{n=1}^\infty$ and $(b_n)_{n=1}^\infty$
in $A$ such that $(\pi(a_n)\,)_{n=1}^\infty \rightarrow a$, \
$(\pi(b_n)\,)_{n=1}^\infty \rightarrow b$,
$(\lafi(\pi(a_n))\,)_{n=1}^\infty \rightarrow \lafi(a)$
and $(\sde^\frac{1}{2} \pi(b_n)\,)_{n=1}^\infty \rightarrow \sde^\frac{1}{2} 
b$.

Choose $c_1,c_2,d_1,d_2 \in A$. Then
\begin{eqnarray*}
& &  \langle (\lafi \ot \io)(\de(\pi(c_1))(1 \ot \pi(d_1))\,) ,
(\lafi \ot \io)(\de(\pi(c_2))(1 \ot \pi(d_2))\,) \rangle      \\
& & \spat =  \langle (\lafi \ot \io)(\,(\pi \od \pi)(\de(c_1)(1 \ot d_1))\,) ,
(\lafi \ot \io)(\,(\pi \od \pi)(\de(c_2)(1 \ot d_2))\,) \rangle    \\
& & \spat = (\vfi \od \pi)(\,(1 \ot d_2^*) \de(c_2^* c_1) (1 \ot d_1)\,)
= \vfi(c_2^* c_1) \pi(d_2^* \sde d_1) \\
& & \spat =  \langle \lafi(\pi(c_1)) , \lafi(\pi(c_2)) \rangle  \, 
\pi(\sde^\frac{1}{2} d_2)^*
\pi(\sde^\frac{1}{2} d_1) .
\end{eqnarray*}
Therefore,
\begin{eqnarray*}
& &  \langle (\lafi \ot \io)(\de(\pi(c_1))(1 \ot \pi(d_1))\,) ,
(\lafi \ot \io)(\de(\pi(c_2))(1 \ot \pi(d_2))\,) \rangle      \\
& & \spat =  \langle \lafi(\pi(c_1)) , \lafi(\pi(c_2)) \rangle  \, 
(\sde^\frac{1}{2} \pi(d_2))^*
(\sde^\frac{1}{2} \pi(d_1)) .  \hspace{3.5cm} \text{a)} \\
\end{eqnarray*}

It is clear that $(\de(\pi(a_n))(1 \ot \pi(b_n))\,)_{n=1}^\infty \rightarrow
\de(a)(1 \ot b)$.
Using equation (a), we get for any $m,n \in \N$ that
\begin{eqnarray*}
& & \| (\lafi \ot \io)(\de(\pi(a_m))(1 \ot \pi(b_m))\,) -
(\lafi \ot \io)(\de(\pi(a_n))(1 \ot \pi(b_n))\,) \|^2 \\
& & \ \  = \|  \langle \lafi(\pi(a_m)) , \lafi(\pi(a_m)) \rangle  \, 
(\sde^\frac{1}{2}
\pi(b_m))^*(\sde^\frac{1}{2} \pi(b_m)) 
-   \langle \lafi(\pi(a_m)) , \lafi(\pi(a_n)) \rangle  \, (\sde^\frac{1}{2}
\pi(b_m))^*(\sde^\frac{1}{2} \pi(b_n)) \\
& & \ \ -   \langle \lafi(\pi(a_n)) , \lafi(\pi(a_m)) \rangle  \, 
(\sde^\frac{1}{2}
\pi(b_n))^*(\sde^\frac{1}{2} \pi(b_m)) 
+   \langle \lafi(\pi(a_n)) , \lafi(\pi(a_n)) \rangle  \, (\sde^\frac{1}{2}
\pi(b_n))^*(\sde^\frac{1}{2} \pi(b_n)) \|.
\end{eqnarray*}
This equation implies that
$\bigl(\,(\lafi \ot \io)(\de(\pi(a_n))(1 \ot 
\pi(b_n))\,)\,\bigr)_{n=1}^\infty$
is a Cauchy sequence and hence convergent in $H \ot A$.
The closedness of $\lafi \ot \io$ implies that
$\de(a)(1 \ot b)$ belongs to ${\cal N}_{\vfi \ot \io}$
and
$\bigl(\,(\lafi \ot \io)(\de(\pi(a_n))(1 \ot 
\pi(b_n))\,)\,\bigr)_{n=1}^\infty$
converges to $(\lafi \ot \io)(\de(a)(1 \ot b)\,)$.

By equation (a), we have that
\begin{eqnarray*}
& &  \langle (\lafi \ot \io)(\de(\pi(a_n))(1 \ot \pi(b_n))\,)   ,
(\lafi \ot \io)(\de(\pi(a_n))(1 \ot \pi(b_n))\,)   \rangle  \\
& & \spat =   \langle \lafi(\pi(a_n)) , \lafi(\pi(a_n)) \rangle  \, 
(\sde^\frac{1}{2}
\pi(b_n))^*(\sde^\frac{1}{2} \pi(b_n))
\end{eqnarray*}
for every $n \in \N$.
So,
$$ \bigl(\, \langle (\lafi \ot \io)(\de(\pi(a_n))(1 \ot \pi(b_n))\,)   ,
(\lafi \ot \io)(\de(\pi(a_n))(1 \ot \pi(b_n))\,)   \rangle   
\bigr)_{n=1}^\infty $$
converges to $ \langle \lafi(a),\lafi(a) \rangle  \, (\sde^\frac{1}{2} b)^* 
(\sde^\frac{1}{2}
b)$.  Hence,
\begin{eqnarray*}
& &  \langle (\lafi \ot \io)(\de(a)(1 \ot b)\,)   ,
(\lafi \ot \io)(\de(a)(1 \ot b)\,)   \rangle  \\
& & \spat =   \langle \lafi(a) , \lafi(a) \rangle  \, (\sde^\frac{1}{2} 
b)^*(\sde^\frac{1}{2} b).
\end{eqnarray*}
This implies that
$$ (\vfi \ot \io)(\,(1 \ot b^*)\de(a^* a)(1 \ot b)\,)
= \vfi(a^* a) \, (\sde^\frac{1}{2} b)^* (\sde^\frac{1}{2} b). $$
\end{demo}

The following lemma is an easy consequence of the previous one.

\begin{lemma}
Consider $a \in \Nfi$ and $b \in {\cal D}(\sde)$.
Then
$$ (\vfi \ot \io)(\,(1 \ot b^*)\de(a^* a)(1 \ot b)\,)
= \vfi(a^* a) \, (\sde b)^* b = \, \vfi(a^* a) \, b^* (\sde b).$$
\end{lemma}

\begin{proposition}
We have that $\si_t(\sde) = \nu^{-t} \, \sde$ for every $t \in \R$.
\end{proposition}
\begin{demo}
Remember from the beginning of this section the existence of a strictly 
positive number
$\gamma$ such that $\si_t(\sde) = \gamma^t \, \sde$ for every $t \in \R$.

We take an element $a \in \Nfi$ such that $\vfi(a^* a) = 1$.
Moreover, there exists a non-zero element $b \in {\cal D}(\sde)$. The strict 
positivity
of $\sde$ implies that $b^* (\sde b) \neq 0$.

\begin{enumerate}
\item
Because $\si_t(\sde) = \gamma^t \sde$, we have that $\si_{-t}(b)$ belongs
to ${\cal D}(\sde)$ and $\si_t(\sde \si_{-t}(b)\,) = \gamma^t \, \sde b$.
The previous lemma
implies that $$(1 \ot \si_{-t}(b)^*)\de(a^* a)(1 \ot \si_{-t}(b))$$ belongs to
${\cal M}_{\vfi \ot \io}$
and
\begin{eqnarray*}
& &(\vfi \ot \io)(\,(1 \ot \si_{-t}(b)^*)\de(a^* a)(1 \ot \si_{-t}(b))\,) \\
& & \spat = \vfi(a^* a) \, \si_{-t}(b)^* (\sde \si_{-t}(b)) = \si_{-t}(b)^* 
(\sde
\si_{-t}(b)).
\end{eqnarray*}
The fact that $\vfi \tau_t = \nu^t \vfi$ implies
that
$$(\tau_t \ot \si_t)(\,
(1 \ot \si_{-t}(b)^*)\de(a^* a)(1 \ot \si_{-t}(b))\,) \text{\ \ \ (a)}$$ 
belongs
to ${\cal M}_{\vfi \ot \io}$
and
\begin{eqnarray*}
& & (\vfi \ot \io)\bigl(\, (\tau_t \ot \si_t)(\,
(1 \ot \si_{-t}(b)^*)\de(a^* a)(1 \ot \si_{-t}(b))\,)\,\bigr) \\
& & \spat = \nu^t \, \si_t\bigl(\,(\vfi \ot \io)(\,
(1 \ot \si_{-t}(b)^*)\de(a^* a)(1 \ot \si_{-t}(b))\,)\,\bigr) \\
& & \spat = \nu^t \, \si_t(\si_{-t}(b)^* (\sde \si_{-t}(b))\,)
= \nu^t \, b^* \, \si_t(\sde \si_{-t}(b)) = \nu^t \gamma^t \, b^* (\sde b).
\end{eqnarray*}
Using the fact that $(\tau_t \ot \si_t)\de = \de \si_t$, we see
that the element (a) equals $(1 \ot b^*)\de(\si_t(a)^* \si_t(a))(1 \ot b)$.

Consequently, we get that
$(1 \ot b^*)\de(\si_t(a)^* \si_t(a))(1 \ot b)$ belongs to
${\cal M}_{\vfi \ot \io}$
and
$$ (\vfi \ot \io)(\, (1 \ot b^*)\de(\si_t(a)^* \si_t(a))(1 \ot b)\,)
= \nu^t \gamma^t \, b^* (\sde b) .$$

\item
We have also that $\si_t(a)$ belongs to $\Nfi$, so the previous lemma gives us
that $(1 \ot b^*)\de(\si_t(a)^* \si_t(a))$ $(1 \ot b)$ belongs to
${\cal M}_{\vfi \ot \io}$
and
$$ (\vfi \ot \io)(\, (1 \ot b^*)\de(\si_t(a)^* \si_t(a))(1 \ot b)\,)
= b^* (\sde b) .$$
\end{enumerate}
Combining these two results, we see that $b^* (\sde b) = \nu^t \gamma^t \,
b^* (\sde b)$, so $\nu^t \gamma^t$ must be equal to 1.
Hence $\gamma^t = \nu^{-t}$.
\end{demo}

This result implies the following result on the $^*$-algebra level.

\begin{proposition}  \label{prop7.2}
For every $z \in \C$, we have that $\rho(\sde^z) = \nu^{i z} \, \sde^z$.
\end{proposition}
\begin{demo}
Choose $a \in A$.
Fix $n \in \N$ and define
$$x_n = \frac{n}{\sqrt{\pi}} \int \exp(-n^2 (s + i z)^2) \, \sde^{is} \pi(a)
\, ds.$$
Take $s \in \R$. By the previous proposition, we have for all $t \in \R$ that
$$\si_t(\sde^{is} \pi(a)) = \nu^{-i s t} \sde^{i s} \si_t(\pi(a)) .$$
Because $\pi(a)$ is analytic with respect to $\si$, this implies that
$\sde^{is} \pi(a)$ belongs to ${\cal D}(\si_{-i})$
and
$$\si_{-i}( \sde^{is} \pi(a) ) = \nu^{-s} \sde^{is} \si_{-i}(\pi(a))
= \nu^{-s} \sde^{is} \pi(\rho(a)).$$
This implies that the function
$$\R \rightarrow  A_r : s \mapsto \exp(-n^2 (s + i z)^2) \,
\si_{-i}(\sde^{is}\pi(a)) $$ is integrable.
Therefore, the closedness of $\si_{-i}$ implies that $x_n$ belongs to
${\cal D}(\si_{-i})$ and
\begin{eqnarray*}
\si_{-i}(x_n) & = &
\frac{n}{\sqrt{\pi}} \int \exp(-n^2 (s + i z)^2) \,
\si_{-i}(\sde^{is}\pi(a)) \, ds \\
& = &
\frac{n}{\sqrt{\pi}} \int \exp(-n^2 (s + i z)^2) \, \nu^{-s} \sde^{is}
\pi(\rho(a)) \, ds.
\end{eqnarray*}
Using the fact that all elements of $\pi(A)$ are analytic with respect to
$\sde$, we get the following expressions for $x_n$ and $\si_{-i}(x_n)$.
\begin{eqnarray*}
x_n & = &  \frac{n}{\sqrt{\pi}} \int \exp(-n^2 s^2) \,
\sde^{is}(\sde^z \pi(a)) \, ds \text{\ \ \ \ \ \ \ \ \ \ \  and } \\
\si_{-i}(x_n) & = &
\frac{n}{\sqrt{\pi}} \int \exp(-n^2 s^2) \, \nu^{-s+iz} \sde^{is}
(\sde^z \pi(\rho(a))\,) \, ds.
\end{eqnarray*}
These two equalities imply that
$(x_n)_{n=1}^\infty \rightarrow \sde^z \pi(a)$
and $(\si_{-i}(x_n)\,)_{n=1}^\infty \rightarrow  \nu^{iz}
\sde^z \pi(\rho(a))$.
Because $\si_{-i}$ is closed, this implies that  $\sde^z \pi(a)$ belongs
to $D(\si_{-i})$ and
$$\si_{-i}(\sde^z \pi(a)) = \nu^{iz} \sde^z \pi(\rho(a)).$$
We can restate this equality in the following form
$$\pi(\rho(\sde^z a)) = \nu^{iz} \pi(\sde^z \rho(a)),$$
which implies that $\rho(\sde^z a) = \nu^{iz} \sde^z \rho(a)$.

From this all, we conclude that $\rho(\sde^z) = \nu^{iz} \sde^z$.
\end{demo}

Remember from section 1 that $\mu$ is by definition the unique complex number
such that $\vfi S^2 = \mu \vfi$. Now, we can make a connection
with our constant $\nu$.

\begin{corollary}
We have that $\mu = \nu^{-i}$.
\end{corollary}
\begin{demo}
From section 1, we know that $\rho(\sde) = \frac{1}{\mu} \sde$,
whereas from the previous proposition, we infer that $\rho(\sde) =
\nu^{i} \sde$. Comparing these two equations, we get our equality.
\end{demo}

\section{The right Haar weight.}

In this section, we are going to introduce the right Haar weight $\psi$. 
First, we will 
define it using the left Haar weight and the anti-unitary antipode. In a 
later stage
we will show that $\psi$ is in some way absolutely continuous with respect to 
$\vfi$ and that the Radon Nikodym derivative is equal to $\sde$. It will also 
be shown that there exists
a non-zero positive right invariant functional on the $^*$-algebra level.

\medskip

Let us start of with the definition of our right Haar weight $\psi$.

\begin{definition}
We define the weight $\psi = \vfi R$ on $A_r$, then $\psi$ is a faithful 
densely defined
lower semi-continuous weight on $A_r$.
\end{definition}

We can also define the modular group of $\psi$:

\begin{definition}
We define the norm-continuous one-parameter group $\si'$ on $A_r$ such
that $\si_t' = R \si_{-t} R$ for every $t \in \R$. Then $\psi$ is a KMS-weight
with modular group $\si'$
\end{definition}

The formula $\de R = \flip(R \ot R)\de$ implies immediately the right 
invariance of $\psi$:

\begin{theorem}
Consider $x \in \Mps$, then $\de(x)$ belongs to $\overline{{\cal M}}_{\psi \ot
\io}$ and  $(\psi \ot \io)\de(x) = \psi(x) \, 1$.
\end{theorem}

\begin{corollary}
Consider $x \in \Mps$ and $\om \in A_r^*$. Then $(\io \ot \om)\de(x)$
belongs to $\Mps$ and $\psi(\,(\io \ot \om)\de(x)\,) = \psi(x) \, \om(1)$.
\end{corollary}

We can also prove easily the following result:

\begin{proposition} \label{prop8.1}
For every $t \in \R$, we have that $\de \si_t' = (\si_t' \ot \tau_{-t})\de$.
\end{proposition}
\begin{demo}
We have that
\begin{eqnarray*}
& & \de \si_t' = \de R \si_{-t} R = \flip(R \ot R)\de \si_{-t} R
= \flip(R \tau_{-t} \ot R \si_{-t})\de R \\
& & \spat = (R \si_{-t} \ot R \tau_{-t})\flip \de R = (R \si_{-t} R \ot R 
\tau_{-t}
R)\de = (\si_t' \ot \tau_{-t}) \de.
\end{eqnarray*}
\end{demo}

We want to use this right Haar weight $\psi$ to define a non-zero positive 
right invariant linear functional on our $^*$-algebra $A$ by the formula 
$\psi \, \pi$. We need for this that
$\pi(A)$ is a subset of $\Mps$ , which is the content of the following lemma.

\begin{lemma}
We have that $\pi(A)$ is a subset of $\Mps$.
\end{lemma}
\begin{demo}
Because $\Nfi$ is a dense left ideal in $A_r$, $R(\Nfi)^*$ is a dense left
ideal in $A_r$.
Choose $x \in \Nfi$ and $\om \in A_r^*$. Using the fact that $(R \ot R)\de
= \flip \de R$, we see that
\begin{eqnarray*}
& & (\io \ot \om)\de(R(x)^*) = (\,(\io \ot \overline{\om})\de(R(x))\,)^* \\
& & \spat =(\,(\io \ot \overline{\om})(\flip(R \ot R)\de(x))\,)^*
= R(\,(\overline{\om} R \ot \io)\de(x)\,)^* .
\end{eqnarray*}
By the left invariance of $\vfi$, we have that $(\overline{\om} R \ot
\io)\de(x)$ belongs to $\Nfi$,  so $(\io \ot \om)\de(R(x)^*)$
belongs to $R(\Nfi)^*$.

By proposition \ref{prop5.2}, we see that $\pi(A) \subseteq R(\Nfi)^*$,
consequently $R(\pi(A)) \subseteq \Nfi^*$.
This implies that
$$R(\pi(A)) = R(\pi(A^* A)) = R(\pi(A)) \, R(\pi(A))^*
\subseteq \Nfi^* \Nfi = \Mfi.$$
Therefore, $\pi(A)$ is a subset of $\Mps$.
\end{demo}

We will now define a faithful positive right invariant linear functional on 
the $^*$-algebra
$A$. For this linear functional on $A$, we will use the same symbol $\psi$ as 
for our right
Haar weight on our \cst\ $A_r$. The context will make clear which viewpoint 
is under consideration. For instance, it is clear that in the following 
definition the weight $\psi$ only turns up in the expression \lq 
$\psi(\pi(a))$ \rq.

\begin{definition}
We define the linear functional $\psi$ on $A$ such that $\psi(a) = 
\psi(\pi(a))$ for every
$a \in A$, then $\psi$ is a faithful positive linear functional on $A$.
\end{definition}

We could prove the right invariance of the linear functional $\psi$ by using 
the right invariance of the weight $\psi$. However, we will prove first 
another interesting formula
which connects $\psi$ with $\vfi S$ and this formula implies also immediately 
the 
right invariance of the linear functional $\psi$.

\begin{proposition}
We have that $\psi(a) = \nu^\frac{i}{2} \, \vfi(S(a))$ for every $a \in A$.
\end{proposition}
\begin{demo}
Choose $a \in A$. By the polar decomposition of the antipode,
we have that $R(\pi(a)) = \tau_\frac{i}{2}(\pi(S(a))\,)$ \ \ (a).
We know already that $\pi(S(a))$ belongs to $\Mfi$. By the previous lemma and 
equality (a),
we see that $\tau_\frac{i}{2}(\pi(S(a))\,)$ also belongs to $\Mfi$.
Because $\vfi \tau_t = \nu^t \vfi$ for every $t \in \R$, these two facts
imply that
$$\vfi(\tau_\frac{i}{2}(\pi(S(a))\,)\,) = \nu^\frac{i}{2} \,
\vfi(\pi(S(a))\,) .$$
Therefore,
$$\psi(a) = \psi(\pi(a)) = \vfi(R(\pi(a))\,) 
= \vfi(\tau_\frac{i}{2}(\pi(S(a))\,)\,) = \nu^\frac{i}{2} \,
\vfi(\pi(S(a))\,) = \nu^\frac{i}{2} \, \vfi(S(a)).$$
\end{demo}

\begin{theorem}
We have that $\psi$ is a faithful positive right invariant functional on $A$.
\end{theorem}

So we have proved a result in the algebraic case, using our \cst ic framework.

\medskip

Another nice formula for $\psi$ on the algebraic level is given in the 
following proposition
(Remember that we can also define $\sde^\frac{1}{2}$ on the algebraic level).

\begin{proposition} \label{prop8.2}
We have that $\psi(a) = \vfi(\sde^\frac{1}{2} \, a \, \sde^\frac{1}{2})$ for 
all $a \in A$.
\end{proposition}
\begin{demo}
We already know that $\psi(a) = \nu^\frac{i}{2} \, \vfi(S(a))$. So, by 
section 1, we get
that $\psi(a) = \nu^\frac{i}{2} \, \vfi(a \, \sde)$. 

By proposition \ref{prop7.2}, we have that $\rho(\sde^\frac{1}{2}) = 
\nu^\frac{i}{2} \,
\sde^\frac{1}{2}$. Hence,
$$\psi(a) = \nu^\frac{i}{2} \, \vfi(a \, \sde^\frac{1}{2} \, \sde^\frac{1}{2})
= \vfi(a \, \sde^\frac{1}{2} \, \rho(\sde^\frac{1}{2}) \,) 
= \vfi(\sde^\frac{1}{2} \, a \, \sde^\frac{1}{2}).$$
\end{demo}

\medskip

Just like in the case of our left Haar weight $\vfi$, we want to show that 
our right Haar weight $\psi$ is completely determined by its values on 
$\pi(A)$. It won't be a great surprise
that both proofs run along the same lines.

Let us take a GNS-triple $(H_\psi,\laps,\pi_\psi)$ for the weight $\psi$.
 
Define $\Gamma_0$ as the closure of the mapping $\pi(A) \rightarrow H_\psi : 
x \mapsto \laps(x)$ and denote the domain of $\Gamma_0$ by $B_0$. It is not 
difficult to see that
$B_0$ is a dense left ideal in $A_r$ which is a subset of $\Nps$. 
Furthermore, $\Gamma_0$ is a closed linear map from $B_0$ into $H_\psi$ which 
is a restriction of $\laps$.

\medskip

The following lemma is the equivalent of lemma \ref{equi1}, but we can 
immediately use the 
right invariance of $\psi$ without proving the preceding lemma's (as in 
section 6).

\begin{lemma}  
Consider $x \in \Nps$ and $\om \in A_r^*$. Then
$(\io \ot \om)\de(x)$ belongs to $\Nps$ and $\|\laps(\,(\io \ot 
\om)\de(x)\,)\|
\leq \|\om\|\,\|\laps(x)\|$.
\end{lemma}
\begin{demo}
We know that there exist an element $\th \in (A_r)_+^*$ with
$\|\th\|=\|\om\|$ such that $$ [(\io \ot \om)\de(y)]^* [(\io \ot \om)\de(y)]
\leq \|\om\| \, (\io \ot \th)\de(y^* y) $$
for all $y \in A_r$. Hence, using the right invariance of $\psi$, we see that
\begin{eqnarray*}
& & \psi( [(\io \ot \om)\de(x)]^* [(\io \ot \om)\de(x)] \,)
\leq \|\om\| \, \psi(\,(\io \ot \th)\de(x^* x)\,)    \\
& & \spat \leq \|\om\| \, \|\th\| \, \psi(x^* x) = \|\om\|^2 \, \psi(x^* x) .
\end{eqnarray*}
The lemma follows.
\end{demo}

The proof of the following lemma is completely analogous to the proof of 
lemma \ref{prop5.3}
and will therefore be left out.

\begin{lemma}  
Consider $\om \in A_r^*$ and $x \in B_0$. Then $(\io \ot \om)\de(x)$
belongs to $B_0$.
\end{lemma}

Similar to lemma \ref{lem5.1}, we have the following lemma.

\begin{lemma}  
Consider $t \in \R$, then $\si_t'(B_0)$ is a subset of $B_0$.
\end{lemma}

The proof is again completely similar to the proof of lemma \ref{lem5.1} with 
the following
slight modifications.
\begin{enumerate}
\item We use proposition \ref{prop5.1} in stead of \ref{prop5.2}.
\item Proposition \ref{prop8.1} is in this situation more useful than 
proposition \ref{prop4.1}.
\end{enumerate}

At last, we get our desired result:

\begin{theorem}
We have that $\pi(A)$ is a core for $\laps$.
\end{theorem}

The proof is analogous to the proof of theorem \ref{thm5.1}, but we have to 
use
the modular group $\si'$ in stead of $\si$.

\bigskip

In the next part, we want to prove that $\psi$ is absolutely continuous with 
respect
to $\vfi$ and that the Radon Nikodym derivative is equal to the modular 
function $\sde$.

\medskip

We know that $\si_t(\sde) = \nu^{-t} \, \sde$ for every $t \in \R$.
Like in the paper of Takesaki \& Pedersen (\cite{Pe-Tak}), it is possible to 
define a
weight $\vfi(\sde^\frac{1}{2}\, . \, \sde^\frac{1}{2})$. However, we cannot
apply their theory because $\sde$ is not invariant under $\si$.
Our construction procedure of this weight
$\vfi(\sde^\frac{1}{2} \, . \, \sde^\frac{1}{2})$ is different from
the one used by Takesaki \& Pedersen. For further details, we refer to 
\cite{JK1}.
Now, we give the basic properties of this weight.

\medskip

We denote $\Upsilon = \vfi(\sde^\frac{1}{2} \, . \, \sde^\frac{1}{2})$, then
$\Upsilon$ is a faithful densely defined lower semicontinuous weight on $A_r$.
Next, we will describe a GNS-construction for $\Upsilon$. In fact, in 
\cite{JK1}, 
the weight $\Upsilon$ is defined via this GNS-construction.

As the GNS-space of $\Upsilon$, we take our Hilbert space $H$.
We define a closed linear map $\la_\Upsilon$ from within $A_r$ into $H$ in the
following way. Put
$$ C_\Upsilon = \{ a \in A_r | \, a \sde^\frac{1}{2} \text{ is bounded and }
\overline{a \sde^\frac{1}{2}} \text{ belongs to } \Nfi \}. $$
Then,  $C_\Upsilon$ is a core for $\la_\Upsilon$ and $\la_\Upsilon(a)
= \lafi(\overline{a \sde^\frac{1}{2}})$ for every $a \in C$.
In this setting, the identity mapping is the GNS-representation of $\Upsilon$

We also define the norm-continuous one-parameter group $\Sigma$ on $A_r$
such that $\Sigma_t(a) = \sde^{it} \si_t(a) \sde^{-it}$ for every $a \in A_r$
and $t \in \R$. We have that $\Upsilon$ is KMS with respect to $\Sigma$.

\medskip

If we denote the modular operator of $\Upsilon$ by $\nabp$, we have that 
$\nabp^{it} = J \sde^{it} J \sde^{it} \nab^{it}$ for every $t \in \R$.

If we denote the modular conjugation of $\Upsilon$ by $J'$, we have that
$J' = \lambda^\frac{i}{4} J$ .

\bigskip

Now, we are going to prove some results about $\Upsilon$ in this specific
case.

\begin{lemma}     \label{lem8.1}
Let $a \in A$. Then $\pi(a)$ belongs to $C_\Upsilon \subseteq {\cal 
N}_\Upsilon$
and $\la_\Upsilon(\pi(a)) = \la(a \sde^\frac{1}{2})$
\end{lemma}
\begin{demo}
We know that $\pi(a^*)$ belongs to ${\cal D}(\sde^\frac{1}{2})$ and
$\sde^\frac{1}{2} \pi(a^*) = \pi(\sde^\frac{1}{2} a^*)$.
This implies that $\pi(a) \sde^\frac{1}{2}$  is bounded and 
$$\overline{\pi(a) \sde^\frac{1}{2}} = (\sde^\frac{1}{2} \pi(a^*))^* =
\pi(a \sde^\frac{1}{2}) \ , $$ which belongs
to $\Nfi$. By definition, we have that $\pi(a)$ belongs to $C_\Upsilon$
and $$\la_\Upsilon(\pi(a)) = \lafi(\pi(a \sde^\frac{1}{2})) =
\la(a \sde^\frac{1}{2}). $$
\end{demo}

\begin{lemma} 
Consider $x \in \pi(A)$, then $x$ belongs to ${\cal M}_\Upsilon$ and
$\Upsilon(x) = \psi(x)$.
\end{lemma}
\begin{demo}
Choose $b,c \in A$. We know from the previous lemma that $\pi(b),\pi(c)$
belong to ${\cal N}_\Upsilon$. This implies that $\pi(c^* b)$ belongs
to $\Mps$ and
$$\Upsilon(\pi(c^* b)) =    \langle \la_\Upsilon(\pi(b)),\la_\Upsilon(\pi(c)) 
\rangle
=  \langle \la(b \sde^\frac{1}{2}),\la(c \sde^\frac{1}{2}) \rangle
= \vfi(\sde^\frac{1}{2} c^* b \sde^\frac{1}{2}).$$
From proposition \ref{prop8.2}, we know that this last term equals 
$\psi(\pi(c^*b))$.

The lemma follows because $A^* A = A$.
\end{demo}

The next lemmma is the last step towards a real equality between $\Upsilon$ 
and $\psi$.

\begin{lemma}
We have that $\pi(A)$ is a core for $\la_\Upsilon$.
\end{lemma}
\begin{demo}
Define $\ga$ as the closure of the map $\pi(A) \rightarrow H: x \mapsto
\la_\Upsilon(x)$ and call $B$ the domain of $\ga$. It is clear that 
$\la_\Upsilon$
is an extension of $\ga$.

Choose $x \in \Nfi \cap {\cal N}_\Upsilon$.

Take $n \in \N$  and put
$$ x_n = \frac{n}{\sqrt{\pi}} \int \exp(-n^2 t^2) \, x
\sde^{it} \, dt .$$
Because $\Sigma_s(\sde) = \nu^{-s} \, \sde$ for all $s \in
\R$, we get that $x_n$  belongs to ${\cal N}_\Upsilon$ and
$$\la_\Upsilon(x_n) =  \frac{n}{\sqrt{\pi}} \int \exp(-n^2 t^2) \, 
\nu^\frac{t}{2}
J' \sde^{-it} J' \la_\Upsilon(x) \, dt   \text{\ \ \ (a)} $$
(cf. lemma \ref{lem5.2}).

Furthermore, we have that $x_n \sde^\frac{1}{2}$ is bounded and
$$\overline{x_n \sde^\frac{1}{2}}
= \frac{n}{\sqrt{\pi}} \int \exp(-n^2 (t+\frac{i}{2})^2) \, x \sde^{it} \, dt
\ .$$
Because $\si_s(\sde) = \nu^{-s} \, \sde$ for all $s \in \R$, this implies
that $\overline{x_n \sde^\frac{1}{2}}$ belongs to $\Nfi$
and
$$ \lafi(\overline{x_n \sde^\frac{1}{2}})
= \frac{n}{\sqrt{\pi}} \int \exp(-n^2 (t+\frac{i}{2})^2) \, \nu^\frac{t}{2}
J \sde^{-it} J \lafi(x) \, dt \ .$$

By definition, we get that
$$\la_\Upsilon(x_n)
= \frac{n}{\sqrt{\pi}} \int \exp(-n^2 (t+\frac{i}{2})^2) \, \nu^\frac{t}{2}
J \sde^{-it} J \lafi(x) \, dt \ .  \text{\ \ \ (b)} $$

\medskip

Choose $\vep \in \R_0^+$ and $y \in \Nfi \cap {\cal N}_\Upsilon$.
From (a), we infer the existence of a natural number $m$ such that
$\|y - y_m \| \leq \frac{\vep}{2}$
and $\|\la_\Upsilon(y) - \la_\Upsilon(y_m)\| \leq \frac{\vep}{2}$.
Because $\pi(A)$ is a core for $\lafi$, (b) implies the existence of an
element $z \in \pi(A)$ such that
$\|z_m - y_m \| \leq \frac{\vep}{2}$
and $\|\la_\Upsilon(z_m) - \la_\Upsilon(y_m)\| \leq \frac{\vep}{2}$.
This implies that
$\|y - z_m \| \leq \vep$ and $\|\la_\Upsilon(y) - \la_\Upsilon(z_m)\| \leq 
\vep$.
Next, we show that $z_m$ belongs to $B$.

\begin{list}{}{\setlength{\leftmargin}{.4 cm}}

\item Remembering proposition \ref{prop6.1}, we get for any $t \in \R$ that
$z \sde^{it}$ belongs to $\pi(A)$
and
$$\ga(z \sde^{it}) = \la_\Upsilon(z \sde^{it})
= \nu^\frac{t}{2} J' \sde^{-it} J' \la_\Upsilon(z). $$
Consequently, the function
$\R \rightarrow H : t \mapsto \exp(-m^2 t^2) \,  \ga(z \sde^{it})$ is
integrable. The closedness of $\ga$ implies that $z_m$ belongs to $B$.

\end{list}

It is not so difficult to prove, in general, that $\Nfi \cap {\cal 
N}_\Upsilon$ is a core
for $\la_\Upsilon$ (see \cite{JK1}).  This implies with the foregoing results 
that $\la_\Upsilon=
\ga$. So, $\pi(A)$ is a core for $\la_\Upsilon$

\end{demo}

So, we arrive at the following nice conclusion:

\begin{theorem}
We have that $\psi$ equals $\Upsilon$
\end{theorem}

The unicity of the modular group implies that $\si'$ must be equal to 
$\Sigma$. So, we get
the following nice proposition.

\begin{proposition}
For every $t \in \R$ and $a \in A_r$, we have that $\si_t'(a) = \sde^{it} 
\si_t(a) \sde^{-it}$.
\end{proposition}

This formula for $\si'$ allows us to prove a continuous version of lemma 
\ref{lem2.1}.

\begin{proposition}
For every $t \in \R$, we have that $\de \tau_t = (\si_t \ot \si_t') \de$.
\end{proposition}
\begin{demo}
Choose $x \in A_r$. Then
\begin{eqnarray*}
& & (\si_t \ot \si_{-t}')\de(x)
= (\sde^{-it} \ot \sde^{-it})(\si_t' \ot \si_{-t})(\de(x))(\sde^{it} \ot
\sde^{it}) \\
& & \spat = (\sde^{-it} \ot \sde^{-it})(\io \ot \si_{-t} 
\tau_t)(\de(\si_t'(x))\,)
(\sde^{it} \ot \sde^{it})  \\
& & \spat = (\sde^{-it} \ot \sde^{-it})(\io \ot \tau_t \si_{-
t})(\de(\si_t'(x))\,)
(\sde^{it} \ot \sde^{it})  \\
& & \spat = (\sde^{-it} \ot \sde^{-it})(\tau_t \ot \tau_t)
(\de(\si_{-t}(\si_t'(x))\,)\,) (\sde^{it} \ot \sde^{it})  \\
& & \spat = (\tau_t \ot \tau_t)(\,(\sde^{-it} \ot \sde^{-it})
\de(\sde^{it} x \sde^{-it})(\sde^{it} \ot \sde^{it})\,)  \\
& & \spat = (\tau_t \ot \tau_t)(\de(\sde^{-it})\de(\sde^{it} x \sde^{-it})
\de(\sde^{it})\,) \\
& & \spat = (\tau_t \ot \tau_t)\de(x) = \de(\tau_t(x)).
\end{eqnarray*}
\end{demo} 

\begin{corollary}
We have that $K_t(a) = \sde^{-it} \tau_{-t}(a)  \sde^{it}$
for every $t \in  \R$ and $a \in A_r$.
\end{corollary}
\begin{demo}
Choose $x \in A_r$ and $\om \in A_r^*$.
By the previous proposition, we have that
$\de(\si_t'(x)) = (\si_t' \ot \tau_{-t})\de(x)$.
Remembering that $\si_t'(y) = \sde^{it} \si_t(y) \sde^{-it}$ for all $y \in
A_r$, we infer from this that
$$ (\sde^{it} \ot \sde^{it}) \de(\si_t(x)) (\sde^{-it} \ot \sde^{-it})
=  (\sde^{it} \ot 1) (\si_t \ot \tau_{-t})(\de(x)) (\sde^{-it} \ot 1) \ ,$$
so,
$\de(\si_t(x))
=(1 \ot \sde^{-it})(\si_t \ot \tau_{-t})(\de(x))(1 \ot \sde^{it})$.
From proposition \ref{prop2.1}, we know already that
$\de(\si_t(x)) = (\si_t \ot K_t)\de(x)$, therefore,
$$(1 \ot \sde^{-it})(\io \ot \tau_{-t})(\de(x))(1 \ot \sde^{it})
= (\io \ot K_t)\de(x).$$
Applying $\om \ot \io$ on this equation, results in the equality
$\sde^{-it} \tau_{-t}\bigl((\om \ot \io)\de(x)\bigr) \sde^{it}
= K_t\bigl((\om \ot \io)\de(x)\bigr)$. The proposition follows now from the
density conditions.
\end{demo}

\medskip

The existence of a non-zero positive right invariant linear functional $\psi$ 
on 
the algebraic quantum group $(A,\de)$ allows us to find other implementations 
of 
the polar decomposition of the antipode. Also, we can prove the manageability 
of $W$
(see \cite{Wor3}). For this, we use the construction procedure of section 4 
with $\eta$ equal
to $\psi$.

For every $a \in A$, we have that $\psi(a) = \nu^\frac{i}{2} \, \vfi(S(a))$. 
Therefore, we have
in this case that $x = \nu^\frac{i}{2} 1$. We take $y = 1$.
Because $\psi(a) = \vfi(\sde^\frac{1}{2} a \sde^\frac{1}{2})$ for all $a \in 
A$, we can choose
our GNS-pair $(K,\ga)$ such that $K=H$ and $\ga(a) = \la(a \sde^\frac{1}{2})$ 
for all $a \in A$. We put $Q=\cp$ and $D=\cj$. Again, we give a summary of 
some of the results of section 4 in this case.

\medskip

We have that $G$ is the closed antilinear map from within $H$ into $H$
such that $\la(A)$ is a core for $G$ and $G \la(a) = \ga(S(a)^*)= \la(S(a)^* 
\sde^\frac{1}{2})$
for every $a \in A$. It is clear that $G$ is involutive in this case.

Therefore,  $D$ is an involutive anti-unitary transformation on $H$,  \
$Q=G^* G$ is an injective positive operator in $H$
such that $G = D Q^\frac{1}{2} = Q^{-\frac{1}{2}} D$.
For $t \in \R$ we have that $Q^{it} =  I Q^{it} I$ and
$Q^t = I Q^{-t} I$.

\medskip

For any $a,b \in A$, we have by definition of $V$ and $W$ that
\begin{eqnarray*}
& & V(\la \od \la)(\de(b)(a \sde^\frac{1}{2} \ot 1)\,)
= V(\ga \od \la)(\de(b)(a \ot 1)\,)  = (\ga \ot \la)(a \ot b) \\
& & \spat\spat (\la \od \la)(a \sde^\frac{1}{2} \ot b)
= W(\la \od \la)(\de(b)(a \sde^\frac{1}{2} \ot 1)\,)  \ .
\end{eqnarray*}

So, we get in this case that $W=V$. This implies that
$$ (Q \ot \nab)W = W(Q \ot \nab) \text{\ \ \ \ and \ \ \ \ }
(D \ot J)W = W^* (D \ot J)  . $$
These two commutation relations allow us to prove the next proposition in an 
analogous
way as proposition \ref{prop6.6}.

\begin{proposition}
\begin{enumerate}
\item For every $t \in \R$ and $x \in A_r$, we have that $\tau_t(x) = Q^{it} 
\, x \, Q^{-it}$.
\item For every $x \in A_r$, we have that $R(x)= D \, x^* \, D$.
\end{enumerate}
\end{proposition}

It is clear that in this case $S(y^*) \, x \, \rho(S(y)\,) \in \C 1$. So the 
conclusion
of proposition \ref{prop3.5} holds. We also have the result of proposition 
\ref{prop3.4}:
\begin{itemize}
\item Consider $u_1 \in D(Q^\frac{1}{2})$, $u_2 \in D(Q^{-\frac{1}{2}})$ and
$v_1,v_2 \in H$. Then
$$  \langle W^* (Q^\frac{1}{2} u_1 \ot J v_1) , Q^{-\frac{1}{2}} u_2 \ot J 
v_2 \rangle
=  \langle W (u_1 \ot v_2), u_2 \ot v_1 \rangle .$$
\item We have that $W(Q \ot Q) = (Q \ot Q)W$.
\end{itemize}

By definition 1.2 of \cite{Wor3}, we arrive at the following conclusion.

\begin{theorem}
We have that $W$ is manageable.
\end{theorem}

\section{Appendix: some information about weights}

In this section, we will collect some necessary information and conventions 
about weights.

\medskip

Consider a $C^*$-algebra $A$  and a densely defined lower semi- continuous 
weight $\varphi$
on $A$. We will use the following notations:
\begin{itemize}
\item ${\cal M}^+_\varphi = \{ a \in A^+ \mid \varphi(a) < \infty  \} $
\item ${\cal N}_\varphi = \{ a \in A \mid \varphi(a^*a) < \infty \} $
\item ${\cal M}_\varphi = \text{span\ } {\cal M}^+_\varphi
= {\cal N}_\varphi^* {\cal N}_\varphi$ . 
\end{itemize}

\medskip

A GNS-construction of $\varphi$ is by definition a triple
$(H_\varphi,\pi_\varphi,\Lambda_\varphi)$ such that
\begin{itemize}
\item $H_\varphi$ is a Hilbert space
\item $\Lambda_\varphi$ is a linear map from ${\cal N}_\varphi$ into
      $H_\varphi$ such that
      \begin{enumerate}
      \item  $\Lambda_\varphi({\cal N}_\varphi)$ is dense in $H_\varphi$
      \item   For every $a,b \in {\cal N}_\varphi$, we have that
              $\langle \Lambda_\varphi(a),\Lambda_\varphi(b) \rangle
              = \varphi(b^*a) $
      \end{enumerate}
      Because $\varphi$ is lower semi-continuous, $\Lambda_\varphi$ is closed.
\item $\pi_\varphi$ is a non-degenerate representation of $A$ on
      $H_\varphi$ such that $\pi_\varphi(a)\,\Lambda_\varphi(b) = 
\Lambda_\varphi(ab)$
      for every $a \in M(A)$ and $b \in {\cal N}_\varphi$.  (The non-
degeneracy of             $\pi_\varphi$ is a consequence of the lower semi-
continuity of $\varphi$.)
\end{itemize}

The following concepts play a central role in the theory of lower
semi-continuous weights.

\medskip

We define ${\cal F}_\varphi = \{ \omega \in A^*_+ \mid
\omega \leq \varphi \}$.
Because $\varphi$ is lower semi-continuous, we have that
$$ \varphi(x) = \sup_{\omega \in {\cal F}_\varphi} \omega(x) $$
for all $x \in A^+$.

Also, we put ${\cal G}_\varphi = \{ \alpha\,\omega \mid
\omega \in {\cal F}_\varphi \, , \, \alpha \in  \,\, ]0,1[ \,\,\}$.
Then ${\cal G}_\varphi$ is a directed subset of ${\cal F}_\varphi$
such that
$$ \varphi(x) = \lim_{\omega \in {\cal G}_\varphi} \omega(x)$$
for all $x \in A^+$.

It follows easily that
$$ \varphi(x) = \lim_{\omega \in {\cal G}_\varphi} \omega(x)$$
for every $x \in {\cal M}_\varphi$.

\medskip

Any lower semi-continuous weight $\vfi$ has a natural extension to a weight 
$\overline{\vfi}$ 
on $M(A)^+$ by putting $\overline{\vfi}(x)= \sup\{ \om(x) \mid \om \in {\cal 
F}_\vfi \}$
for every $x \in M(A)^+$.

We define ${\overline{\cal M}}_\vfi = {\cal M}_{\overline{\vfi}}$ and 
${\overline{\cal N}}_\vfi = {\cal N}_{\overline{\vfi}}$. For any $x \in 
{\overline{\cal M}}_\vfi$, we put $\vfi(x)
= \overline{\vfi}(x)$.

Then, we have that
$$ \varphi(x) = \lim_{\omega \in {\cal G}_\varphi} \omega(x)$$
for every $x \in {\overline{\cal M}}_\varphi$.

\bigskip

Now, we will introduce the class of KMS-weights. All weights used in this 
paper, will
belong to this class.

\begin{definition}
Consider a \cst\ $A$ and a weight $\vfi$ on $A$, we say that $\vfi$ is a KMS-
weight on A
if and only if $\vfi$ is a faithful densely defined lower semi-continuous 
weight on A
such that there exists  a norm-continuous one-parameter group $\si$ on 
$A$ satisfying the following poperties:
\begin{enumerate}
\item $\vfi$ is invariant under $\si$: $\vfi \si_t = \vfi$ for every $t \in 
\R$.
\item For every $a \in {\cal D}(\si_\frac{i}{2})$, we have that
      $\vfi(a^* a) = \vfi(\si_\frac{i}{2}(a) \si_\frac{i}{2}(a)^*)$.
\end{enumerate}
\end{definition}

If a weight is KMS in the sense of the previous definition, then the one-
parameter group $\si$
is unique and is called the modular group of $\si$. 

\medskip

It is possible to replace condition 2) with a weaker condition like:

There exist a core $K$ for $\lafi$ such that
\begin{itemize}
\item $\si_t(K) \subseteq K$ for every $t \in \R$.
\item We have that $K \subseteq {\cal D}(\si_\frac{i}{2})$,             
$\si_\frac{i}{2}(K)^* 
      \subseteq \Nfi$ and
      $$ \| \lafi(x) \| = \| \lafi(\si_\frac{i}{2}(x)^*) \| $$ for every 
      $x \in K$.
\end{itemize}

\medskip

We are going to give some properties about KMS-weights.

\begin{lemma}  \label{lem5.2}
Consider a \cst\ $A$, and a KMS-weight $\vfi$ on $A$ with modular group 
$\si$. Then:
\begin{enumerate}
\item There exists a unique anti-unitary operator $J$ on $H_\vfi$ such
      that $J \lafi(x) = \lafi(\si_\frac{i}{2}(x)^*)$ for every $x \in        
                 \Nfi \cap D(\frac{i}{2})$.       
\item Let $a \in {\cal D}(\si_\frac{i}{2})$ and $x \in \Nfi$. Then $x a$      
             belongs to $\Nfi$ and $\lafi(x a) = J 
\pi_\vfi(\si_\frac{i}{2}(a))^* J \, \lafi(x)$.
\item Let $a \in {\cal D}(\si_{-i})$ and $x \in \Mfi$. Then $a x$ and $x 
\si_{-i}(a)$ belong       to $\Mfi$ and $\vfi(a x) = \vfi(x \si_{-i}(a))$.
\end{enumerate}
\end{lemma}

We want to mention that in the definitition of a KMS-weight the condition 2) 
can be replaced by the existence of a anti-unitary $J$ on $H_\vfi$ such that 
condition 2) of the previous lemma
is satisfied (we can even allow this equality to hold for fewer (but enough) 
elements).

\medskip

Define $\cu=\lafi(\Nfi \cap \Nfi^*)$, then $\cu$ is turned into a $^*$-
algebra in the usual way:
\begin{itemize}
\item For every $a,b \in \Nfi \cap \Nfi^*$, we have that $\lafi(a) \lafi(b) = 
\lafi(ab)$.
\item For every $a \in \Nfi \cap \Nfi^*$, we have that $\lafi(a)^* = 
\lafi(a^*)$.
\end{itemize}
It is not so difficult to prove that $\cu$ is a left Hilbert algebra on 
$H_\vfi$.

\medskip

For every element $a \in A$, we have that $\lafi(a)$ is left bounded with 
respect to $\cu$
and $L_{\lafi(a)}= \pi_\vfi(a)$.

Define $T$ as the closed antilinear mapping from within $H_\vfi$ into 
$H_\vfi$ such that
$\cu$ is a core for $T$ and $T v = v^*$ for every $v \in \cu$. Denote the 
modular operator
by $\nab=T^* T^*$. Then $J$,$\nab^\frac{1}{2}$ is the polar decomposition of 
$T$. We call $\nab$ the modular operator of $\vfi$ and 
$J$ the modular conjugation of $\vfi$.

\medskip

We also have that $\nab^{it} \lafi(a) = \lafi(\si_t(a))$ for every $t \in \R$ 
and $a \in \Nfi$. 
As a consequence, $\pi_\vfi(\si_t(a)) = \nab^{it} \pi(a) \nab^{-it}$ for 
every $t \in \R$
and $a \in A$.

As a last remark, we want to say that $\vfi$ satisfies the so-called KMS-
condition with 
respect to $\si$.

\medskip

We will need the following generalization of theorem 3.6 of \cite{Pe-Tak}.

\begin{lemma} \label{lemA1}
Consider a \cst\ $A$ and a KMS-weight $\vfi$ on $A$ with modular group $\si$.
Let $a$ be an element in $M(A)$ such that $a \Mfi \subseteq \Mfi$ and $\Mfi a 
\subseteq \Mfi$
and such that there exists a strictly positive number $\lambda$ such
that $\vfi(a x) = \lambda \vfi(x a)$ for every $x \in \Mfi$. Then we have that
$\si_t(a) = \lambda^{it} a$ for every $t \in \R$.
\end{lemma}

\bigskip

In a last part we say something about slicing with weights.
Therefore, we fix a \cst\ A and a KMS-weight $\vfi$ on $A$.

\medskip

\begin{definition}
We define the map $\io \ot \vfi$ from within $(A \ot A)^+$ into $A^+$ as 
follows:
\begin{itemize}
\item We define the set  ${\cal M}_{\io \ot \vfi}^+ = \{ a \in (A \ot A)^+ 
\mid
      \text{\ \ the net \ \ } \bigl( \,(\io \ot \om)(a)\,\bigr)_{\om \in 
{\cal G}_\vfi}
      \text{\ \ is norm convergent in \  } A \}$.
\item The mapping $\io \ot \vfi$ will have as domain the set ${\cal M}_{\io 
\ot \vfi}^+$
      and for any $a \in {\cal M}_{\io \ot \vfi}^+$, we have by definition 
that the net
      $\bigl( \,(\io \ot \om)(a)\,\bigr)_{\om \in {\cal G}_\vfi}$ converges
      to $(\io \ot \vfi)(a)$.
\end{itemize}
\end{definition}

It is clear that ${\cal M}_{\io \ot \vfi}^+$ is a dense hereditary cone in 
$(A \ot A)^+$.
Furthermore we define the following sets:
\begin{enumerate}
\item We define ${\cal N}_{\io \ot \vfi} = \{ a \in A \ot A \mid a^* a 
\text{\ \ belongs
      to \ \ } {\cal M}_{\io \ot \vfi}^+ \}$. 
\item Also, we define ${\cal M}_{\io \ot \vfi} = \text{span \ } {\cal M}_{\io 
\ot \vfi}^+
      = {\cal M}_{\io \ot \vfi}^* {\cal N}_{\io \ot \vfi}$.
\end{enumerate}
Of course, there exist a unique linear map $\psi$ from ${\cal M}_{\io \ot 
\vfi}$
to $A$ which extend $\io \ot \vfi$. For any $a \in {\cal M}_{\io \ot \vfi}$, 
we put
$(\io \ot \vfi)(a) = \psi(a)$.

It is then clear that $\bigl( \,(\io \ot \om)(a)\,\bigr)_{\om \in {\cal 
G}_\vfi}$ converges
to $(\io \ot \vfi)(a)$ for every $a \in {\cal M}_{\io \ot \vfi}$.

We also have for any $a \in {\cal M}_{\io \ot \vfi}$ and $\th \in A^*$ that
$(\th \ot \io)(a)$ belongs to $\Mfi$ and $\vfi(\,(\th \ot \io)(a)\,) = 
\th(\,(\io \ot \vfi)(a)\,)$.

\medskip

We are now going to describe a GNS-construction for $\io \ot \vfi$.
It is possible to prove that the mapping $\io \od \lafi$
from $A \od \Nfi$ into $A \od H$ is closable (as a mapping from the \cst\ $A 
\ot A$ into
the Hilbert-C$^*$-module $A \ot H$). We define $\io \ot \lafi$ to be the 
closure
of this mapping $\io \od \lafi$.

It is also possible to prove that $D(\io \ot \lafi)$ is a subset of ${\cal 
N}_{\io \ot \vfi}$ and that 
$$(\io \ot \vfi)(b^* a) = \langle (\io \ot \lafi)(a) , (\io \ot \lafi)(b) 
\rangle $$
for every $a,b \in D(\io \ot \lafi)$.

All that is said about $\io \ot \vfi$ untill now goes through for lower semi-
continuos weights  
(and can be found in \cite{Q-V} and \cite{Ver}). 
Because we assumed that $\vfi$ is also KMS, we have also the rather non-
trivial result
that $D(\io \ot \lafi) = {\cal N}_{\io \ot \vfi}$. This last result is also 
true if
$\vfi$ would obey a weaker condition called regularity (see \cite{JK2}).

\bigskip

We would like to have an extension of $\io \ot \vfi$ to $M(A \ot A)$. This is 
done in
the following way:

\begin{definition}
We define the map $\overline{\io \ot \vfi}$ from within $M(A \ot A)^+$ into 
$M(A)^+$ as follows:
\begin{itemize}
\item We define the set  ${\overline{\cal M}}_{\io \ot \vfi}^+ = \{ a \in M(A 
\ot A)^+ \mid
\text{\ \ the net \ \ } \bigl( \,(\io \ot \om)(a)\,\bigr)_{\om \in {\cal 
G}_\vfi}
\text{\ \ is strictly convergent in }$ $M(A) \}$.
\item The mapping $\overline{\io \ot \vfi}$ will have as domain the set 
${\overline{\cal M}}_{\io \ot \vfi}^+$
      and for any $a \in {\overline{\cal M}}_{\io \ot \vfi}^+$, we have by 
definition that the net $\bigl( \,(\io \ot \om)(a)\,\bigr)_{\om \in {\cal 
G}_\vfi}$ converges strictly
      to $(\overline{\io \ot \vfi})(a)$.
\end{itemize}
\end{definition}

This definition is in fact not entirely correct because it depends on $\vfi$ 
and not on 
$\io \ot \vfi$. It is possible to give a definition in terms of the mapping
$\io \ot \vfi$ and our definition would then be a proposition.

\medskip

The next proposition reveals a nice feature about $\overline{\io \ot \vfi}$.

\begin{proposition}
Consider $a \in M(A \ot A)^+$. Then $a$ belongs to ${\overline{\cal M}}_{\io 
\ot \vfi}^+ $
if and only if the net \newline $\bigl( \,b^* (\io \ot \om)(a) b 
\,\bigr)_{\om \in {\cal G}_\vfi}$ 
is norm convergent for every $b \in A$.
\end{proposition}

It is clear that ${\overline{\cal M}}_{\io \ot \vfi}^+$ is a hereditary cone 
in $M(A \ot A)^+$.
Furthermore we define the following sets:
\begin{enumerate}
\item We define ${\overline{\cal N}}_{\io \ot \vfi} = \{ a \in A \ot A \mid 
a^* a \text{\ \                         belongs to \ \ } {\overline{\cal 
M}}_{\io \ot \vfi}^+ \}$. 
\item Also, we define ${\overline{\cal M}}_{\io \ot \vfi} = \text{span \ } 
{\overline{\cal                M}}_{\io \ot \vfi}^+ = {\overline{\cal 
N}}_{\io \ot \vfi}^* {\overline{\cal N}}_{\io \ot                   \vfi}$.
\end{enumerate}
Of course, there exist a unique linear map $\overline{\psi}$ from 
${\overline{\cal M}}_{\io \ot \vfi}$
to $M(A)$ which extend $\overline{\io \ot \vfi}$. For any $a \in 
{\overline{\cal M}}_{\io \ot \vfi}$, we put
$(\io \ot \vfi)(a) = \overline{\psi}(a)$.

It is then clear that $\bigl( \,(\io \ot \om)(a)\,\bigr)_{\om \in {\cal 
G}_\vfi}$ converges
strictly to $(\io \ot \vfi)(a)$ for every $a \in {\overline{\cal M}}_{\io \ot 
\vfi}$.

We also have for any $a \in {\overline{\cal M}}_{\io \ot \vfi}$ and $\th \in 
A^*$ that
$(\th \ot \io)(a)$ belongs to ${\overline{\cal M}}_\vfi$ and $\vfi(\,(\th \ot 
\io)(a)\,) = \th(\,(\io \ot \vfi)(a)\,)$.  

\medskip

It is clear that we can do the same things for the slice-mapping $\vfi \ot 
\io$.


\begin{thebibliography}{VD}

\bibitem{Abe} {\sc E. Abe},
Hopf Algebras. {\it Cambridge University Press} (1977).


\bibitem{B-S}  {\sc S. Baaj  \&  G. Skandalis},
Unitaires multiplicatifs et dualit\'e pour les produits crois\'es de
C$^*$-alg\`ebres. {\it Ann. scient. \'{E}c. Norm. Sup., $4{}^e$ s\'{e}rie, t. 
26}
(1993), 425--488.

\bibitem{ER} {\sc  E.G. Effros \&  Z.-J. Ruan},
Discrete Quantum Groups I. The Haar Measure. {\it Int. J. of Math.} (1994).

\bibitem{E}  {\sc  M. Enock \&  J.-M. Schwartz},
Kac Algebras and Duality of Locally Compact Groups.
{\it Springer-Verlag, Berlin}  (1992).

\bibitem{Drab}{\sc B. Drabant \& A. Van Daele}, Pairing and Quantum Double of 
Multiplier Hopf
Algebras. {\it Preprint K.U.Leuven}  (1996)


\bibitem{GL} {\sc  E.C. Gootman \& A.J. Lazar},
Quantum Groups and Duality.
{\it Reviews in Math. Physics} {\bf 5} No. 2 (1993),
417--451

\bibitem{JK1} {\sc J. Kustermans}, A construction procedure for KMS-weights 
on C$^*$-algebras.
In preparation.

\bibitem{JK2} {\sc J. Kustermans}, Regular C$^*$-valued weights on C$^*$-
algebras. In preparation.

\bibitem{Lan} {\sc C. Lance},
Hilbert $C^*$-modules, a toolkit for operator algebraists. Leeds. (1993).

\bibitem{MasNak} {\sc T. Masuda \& Y. Nakagami} A von Neumann Algebra 
Framework for the Duality
of Quantum Groups. {\it Publications of the RIMS Kyoto University} {\bf 30} 
(1994) , 799--850

\bibitem{Pe-Tak} {\sc G.K. Pedersen \& M. Takesaki},
The Radon-Nikodym theorem for von Neumann algebras.
{\it Acta Math.} {\bf 130} (1973), 53--87.


\bibitem{PW}  {\sc  P. Podle\'s \&  S.L. Woronowicz}, Quantum Deformation of 
the Lorentz Group.
{\it Commun. Math. Phys.} {\bf 130} (1990), 381--431.


\bibitem{Q-V} {\sc J. Quaegebeur \& J. Verding},
A construction for weights on $C^*$-algebras. Dual weights for $C^*$-crossed
products. {\it Preprint K.U. Leuven} (1994).


\bibitem{Ver} {\sc J. Quaegebeur \& J. Verding},
Left invariant weights and the left regular corepresentation
for locally compact quantum semi-groups. {\it Preprint K.U. Leuven} (1994).


\bibitem{Stra} {\sc S. Stratila \& L. Zsido},
Lectures on von Neumann algebras. {\it Abacus Press, Tunbridge Wells, 
England} (1979).

\bibitem{Tay}  {\sc D.C. Taylor}, The Strict Topology for Double Centralizer 
Algebras.
{\it Trans. Am. Math. Soc}{\bf 150} (1970), 633 -- 643

\bibitem{Tak}  {\sc  M. Takesaki},
Theory of Operator Algebras I.
{\it Springer-Verlag, New York} (1979).



\bibitem{VD1}   {\sc A. Van Daele}, An Algebraic Framework for Group Duality.
{\it Preprint K.U.Leuven} (1996).

\bibitem{VD2}   {\sc  A. Van Daele}, Dual Pairs of Hopf ${}^*$-algebras.
{\it Bull. London Math. Soc.} {\bf 25} (1993), 209--230.


\bibitem{VD4}   {\sc  A. Van Daele}, Discrete Quantum Groups.
{\it Journal of Algebra}{\bf 180} (1996), 431--444.

\bibitem{VD5}   {\sc   A. Van Daele}, The Haar Measure on a Compact
Quantum Group. {\it Proc. Amer. Math. Soc.}{\bf 123}(1995), 3125-3128

\bibitem{VD6}   {\sc   A. Van Daele}, Multiplier Hopf Algebras. {\it Trans. 
Am. Math. Soc.}{\bf 342} (1994), 917--932.

\bibitem{Verd} {\sc J. Verding}, Weights on C$^*$-algebras. {\it Phd-thesis.  
K.U. Leuven} (1995)

\bibitem{Wor1}   {\sc  S.L. Woronowicz},  Compact matrix pseudogroups.
{\it Commun. Math. Phys.} {\bf 111}  (1987),  613--665.

\bibitem{Wor2}   {\sc  S.L. Woronowicz},  Compact quantum groups.
{\it Preprint Warszawa} (1993).

\bibitem{Wor3}  {\sc  S.L. Woronowicz}, From multiplicative unitaries to 
quantum
groups. {\it Preprint Warszawa} (1995).

\bibitem{Wor4}  {\sc  S.L. Woronowicz}, Pseudospaces, pseudogroups and 
Pontriagin duality.
{\it Proceedings of the International Conference on
Mathematical Physics, Lausanne} (1979),  407--412.


\bibitem{Wor6}  {\sc  S.L. Woronowicz}, Unbounded elements affiliated with 
$C^*$-algebras
and non-compact quantum groups.
{\it Commun. Math. Phys.} {\bf  136} (1991),  399--432.


\end{thebibliography}
\end{document}